\documentclass{amsart} 
\usepackage[margin=1in]{geometry}

\usepackage[mathletters]{ucs}
\usepackage[utf8x]{inputenc}

\usepackage{graphicx}
\usepackage{amssymb, amsmath, amsthm, graphics}
\usepackage{mathtools}
\usepackage{adjustbox}

\usepackage[mathscr]{euscript}
\DeclareMathAlphabet{\mathpzc}{OT1}{pzc}{m}{it}

\usepackage{tikz}
\usepackage{tikz-cd}
\usetikzlibrary{calc}
\usetikzlibrary{matrix}
\usetikzlibrary{positioning}
\usetikzlibrary{shapes.geometric}
\usetikzlibrary{shapes.misc}
\usetikzlibrary{shapes,arrows}
\usepackage{pgflibraryarrows}
\usepackage{pgflibrarysnakes}
\usepackage{tikz-3dplot}

\usepackage[all]{xy}

\usepackage{enumitem}

\usepackage{color}

\newtheorem{theorem}{Theorem}
\newtheorem{definition}{Definition}
\newtheorem{lemma}[theorem]{Lemma}
\newtheorem{corollary}[theorem]{Corollary}
\newtheorem{proposition}[theorem]{Proposition}

\newtheorem{remark}{Remark}

\newtheorem{example}{Example}
\newtheorem{notation}{Notation}
\newtheorem{warning}{Warning}

\numberwithin{equation}{section}

\renewcommand{\(}{\begin{equation*}}
\renewcommand{\)}{\end{equation*}}
\newcommand{\bea}{\begin{eqnarray*}}
\newcommand{\eea}{\end{eqnarray*}}

\def\proof {\emph{Proof.}\hspace{7pt}}

\def\endofproof {\hfill{$\Box$}\\}

\newcommand{\onto}{\twoheadrightarrow}

\newcommand{\into}{\hookrightarrow}

\def\hom{{\rm Hom}}
\def\map{{\rm Map}}

\newcommand{\RR}{\ensuremath{\mathbb R}}
\newcommand{\NN}{\ensuremath{\mathbb N}}
\newcommand{\ZZ}{\ensuremath{\mathbb Z}}
\newcommand{\QQ}{\ensuremath{\mathbb Q}}

\newcommand{\CC}{\ensuremath{\mathbb C}}

\newcommand{\sset}{\ensuremath{s\mathscr{S}\mathrm{et}}}

\newcommand{\sh}{\ensuremath{\mathscr{S}\mathrm{h}}}

\newcommand{\cartsp}{{\sf OCart}}

\newcommand{\D}{\ensuremath{\mathscr{D}}}
\newcommand{\E}{\ensuremath{\mathscr{E}}}

\newcommand{\Ab}{\mathscr{A}{\rm b}}

\def\MT{{\sf MT}}

\def\smcatdual{{\rm C^{\infty}Cat}^{\otimes,\vee}}

\def\sm{{\rm C}^{\infty}}
\def\hocolim{\mathop{\rm hocolim}}
\def\holim{\mathop{\rm holim}}
\def\colim{\mathop{\rm colim}}
\def\sing{\mathord{\rm sing}}

\def\id{{\rm id}}
\def\Sp{\mathscr{S}{\rm p}}
\def\sing{{\rm Sing}}

\def\smcat{{\rm C}^{\infty}{\sf Cat}^{\otimes}}

\def\PSh{\mathscr{P}{\rm Sh}}

\def\op{{\rm op}}

\def\GL{{\rm GL}}
\def\deloop{{\rm B}}

\def\cC{{\sf C}}

\def\FEmb{{\sf FEmb}}

\def\Bord{\mathfrak{Bord}}
\def\Riem{\mathrm{Riem}}
\def\RelCat{\mathscr{R}{\rm el}\mathscr{C}{\rm at}}

\def\SCat{\mathscr{SC}{\rm at}}
\def\sS{{\rm s}\mathscr{S}}

\def\downv{\rotatebox{270}{$\dashv$}}
\def\E{{\sf Isot}}
\def\Ei{\mathcal{I}{\rm sot}}
\def\csp{{\sf B}}
\def\D{{\sf Def}}

\def\Icsp{\mathcal{B}}
\def\Ho{\mathop{\rm Ho}}
\def\dg{\mathfrak{d}}
\def\Con{\mathcal{C}}
\def\R{{\sf Res}}

\title{Deformation classes of invertible field theories and the Freed--Hopkins conjecture}
\author{Daniel Grady}
\date{\today}

\begin{document}

\begin{abstract}
We prove a conjecture of Freed and Hopkins, which relates deformation classes of reflection positive, invertible, $d$-dimensional extended field theories with fixed symmetry type to a certain generalized cohomology of a Thom spectrum. Along the way, we establish several results, including the construction of a smooth variant of the Brown--Comenetz dual of the sphere spectrum and a calculation of the ``deformation type'' of the extended geometric bordism category.
\end{abstract}

\maketitle
\tableofcontents

\section{Introduction}
This paper is concerned with the classification of certain deformation classes of invertible field theories. In the seminal work of Freed--Hopkins \cite{FreedHopkins}, a mechanism was proposed for translating between certain lattice systems and corresponding low energy effective field theories. What is so useful about this dictionary, is that the deformation classes of these low energy effective field theories (or rather their Wick-rotated counterparts) are calculable by means of stable homotopy theory. 

The main result of \cite{FreedHopkins} is the classification of deformation classes of reflection positive invertible extended $d$-dimensional \emph{topological} field theories with fixed symmetry type $H_d$. The classification is in terms of the torsion in some generalized cohomology of the Thom spectrum associated to the stable symmetry type $H=\lim_{d\to \infty}H_d$. Such field theories do not depend on any continuously varying background fields, such as metrics, conformal structures, or connections on principal bundles. It was conjectured in \cite{FreedHopkins}, that passing from torsion classes to all classes should yield a classification that holds for a more general class of field theories, which include Riemannian metrics and connections on principal bundles as background fields.  A positive proof of this conjecture would justify the computations of topological phases for various symmetry types in \cite[Section 9.3]{FreedHopkins}, which have already had a considerable impact in condensed matter physics \cite{ThorngrenElse,KapustinThorngren,WangGu,GaiottoJohnson-Freyd}. The main result of this paper is a proof of  \cite[Conjecture~8.37]{FreedHopkins}.
 \begin{theorem}\label{freedhopconj}
Let $H=\colim_{d\to \infty}H_d$ be a stable symmetry type in the sense of Freed–Hopkins \cite{FreedHopkins}.
There is a bijective correspondence
\begin{equation}\left\{
\parbox{6.5cm}{\itshape deformation classes of reflection positive invertible $d$-dimensional extended field theories with symmetry type $(H_d,\rho_d)$}
\right\}\cong [MTH,\Sigma^{d+1}I_{\ZZ(1)}] \label{fhconj}
\end{equation}
where $MTH$ is the Thom spectrum associated to the stable symmetry type and $\Sigma^{d+1}I_{\ZZ(1)}$ is the $(d+1)$-fold suspension of the Anderson dual of the sphere.
\end{theorem}

What allows us to prove this conjecture is our recent work with Pavlov \cite{GradyPavlov}, which provides both a construction of the geometric bordism category and a proof that functorial field theories defined using this bordism category are fully local. Our fully-extended bordism category is capable of handling geometric structures like metrics, principal bundles with connections, conformal structures, and much more. In the present work, our main contribution is the following theorem, which we show implies Theorem \ref{freedhopconj}. 

\begin{theorem}\label{hopfredintro}
Let $H_d$ be a compact Lie group and let $\rho_d:H_d\to {\rm O}(d)$ be a homomorphism. 
In particular, we can take $(H_d,\rho_d)$ to be a symmetry type in the sense of Freed–Hopkins \cite[Definition~2.4]{FreedHopkins}. Consider the geometric structures, $\mathcal{H}_d$, $\mathcal{H}_d^{\nabla}$ and $\mathcal{H}_d^{\rm fl}$, of smooth (Definition \ref{smoothstr}), differential (Definition \ref{diffstr}), and flat (Definition \ref{flatstr}) $(H_d,\rho_d)$-structures (resp.). 
The following spaces are isomorphic in the homotopy category $\Ho \mathscr{S}$:
\begin{enumerate}
\item Smooth deformations of field theories with smooth $(H_d,\rho_d)$-structure:
$$\mathscr{I}_{d}(\mathcal{H}_d):={\mathcal{F}}{\rm un}^{\otimes}(\Bord_d^{\mathcal{H}_d},\Sigma^dI_{\CC^{\times}})$$
\item Smooth deformations of field theories with differential $(H_d,\rho_d)$-structure:
$$\mathscr{I}_{d}(\mathcal{H}^{\nabla}_d):= {\mathcal{F}}{\rm un}^{\otimes}(\Bord_d^{\mathcal{H}^{\nabla}_d},\Sigma^dI_{\CC^{\times}})$$
\item Smooth deformations of field theories with flat $(H_d,\rho_d)$-structure:
$$ \mathscr{I}_{d}(\mathcal{H}^{\rm fl}_d):={\mathcal{F}}{\rm un}^{\otimes}(\Bord_d^{\mathcal{H}^{\rm fl}_d},\Sigma^dI_{\CC^{\times}})$$
\item The space of morphisms of spectra:
$$\map(\Sigma^d MTH_d, \Sigma^{d+1}I_{\ZZ(1)}),$$
\end{enumerate}
Here $MTH_d$ denotes the Madsen--Tillmann spectrum \cite[Section~3]{GMTW}, $I_{\ZZ(1)}$ is the Anderson dual of the sphere spectum, and $I_{\CC^{\times}}$ denote the smooth Brown--Comenetz dual of the sphere spectrum  (Definition \ref{smanddu}). The bordism category that appears was defined in Grady--Pavlov \cite{GradyPavlov}, which is briefly recalled in Section 4. The object $\mathcal{F}{\rm un}^{\otimes}$ is the deformation space of all symmetric monoidal functors.
\end{theorem}

In the above theorem, the geometric structure $\mathcal{H}_d^{\nabla}$ is our formulation of what is called a \emph{differential} $H_d$-structure in Freed--Hopkins \cite{FreedHopkins}. The geometric structure $\mathcal{H}_d$ encodes smooth structure, without connections. The geometric structure $\mathcal{H}_d^{\rm fl}$ encodes flat (or integrable) $H_d$-structures. 

The object $\mathcal{F}{\rm un}^{\otimes}$ encodes deformations of functorial field theories. For two smooth symmetric monoidal $(\infty,d)$-categories $\mathscr{C}$ and $\mathscr{D}$, the object $\mathcal{F}{\rm un}^{\otimes}(\mathscr{C},\mathscr{D})$ is obtained by taking the mapping space between the corresponding \emph{deformation objects} $\csp(\mathscr{C})$ and $\csp(\mathscr{D})$, which are symmetric monoidal $(\infty,d)$-categories. The deformation object  $\csp(\mathscr{C})$ is constructed by evaluating $\mathscr{C}$ on smooth simplices and gluing these simplices together by maps of symmetric monoidal $(\infty,d)$-categories (see Definition \ref{struts}). It encodes the geometric homotopy type of a sheaf of symmetric monoidal $(\infty,d)$-categories. 

What is somewhat unexpected is that field theories of type (3) in Theorem \ref{hopfredintro} are also equivalent to the mapping space in (4). The types of structures encoded by $\mathcal{H}_d^{\rm fl}$ are families of \emph{flat} (or \emph{integrable}) $H_d$-structures. The reason that this is unexpected is that one might guess that the relationship between the spaces (3) and (2) should be similar to the relationship between the moduli space of flat $H_d$-bundles and the moduli space of all $H_d$-bundles with connection. Of course, the moduli space of flat $H_d$-bundles is quite different than the moduli space of all principal $H_d$-bundles with connection. What makes this identification work is the presence of isotopies in the bordism category and the presence of parametrizing families of geometric structures. Allowing pullbacks by isotopies for geometric structures has a surprisingly nontrivial effect on the homotopy type of the smooth bordism category (see Example \ref{flatbunexample}). 


 One can regard the equivalence of (3) and (4) as a manifestation of Gromov's h-principle. In fact, the proof of Proposition \ref{shapestr}, which is the crucial ingredient in this equivalence, uses a common technique to prove the h-principle by showing that a map of sheaves induces an isomorphism on relative concordance classes.

The outline of the paper is as follows. In Section \ref{relcat}, we provide the categorical background for the present work, including a brief review of relative categories and simplicial presheaves. In Section \ref{geometric.struct}, we recall our notion of geometric structures from Grady--Pavlov \cite{GradyPavlov}. We also introduce new geometric structures which will be needed in the present work. In Section \ref{bordrefl}, we  briefly recall the bordism category constructed in \cite{GradyPavlov} and discuss reflection positive field theories and the geometric structures that are necessary to encode them. In Section \ref{structdeform} we introduce several deformation space constructions, which produce spaces from various types of geometric data. In Section \ref{stablehom}, we discuss modifications of the stable homotopy theory used in \cite{FreedHopkins} to the setting of sheaves. Finally, in Section \ref{freedhopsec}, we prove the main theorem and its corollary Theorem \ref{freedhopconj}.

\vspace{.2cm}

\paragraph{\bf Acknowledgments} The author would like to thank Arun Debray, for useful conversations on reflection positivity and for feedback on the first draft of this article. The author would also like to thank David Reutter for correcting a misunderstanding about how positivity enters the proof of Freed and Hopkins. Finally, the author would like to thank Dmitri Pavlov for pointing out the recent work of Patchkoria--Pstragowski \cite{PatchkoriaPstragowski}, whose version of the Brown representability theorem is used in the present work, and for suggesting a much simpler proof of Proposition \ref{w.e.bords}, which no longer uses the geometric cobordism hypothesis.

\section{Homotopy sheaves via relative categories}\label{relcat}

In this section, we review the notion of a  $(\infty,1)$-sheaf (a.k.a $\infty$-sheaf). We will consider $(\infty,1)$-sheaves with values in spaces and in spectra. 

\begin{warning}
Throughout the paper, we use the notation $(\infty,1)$ as a subscript when referring to homotopical structures on a category (see Notation \ref{notation}). We emphasize that this \emph{does not} indicate that we are using quasi-categories as presentations. Instead, all of our $(\infty,1)$-categories are relative categories. The notation $(\infty,1)$ is justified, since the category of relative categories, equipped with the Barwick--Kan model structure, indeed serves as a model for $(\infty,1)$-categories (as was shown in Barwick--Kan \cite{BarwickKan}). 
\end{warning} 

We assume the reader has some basic familiarity with Quillen's theory of model categories. Some excellent references for model categories are Hovey \cite{Hovey}, Dwyer--Spalinski \cite{DwyerSpalinski}, Hirschhorn \cite{Hirschhorn}. See also Riehl \cite[Chapter 3]{Riehl} and Balchin \cite{Balchin}. We will also often use Bousfield localization to obtain new model structures from existing ones. An excellent reference for Bousfield localization of model categories is Hirschhorn \cite{Hirschhorn}.

Roughly speaking, an $(\infty,1)$-sheaf is a presheaf with values in spaces that satisfies a sheaf condition, but only up to higher homotopy coherence. More precisely, if ${\sf C}$ is a small category equipped with a coverage and $\{U_a\to U\}_{a\in A}$ is a covering family of some object $U\in {\sf C}$, then an $(\infty,1)$-sheaf $X$ is a presheaf on ${\sf C}$ with values in simplicial sets (equipped with the usual Kan--Quillen model structure) such that the canonical map
$$X(U)\overset{\simeq}{\to} \holim_{[n]\in \Delta} \prod_{\alpha:[n]\to A}X(U_{\alpha})$$
is an equivalence. Here we are using multi-index notation $U_{\alpha}:=U_{\alpha_0}\cap U_{\alpha_1}\cap \hdots U_{\alpha_n}$. The face maps on the right are given by restricting to higher-fold intersections in the obvious way.  

This homotopical sheaf gluing condition can be rephrased by saying that $X$ is a simplicial presheaf that is local (in the sense of Bousfield localization) with respect to the maps
$$\hocolim_{[n]\in \Delta^{\rm op}}\coprod_{\alpha:[n]\to A}U_{\alpha}\to U.$$
Of course, one needs to choose a specific presentation of the homotopy colimit, but the $(\infty,1)$-sheaf property is independent of such a choice. 

Very often we will need to prove some statement about homotopy limits or colimits, or about functors that preserve these. Here we observe that (up to weak equivalence) homotopy limits and colimits do not depend on a specific choice of model structure on functors, but only on the underlying relative category. Indeed, deriving the colimit functor with respect to two different model structures, both having the same weak equivalences, produces weakly equivalent results. Both of these deserve to be called the homotopy colimit.  

For this reason, our $(\infty,1)$-categories will merely be relative categories that admit a model structure whose weak equivalences agree with the given ones. Such categories are known to be fibrant objects in Barwick--Kan model structure on relative categories (see Meier \cite{Meier}). That relative categories indeed model all $(\infty,1)$-categories was proved in Barwick--Kan \cite{BarwickKan}. This gives us some flexibility when proving claims about homotopy limits or colimits in that we are free to fix any model structure we like, provided the underlying relative category agrees with the given one.  

\begin{definition}
A \emph{relative category} is a pair $({\sf C},{\sf W})$ where ${\sf W}$ is a wide subcategory of ${\sf C}$, meaning it contains all objects of ${\sf C}$. We call the morphisms in ${\sf W}$ \emph{weak equivalences}. A morphism of relative categories $f:({\sf C},{\sf W})\to ({\sf D},{\sf V})$ is a functor $f:{\sf C}\to {\sf D}$ that sends weak equivalences to weak equivalences. 
\end{definition}  

The homotopy category of a relative category $({\sf C},W)$ is defined simply as the localization (as a 1-category) 
$${\rm Ho}({\sf C}):={\sf C}[W^{-1}].$$ 
When $({\sf C},W)$ admits a calculus of fractions, then the localization can be computed by taking morphisms to be equivalence classes of spans. It is easy to see that an adjunction between two relative categories, where both functors preserve weak equivalences, induces a corresponding adjunction at the level of homotopy categories. However, this is not quite enough in practice. We will often need adjunctions that are induced by Quillen adjunctions of model categories, and these need not preserve all weak equivalences. This is easily circumvented by the following observation. 

For any model category ${\sf C}$ with cofibrant and fibrant replacement functors $Q$ and $R$, both $Q$ and $R$ induce functors that are naturally isomorphic to the identity at the level of the homotopy category ${\sf C}[W^{-1}]$. Let $f\dashv g$ be a Quillen equivalence between two model categories ${\sf C}$ and ${\sf D}$. Let $\mathbb{L}f:=f\circ Q$ and $\mathbb{R}g:=g\circ R$ be the left and right derived functors. Then since both $\mathbb{L}f$ and $\mathbb{R}g$ preserve all weak equivalences, there is an induced adjunction at the level of homotopy categories
\begin{equation}\label{adjfuns}
\xymatrix{
{\sf C}[W^{-1}]\ar@<.15cm>[r]^{\mathbb{L} f}\ar@{}[r] & \ar@<.15cm>[l]^{\mathbb{R}g} {\sf D}[V^{-1}].
}
\end{equation}

Recall that every relative category can be simplicially enriched via hammock localization \cite{DwyerKan}. In general, if $c,d\in ({\sf C},W)$, we will denote the simplicial set of morphisms obtained by hammock localization as $\map(c,d)$. 

\begin{definition}
Given two morphisms of relative categories $f:({\sf C},W)\to ({\sf D},V)$ and $g:({\sf D},V)\to ({\sf C},W)$, we say that $f$ and $g$ form a \emph{homotopy adjunction}, with $f$ the \emph{left homotopy adjoint} of $g$, if there is an isomorphism in $\Ho \mathscr{S}$:
$$
\map(fc,d)\simeq \map(c,gd)
$$
which is natural in both $c$ and $d$.
\end{definition}

The adjunction \eqref{adjfuns} can be improved to the following.
\begin{lemma}
Let ${\sf C}$ and ${\sf D}$ be model categories and let $({\sf C},W)$ and $({\sf D},V)$ denote the corresponding underlying relative categories. Let
$$
\xymatrix{
{\sf C}\ar@<.15cm>[r]^{f}\ar@{}[r]|{\downv} & \ar@<.15cm>[l]^{g} {\sf D}
}
$$
be a Quillen adjunction. Then $f$ and $g$ induce a corresponding homotopy adjunction $\mathbb{L}f\dashv \mathbb{R} g$.
\end{lemma}
\proof
This follows from Mazel-Gee \cite[Theorem 2.1]{Mazel-Gee}.
\endofproof

Let us denote the category of all relative categories by $\RelCat$. We denote the category of simplicial objects in simplicial sets by $\sS$, and the category of simplicial categories by $\SCat$. In Barwick--Kan \cite{BarwickKan}, it was shown that $\RelCat$ admits a model structure that is transferred from Rezk's model structure on complete Segal spaces via a functor $N_{\xi}:\RelCat\to \sS$. The functor $N_{\xi}$ is equivalent to Rezk's classifying diagram construction $N:\RelCat\to \sS$ (Rezk \cite{Rezk}). In Dwyer--Kan \cite{DwyerKan}, the hammock localization functor
$$L_H:\RelCat\to \SCat$$
was defined. The functor $L_H$ sends a relative category $({\sf C},{\sf W})$ to a simplicial category $L_H{\sf C}$, whose objects are the same as those of ${\sf C}$. The simplicial set of morphisms is given by the hammock construction \cite{DwyerKan}. Finally, in Barwick--Kan \cite{BarwickKan2}, it was shown that the weak equivalences in $\RelCat$ are precisely those morphisms that are sent to Dwyer--Kan equivalences under $L_H$. Hence, the weak equivalences in the Barwick--Kan model structure on $\RelCat$ are precisely the Dwyer--Kan weak equivalences. 

Combining these facts together, we have the following proposition.
\begin{proposition}
Let $L_H:\RelCat\to \SCat$ and $N:\RelCat\to \sS$ be as above. Let $\sset_{J}$ denote the category of simplicial sets, equipped with the Joyal model structure. Let $N^{hc}:\SCat\to \sset_{J}$ denote the homotopy coherent nerve functor \cite[Section 1.51.]{Lurie.HTT} and let $\RR N^{hc}$ denote the corresponding right derived functor. 

Let $f:{\sf C}\to {\sf D}$ be a left Quillen equivalence between model categories ${\sf C}$ and ${\sf D}$. Let $\mathbb{L}f:=f\circ Q$ denote the left derived functor. Then the morphism of relative categories
$$\mathbb{L}f:{\sf C}\to {\sf D},$$
satisfies the following equivalent conditions:
\begin{enumerate}
\item $\mathbb{L}f:{\sf C}\to {\sf D}$ is a Dwyer--Kan equivalence of relative categories;
\item $L_H(\mathbb{L}f):L_H{\sf C}\to L_H{\sf D}$ is a Dwyer--Kan equivalence of simplicial categories;
\item $N(\mathbb{L}f):N({\sf C})\to N({\sf D})$ is a weak equivalence of complete Segal spaces;
\item $\RR N^{hc}L_H(\mathbb{L}f): \RR N^{hc} L_H{\sf C}\to \RR N^{hc}L_H{\sf D}$  is a weak equivalence of quasi-categories.
\end{enumerate}
\end{proposition} 
\proof
First observe that if $f$ is a left Quillen equivalence, then so is $\mathbb{L}f$. By Mazel-Gee \cite{Mazel-Gee}, the functor $\mathbb{L}f$ induces an equivalence of quasi-categories $\RR N^{hc}L_H(\mathbb{L}f):\RR N^{hc}L_H{\sf C}\to \RR N^{hc}L_H{\sf D}$. Let $R$ denote the fibrant replacement functor in the Bergner model structure on $\SCat$. Since the functor $N^{hc}:\SCat\to \sset_{J}$ is a right Quillen equivalence, and all objects in the Joyal model structure are cofibrant, it follows the the components of the counit $\epsilon$ in the diagram 
$$
\xymatrix{
N^{hc} R(L_H{\sf C})\ar[rr]^{\RR N^{hc}L_H(\mathbb{L}f)}_-{\simeq} \ar[d]^-{\epsilon} &&  N^{hc} R(L_H{\sf D})\ar[d]^-{\epsilon}
\\
R(L_H{\sf C})\ar[rr]^-{R(L_H(\mathbb{L}f))} && R(L_H{\sf D})
\\
L_H{\sf C}\ar[rr]\ar[u]^-{\simeq} && L_H{\sf D}\ar[u]^-{\simeq}
}
$$
are derived, hence are weak equivalences. By 2-out-of-3, it follows also that $L_H(\mathbb{L}f)$ is an equivalence, so that $(4)\Rightarrow (2)$. The implication $(2)\Rightarrow (4)$ follows immediately from the fact that $N^{hc}$ is right Quillen. The implications $(2)\Leftrightarrow (1)$ are a tautology. The implications $(1)\Leftrightarrow (3)$ were proved in \cite{BarwickKan2}.
\endofproof

We now introduce some notation, which will be used throughout the paper.

\begin{notation}\label{notation}
We adopt the following notation and conventions. 

\begin{enumerate}
\item[(i)] Let ${\sf C}$ be a relative category and let $c,d\in {\sf C}$.  We denote the simplicial mapping space, given by hammock localization, as $\map(c,d)$. In the case where ${\sf C}$ admits the structure of a simplicial model category, this mapping space is weakly equivalent to the derived enriched mapping space. 
\item[(ii)]
We denote the category of simplicial sets, equipped with its usual weak equivalences, by $\mathscr{S}$. We call the objects of $\mathscr{S}$ \emph{spaces}.
\item[(iii)] Let ${\sf C}$ be a small category. We let $\PSh_{(\infty,1)}({\sf C})$ denote the relative category of simplicial presheaves on ${\sf C}$, with weak equivalences given by those natural transformations that are componentwise weak equivalences in $\mathscr{S}$.
\item[(iv)] Let $({\sf C},J)$ be a Grothendieck site. We let $\sh_{(\infty,1)}({\sf C})$ denote the relative category of simplicial presheaves on ${\sf C}$ with weak equivalence given by $S$-local weak equivalences, where $S$ is the set of morphisms of the form \eqref{cechmaps}, obtained as follows:
 
 \begin{itemize}
 \item Fix a cover $\{U_{a}\to U\}_{a\in A}$ of an object $U\in {\sf C}$. Let $\check{C}(\{U_{a}\})$ denote the simplicial presheaf in $\PSh_{(\infty,1)}({\sf C})$ given by 
$$\check{C}(\{U_{a}\}):[n]\mapsto \coprod_{\alpha:[n]\to A} U_{\alpha},$$
where $U_{\alpha}=U_{\alpha_0}\times_{U} U_{\alpha_1}\times_{U} \hdots U_{\alpha_n}$ (the $n$-fold pullback). The face and degeneracy maps are the obvious ones, induced by the canonical maps defining the pullback. There is a canonical morphism of simplicial presheaves
\begin{equation}
\check{C}(\{U_{a}\})\to U,\label{cechmaps}
\end{equation}
induced by the the  maps $U_{\alpha_0}\times_{U} U_{\alpha_1}\times_{U} \hdots U_{\alpha_n}\to U_{\alpha_0}\to U$. 

More generally, if ${\sf V}$ is a closed symmetric monoidal category, then ${\sf C}$ can be enriched in ${\sf V}$ by ${\sf C}_{\sf V}(c,d):={\sf C}(c,d)\otimes {\bf 1}$, where ${\bf 1}$ is the monoidal unit. In this case, we let $\mathscr{S}{\rm h}_{(\infty,1)}({\sf C},{\sf V})$ denote the relative category with weak equivalences given by $S$-local weak equivalences, where $S$ is the set of morphisms of the form \eqref{cechmaps}, where now $U$ and $U_{\alpha}$ are embedded via the enriched Yoneda embedding. 
\end{itemize}
\item[(v)] We denote the relative category of simplicial symmetric spectra, equipped with stable weak equivalences, by $\Sp$.
\item[(vi)] We denote by $\Sp_{\geq 0}$ the relative category of connective spectra, equipped with local weak equivalences in the Bousfield--Friedlander model structure on $\Gamma$-spaces.
\item[(vi)] Let ${\sf C}$ be a small category.  We let $\PSh_{(\infty,1)}({\sf C};\Sp)$ denote the relative category of presheaves on ${\sf C}$ with values in $\Sp$, with weak equivalences given by those natural transformations that are componentwise weak equivalences in $\Sp$ .
\item[(vii)] Let $({\sf C},J)$ be a Grothendieck site. We let $\sh_{(\infty,1)}({\sf C};\Sp)$ denote the relative category of simplicial presheaves on ${\sf C}$ with values in $\Sp$, with weak equivalence given by $S$-local weak equivalences, where $S$ is the set of morphisms of the form
\begin{equation}\label{cechmorph}
\{F_n(\check{C}(\{U_{a}\}))\to F_n(U) \mid n\in \NN\}
\end{equation}
where $F_n$ is left adjoint to the evaluation functor ${\rm Ev}_n:\PSh_{(\infty,1)}({\sf C};\Sp)\to \PSh_{(\infty,1)}({\sf C})$, which sends a presheaf of spectra $X$ to the presheaf of spaces given by taking its $n$-th level. Morally, the localization inverts the tensoring of {\v C}ech morphisms \eqref{cechmaps} with all desuspensions of the sphere (see \cite[Appendix B]{Dugger2}).
\end{enumerate}
\end{notation}

There are two variants of $(\infty,1)$-sheaves that will be used in this paper. The first is $(\infty,1)$-sheaves on the category of open subsets of cartesian spaces, which will be responsible for keeping track of smooth families. The second is $(\infty,1)$-sheaves on fiberwise open embeddings (see Section \ref{geometric.struct}), which will be responsible for encoding geometric structures on bordisms. We begin by defining the category of open subsets of cartesian spaces.

\begin{definition}
We denote by $\cartsp$ the category whose objects are open subsets $U\subset \RR^n$, for some $n\in \NN$. Morphisms are given by smooth maps. The category $\cartsp$ admits the structure of a Grothendieck site, by taking covering families to be the usual covers by open subsets $\{U_{a}\into U\}_{a\in A}$ where $U,U_{a}\in \cartsp$, for all $a\in A$.
\end{definition}

We will usually refer to $(\infty,1)$-sheaves on the above site as \emph{smooth $(\infty,1)$-sheaves}, or occasionally just \emph{smooth sheaves}. Sometimes we will refer to certain objects as being smooth, e.g., a \emph{smooth symmetric monoidal $(\infty,d)$-category with duals}. In all such cases, the object is a presheaf on $\cartsp$ with values in some combinatorial model category.

\section{Geometric structures for bordisms}\label{geometric.struct}

\subsection{Recollections on reflection positivity}

Here we recall the notion of reflection positivity for invertible field theories. We defer the reader to Freed-Hopkins \cite{FreedHopkins} for a much more thorough treatment. In quantum mechanics, one has a complex separable Hilbert space $\mathcal{ H}$, equipped with a self-adjoint operator $H$ (the Hamiltonian). Time evolution of the quantum system is provided by the 1-parameter group of unitary operators 
\begin{equation}\label{unitarytheory}
U(t)=e^{-itH}, 
\end{equation}
where we set Planck's constant $\hbar=1$. If $H$ is positive definite, then we can extend the group $U(t)$ to the lower half plane of the complex plane $\mathscr{T}=\RR-i\RR_{>0}\subset \CC$. Restricting to the imaginary time axis gives a semigroup
\begin{equation}\label{reflectionpositive}
R(\tau)=e^{-\tau H}, ~~\tau>0.
\end{equation}
Going from \eqref{unitarytheory} to \eqref{reflectionpositive} is called \emph{Wick rotation}.

Recall that for a complex vector space $\mathcal{H}$ over $\CC$, the \emph{conjugate vector space} $\overline{\mathcal{H}}$ of $\mathcal{H}$ is the vector space whose underlying abelian group is the same as that of $\mathcal{H}$, but scalar multiplication is given by conjugated scalar multiplciation $\lambda \cdot v:=\overline{\lambda} v$. Recall that a real structure on $\mathcal{H}$ is a linear map $\sigma:\mathcal{H}\to \overline{\mathcal{H}}$ satisfying $\overline{\sigma}\sigma=1$. A canonical example is given by taking $\mathcal{H}$ to be complex-valued $L^2$ functions with  pointwise complex conjugation as the real structure. A bounded operator $R$ on $\mathcal{H}$ canonically gives rise to a bounded operator $\overline{R}$ on $\overline{\mathcal{H}}$ by defining $\overline{R}(v)=R(v)$ for all $v\in \overline{\mathcal{H}}$. A bounded operator $R$ on $\mathcal{H}$ is called \emph{real} if and only if $\sigma R=\overline{R}\sigma$.  

Under Wick rotation, time reversal corresponds to conjugation of the Wick rotated operators and the unitarity of the semigroup \eqref{unitarytheory} manifests itself as the \emph{reality} of the semigroup \eqref{reflectionpositive} (see Freed–Hopkins \cite[Section 3]{FreedHopkins}). 
Time reversal manifests itself in the bordism category as change of orientation. A field theory with \emph{reflection structure} is thus a field theory that sends orientation reversal to complex conjugation of vector spaces and linear maps, hence a $\ZZ/2$-equivariant functor on the bordism category. 
The word \emph{positive} in \emph{reflection positive} refers to positive definiteness of the Hamiltonian.

As in \cite{FreedHopkins}, we use Segal's axioms of functorial field theories to formalize the Wick rotated theory. Recall that a \emph{bordism} is given by an $n$-dimensional smooth manifold $\Sigma_n$ with boundary $\partial \Sigma_n$, possibly equipped with some additional structure, such as a Riemannian metric, conformal structure, etc. An additional datum is required to specify which connected components of $\partial \Sigma_n$ should be regarded as the source of the bordism, and which should be regarded as the target. In the fully extended setting, bordisms are manifolds with corners: objects are 0-dimensional manifolds, morphisms are 1-dimensional manifolds with boundary, 2-morphisms are 2-dimensional manifolds with corners, etc. 

In Grady--Pavlov \cite{GradyPavlov}, a fully extended bordism category was defined. This bordism category is capable of encoding a very large class of geometric structures. In particular, we can encode the differential $H_n$-structures of \cite{FreedHopkins}. We first recall the notion of a symmetry type for a quantum field theory \cite[Definition 2.4]{FreedHopkins}.

\begin{definition}\label{symmtype} A \emph{symmetry type} is a pair $(H_n,\rho_n)$, where $H_n$ is a compact Lie group and $\rho_n:H_n\to {\rm O}(n)$ is a homomorphism whose image contains ${\rm SO}(n)\subset {\rm O}(n)$. 

\end{definition}

The motivation for this definition comes from quantum field theory, where it is often assumed that the Poincar\'e group is a subgroup of some larger (unbroken) symmetry group. A more natural assumption is that there is a homomorphism from this larger group to the Poincar\'e group, whose image contains the connected component of the identity. The hypothesis that the image of $\rho_n$ contains ${\rm SO}(n)$ is a reformulation of this for the Wick-rotated theory (see \cite[Section 2.1]{FreedHopkins}).

For the types of field theories we are considering here, the background is not fixed, but is dynamical. The fields are ``coupled to background gravity", which means that the symmetries of the theory are not global symmetries, but are local on spacetime. Such local symmetries are encoded as tangential structures, i.e., a reduction of the structure group of the bundle of orthonormal frames (recall we are in the Euclidean signature) to $H_n$. 

\subsection{Geometric structures as sheaves}
We now recall some basic definitions from Grady--Pavlov \cite{GradyPavlov}, which are necessary to define our notion of smooth and differential $H_d$-structures. We begin be recalling the definition of geometric structures on our bordisms \cite{GradyPavlov}. We then recall the bordism category $\mathfrak{Bord}^{\mathcal{S}}_d$, defined in \cite{GradyPavlov}. 

A geometric structure on a bordism, for example a Riemannian metric, must satisfy a basic set of axioms in order to be regarded as a morphism in the bordism category. In any category a morphism must have a source and a target. Moreover, morphisms can be composed to obtain new morphisms. Composition of bordisms with geometric structures requires some care, as one cannot glue arbitrary geometric structures, such as a Riemannian metric, along a codimension 1 manifold. 

As realized by Stolz--Teichner \cite{StolzTeichner}, one of the most efficient ways to circumvent this issue is to keep track of the germ of a codimension 1-submanifold. Geometric structures should be required to agree on an open neighborhood of the codimension 1 manifold, so that gluing is possible. With these observations in mind, we see that geometric structures should satisfy the following two properties:
\begin{enumerate}\label{bordprops}
\item Geometric structures can be restricted along open embeddings, e.g., an open neighborhood of a codimension 1 submanifold. 
\item Given an open cover $\{U,V\}$ of a bordism $\Sigma$ and geometric structures $g_{U}$ on $U$ and $g_{V}$ on $V$, if the restrictions of $g_{U}$ to $U\cap V$ and $g_{V}$ to $U\cap V$ agree, then there is a geometric structure $g_{\Sigma}$ on $\Sigma$ whose restrictions to $U$ and $V$ agree with $g_{U}$ and $g_{V}$, respectively.
\end{enumerate} 
After a moment's thought, one realizes that this is saying precisely that a geometric structure should form a sheaf on the site of $d$-dimensional smooth manifolds with open embeddings as morphisms, where covering families are given by ordinary open covers. 

In order to encode smooth families, a morphism in the bordism category should not just be a $d$-dimensional manifold with boundary, but rather a smooth family of such manifolds \cite{StolzTeichner}. Keeping track of smooth families, both on the bordism side and in the target, forces field theories to be smooth in the sense that a smooth family of bordisms must be sent to a smooth family of values. 

This motivates the following definition. 

\begin{definition}\label{femb}
Fix $d\geq 0$. Let $\FEmb_d$ be the site with objects submersions $p:M\to U$, with $d$-dimensional fibers and $U$ an object in $\cartsp$. 
A morphism $(f,g):(p:M\to U)\to (q:N\to V)$ is a pair of smooth maps $f:M\to N$ and $g:U\to V$ making the diagram 
$$
\xymatrix{
M\ar[r]^-{f}\ar[d]^-{p} & N\ar[d]^-{q}
\\
U\ar[r]^-{g} & V
}
$$
commute. Moreover, we require that for all $u\in U$, the restricted map $f\vert_{p^{-1}(u)}:p^{-1}(u)\to q^{-1}(g(u))$ is an open embedding (note that this condition is vacuous if the fiber $p^{-1}(u)$ is empty). 

Covering families are given by a collection of morphisms
$$\left\{\vcenter{
\xymatrix{
M_{a}\ar[r]\ar[d] & M\ar[d]^-{p}\\
U_{a}\ar[r] & U
}}
\right\}_{a\in A}
$$ 
such that $\{M_{a}\}_{a\in A}$ is an open cover of $M$.
\end{definition}

One should think that an object $(M\to U)\in \FEmb_d$ will play the role of a smooth family of bordisms, or rather the ambient manifolds containing these bordisms, as parametrized by points of $U$. The two properties (1) and (2) are sufficient to encode structures like Riemannian metrics. However, they are not sufficient to encode more exotic structures, such as principal $G$-bundles with connections, since these can be glued from local data, provided that a compatible collection of  \emph{isomorphisms} is specified on the intersections. This weaker sheaf gluing condition is precisely encoded by the $(\infty,1)$-sheaf condition. This motivates the following.

\begin{definition}
\label{geometric.structure}
Fix $d\geq 0$. A {\it fiberwise $d$-dimensional geometric structure\/} is an $(\infty,1)$-sheaf on $\FEmb_d$, i.e., it is a local object in $\sh_{(\infty,1)}(\FEmb_d)$.
\end{definition}

\begin{example}\label{targetman}
Let $X$ be a smooth manifold, of any dimension. We can regard $X$ as a simplicial presheaf on ${\sf FEmb}_d$, for any fixed $d\geq 0$, as follows. For a submersion $p:M\to U$ with $d$-dimensional fibers, define
$$X(p:M\to U):=C^{\infty}(M,X).$$
Since smooth functions can be pulled back along fiberwise open embeddings, this indeed defines a presheaf on ${\sf FEmb}_d$. It is a sheaf, since the property of being a smooth function on $M$ is a local property, with respect to open covers of $M$
\end{example}


Recall that a point of an $(\infty,1)$-topos $\mathscr{X}$ is a geometric morphism $x:\mathscr{S}\to \mathscr{X}$. An $(\infty,1)$-topos has enough points if equivalences in $\mathscr{X}$ can be detected stalkwise (at each point $x$). A classic example is when $\mathscr{X}$ is $(\infty,1)$-sheaves on a sober topological space $X$. In this case, the points of $\sh_{(\infty,1)}(X)$ are indexed by the points of $X$. When $X$ is paracompact and has finite covering dimension, then $\sh_{(\infty,1)}(X)$ has enough points (this is Lurie \cite[Theorem 7.2.3.6]{Lurie.HTT}).

We now examine the points of $\sh_{(\infty,1)}(\FEmb_d)$. 

\begin{proposition}\label{stalks}
Fix $d\geq 0$ and let $\sh_{(\infty,1)}(\FEmb_d)$ be the relative category of simplicial presheaves with {\v C}ech local weak equivalences (as usual). Then the following hold.
\begin{itemize}
\item[(i)] The $(\infty,1)$-topos $\sh_{(\infty,1)}(\FEmb_d)$ has distinct points indexed by $n\in \NN$. For each $n$, the corresponding left adjoint $p_n^*$ in the geometric morphism is given by 
$$p_n^*(X)=\colim_{\delta\to 0}X(p_{B^n_{\delta}(0)}:B^d_{\delta}(0)\times B^n_{\delta}(0)\to  B^n_{\delta}(0),)$$
where $B^n_{\delta}(0)\subset \RR^n$ is the open ball of radius $\delta$ centered at $0$ and $p_{B_{\delta}(0)^n}$ is the projection. The colimit is taken over the filtered diagram whose transition maps are induced by the pair $(j^d_{\delta,\delta'}\times j^n_{\delta,\delta'},j^n_{\delta,\delta'})$, where $\delta'>\delta$ and $j^n_{\delta,\delta'}: B_{\delta}^n(0)\into B_{\delta'}^n(0)$ is the canonical inclusion. 

\item[(ii)] The $(\infty,1)$-topos $\sh_{(\infty,1)}(\FEmb_d)$ is hypercomplete. Moreover, the points $\{p_n\}$ are sufficient to detect local weak equivalences. That is, if $f:X\to Y$ induces an equivalence after applying $p^*_n$, for all $n\in \NN$, then $f$ is a weak equivalence. 
\end{itemize}
\end{proposition}
\proof
That each $p^*_n$ is the left adjoint in a geometric morphism is clear, since the colimit defining $p^*_n$ is a filtered colimit of simplicial sets, hence a homotopy filtered colimit, and therefore also homotopy cocontinuous and homotopy left exact. 

We first show that the points $p_n$ detect local isomorphisms in the underlying 1-topos $\sh(\FEmb_d)$. This is a standard exercise, but we will spell out the details. Let $f$ be a morphism inducing an isomorphism on corresponding stalks. We claim that $f$ is a local epimorphism. Fix an arbitrary object $p:M\to U\in \FEmb_d$ and let $n={\rm dim}(U)$. Let $B_{\epsilon}^n(0)$ denote the open $n$-ball of radius $\epsilon$ centered at $0$ and let $q_{\epsilon}:B_{\epsilon}(0)^d\times B_{\epsilon}^n(0)\to B_{\epsilon}^n(0)$ denote the projection. By the constant rank theorem, there exists a covering of $p:M\to U$ of the form 
$$\left\{\vcenter{
\xymatrix{
B^d_{\epsilon_x}(0)\times B^n_{\epsilon_x}(0)\ar[r]^-{\tilde j_x}\ar[d]_{q_{\epsilon_x}} & M\ar[d]^-{p}\\
B^n_{\epsilon_x}(0)\ar[r]^-{j_x} & U
}
}
\right\}_{x\in M}
$$ 
where each $\epsilon_x$ can be made arbitrarily small. Fix a section $s\in Y(p:M\to U)$. Since $f$ is surjective on stalks, it follows that for each $x\in M$, there is $\epsilon_x>0$ and $t_{x}\in X(q_{\epsilon_x})$ such that $f(t_x)=(\tilde j_x,j_x)^*(s)$. Since $\{q_{\epsilon_x},(\tilde j_x,j_x)\}$ form a cover of $p:M\to U$, the sheaf gluing condition gives the existence of $t\in X(p:M\to U)$ that locally maps to $s$ under $f$. Since $s$ was arbitrary, $f$ is a local epimorphism. It is easy to verify that $f$ is a monorphism, using the locality property of sheaves. Hence, $f$ is a local isomorphism.

In the following we switch to quasi-categories and then compare to relative categories in the next paragraph. We claim that the $\infty$-topos $\sh_{\infty}(\FEmb_d)$ is hypercomplete. For this, it suffices to prove that for every $M\to U\in \FEmb_d$, the $\infty$-topos $\sh_{\infty}((\FEmb_d)_{/M\to U})$ has finite homotopy dimension \cite[Proposition 7.2.1.10]{Lurie.HTT}. Consider the morphism of sites $\mathcal{O}{\rm pen}(M)\to (\FEmb_d)_{/M\to U}$ that sends an open subset $V\subset M$ to the diagram 
$$
\xymatrix{
V~\ar[d]^-{p}\ar@{^{(}->}[r] & M\ar[d]^-{p} 
\\
U\ar[r]^{=} & U\;.
}
$$
The corresponding restriction $i^*:\PSh_{\infty}((\FEmb_d)_{/M\to U})\to \PSh_{\infty}(\mathcal{O}{\rm pen}(M))$ sends $\infty$-sheaves to $\infty$-sheaves and commutes with sheafification. Let $F$ be an $n$-connected $\infty$-sheaf on $(\FEmb_d)_{/M\to U}$. Then $i^*F$ is $n$-connected and $i^*F(M)$ is nonempty, by \cite[Theorem 7.2.3.6]{Lurie.HTT}. Hence $F(M\to U)=i^*F(M)$ is also nonempty. Therefore, $\sh_{\infty}((\FEmb_d)_{/M\to U})$ has finite homotopy dimension, as claimed. 

Finally, let $\sh_{(\infty,1)}(\FEmb_d)_{\text{\v Cech}}$ be the relative category where weak equivalences are {\v C}ech local weak equivalences. Let $\sh_{(\infty,1)}(\FEmb_d)_{{\rm Jar}}$ be the relative category with Jardine local weak equivalences. Then by \cite[Proposition 6.5.2.14]{Lurie.HTT}, hypercompleteness of the $\infty$-topos $\sh_{\infty}(\FEmb_d)$ implies that the identity functor 
\begin{equation}\label{cechtojardine}
\sh_{(\infty,1)}(\FEmb_d)_{\text{{\v C}ech}}\to \sh_{(\infty,1)}(\FEmb_d)_{{\rm Jar}}
\end{equation}
is a Dwyer--Kan equivalence. Since the 1-topos $\sh(\FEmb_d)$ has enough points, the weak equivalences in the Jardine model structure are precisely stalkwise weak equivalences. Since local isomorphisms are detected by $\{p_n\}$ in the underlying 1-topos $\sh(\FEmb_d)$, as proved above, it follows from the equivalence \eqref{cechtojardine} that stalkwise equivalences are also weak equivalences in the {\v C}ech model structure. 
\endofproof

\def\Struct{{\sf Struct}}
\begin{notation}\label{fembshvnot}
 Given $d\geq 0$, we denote by 
$$\Struct_{d}=\PSh(\FEmb_d,\sset)_{{\rm Cech}}$$
the relative category whose weak equivalences are given by {\v C}ech local weak equivalences. This relative category has a canonical model category presentation given by taking left Bousfield localization of the injective (or projective) model structure at {\v C}ech covers. 
The existence of the enriched left Bousfield localization exists by Grady--Pavlov \cite[Proposition 2.2.2]{GradyPavlov}.
\end{notation}
 
 \subsection{Differential $H_d$-structures} 
 
Following Grady--Pavlov \cite{GradyPavlov}, we let $\mathbf{B}{\rm GL}(d)$ denote the simplicial presheaf on $\FEmb_d$ that sends a submersion $p:M\to U$ to the groupoid of vector bundles on $M$. A priori, this defines only a pseudo-functor, but it can be made strict (see \cite[Definition 3.2.2]{GradyPavlov} and the construction in Definition \ref{strictconn} below). For a Lie group $G$, we let $\mathbf{B} G$ be the simplicial presheaf that sends a submersion $M\to U$ to 
$$\mathbf{B}G(M\to U):=N\left(\vcenter{\xymatrix{
C^{\infty}(M,G)\ar@<.1cm>[r]\ar@<-.1cm>[r] & \ast
}}\right).$$
Here $N:\mathscr{G}{\rm pd}\to \sset$ denotes the nerve functor. The groupoid that is denoted by the two horizontal arrows (source and target) is the groupoid with a single object and morphisms given by the group $C^{\infty}(M,G)$. This is an objectwise Kan-complex, but it does not satisfy descent on $\FEmb_d$, i.e., we are only taking the trivial bundle and its automorphisms on $M$. To get all principal $G$-bundles from $\mathbf{B}G$, we fibrantly replace in the local projective model structure on simplicial presheaves. 

\begin{remark}
 Alternatively, one could just work with a fibrant model for $\mathbf{B}G$ from the beginning. Our only reason for working with this non-fibrant model is that it allows us to define maps between such moduli stacks in a fairly canonical way. 
\end{remark}

Given a representation $\rho:G\to \GL(d)$, there is a corresponding canonical map $\rho:\mathbf{B} G\to \mathbf{B} \GL(d)$ that sends the trivial $G$-bundle to the trivial vector bundle $\RR^d\times M\to M$. It sends an automorphism $\psi:M\to G$ of the trivial bundle to the automorphism $\phi(t,x)=(\rho(\psi(x))t,x)$. 

\begin{definition}\label{smoothstr}
Fix $d\geq 0$. Let $(H_d,\rho_d)$ be a symmetry type in the sense of Definition \ref{symmtype}. We define a corresponding simplicial presheaf on $\FEmb_d$, which we denote by $\mathcal{H}_d$, whose value on a submersion $p:M\to U$ is given by the mapping space of sections
\begin{equation}\label{smhnstr}
\xymatrix{
 & \mathbf{B}H_d\ar[d]^-{\rho_d}
\\
M\ar[r]_-{\tau_p}\ar@{-->}[ru]^{\ell_p} & \mathbf{B}{\rm GL}(d),
}
\end{equation}
i.e., the mapping space \footnote{We remind the reader that ``mapping space" always means the one provided by hammock localization, unless otherwise stated. Alternatively, one could take the derived mapping space in the canonical model structure on the slice.} in the slice over $\mathbf{B}\GL(d)$:
$$\mathcal{H}_d(p:M\to U):=\map_{\mathbf{B}\GL(d)}(\tau_p,\rho_d),$$
where the slice category is equipped with weak equivalences given by maps that become weak equivalences after applying the forgetful functor.  
The map $\tau_p$ is the fiberwise tangent bundle map. We call a vertex in the simplicial set of sections a \emph{fiberwise $(H_d,\rho_d)$-structure on the submersion $p:M\to U$}.
\end{definition}

One can imagine several ``differential" enhancements of fiberwise $(H_d,\rho_d)$-structures. For example, we can add a connection to the principal $H_d$-bundle classified by the lift $\ell_p$ in \eqref{smhnstr}. In the context of Wick rotated relativistic theories, the bordisms should (at the very least) be equipped with a Riemannian metric. In Freed--Hopkins \cite{FreedHopkins}, a differential $H_d$-structure is defined as an $H_d$-structure, along with a connection on the corresponding principal $H_d$-bundle that is compatible with the Levi-Civita connection. Here we make a similar definition, adapted to suit our context.

\begin{definition}\label{riem}
Fix $d\geq 0$. We define the sheaf of \emph{fiberwise Riemannian metrics}
$$\Riem\in \Struct_d$$
as the sheaf that sends a submersion $p:M\to U$ to the set of metrics on the fiberwise tangent bundle $\tau_p$. We regard this set as a simplicial set by taking the nerve of the corresponding discrete category (with only identity morphisms). 

The structure maps of the presheaf are given by pulling back smooth families of metrics along fiberwise embeddings.
\end{definition}

Let $p:M\to U$ be a submersion with $d$-dimensional fibers. Recall that we have the notion of \emph{horizontal} and \emph{vertical} differential forms on $M$. A \emph{horizontal} differential form is a differential form that vanishes on vector fields that are tangent to the nonempty fibers of $p$. Horizontal differential forms are an ideal in the graded algebra of all differential forms on $M$. The quotient by this ideal is the graded algebra of vertical differential forms. This graded algebra admits the structure of a differentially graded algebra when equipped with the \emph{vertical} exterior derivative $d_p$. We denote by $\Omega^n_{p}(M)$ the underlying group of vertical differential $n$-forms on $M$, with respect to $p$. 

Let $G$ be a Lie group. For a submersion $p:M\to U$ with $d$-dimensional fibers, we let $\Omega^1_p(M;\mathfrak{g})$ denote the set of vertical $\mathfrak{g}$-valued 1-forms. We let $\mathbf{B}_{\nabla}G$ denote the simplicial presheaf on $\FEmb_d$ that sends a submersion $M\to U$ to the nerve of the groupoid 
$$\mathbf{B}_{\nabla}G(M\to U):=N\left(\vcenter{\xymatrix{
C^{\infty}(M,G)\times \Omega_p^1(M;\mathfrak{g})\ar@<.1cm>[r]^-{t}\ar@<-.1cm>[r]_-{s} & \Omega_p^1(M;\mathfrak{g})
}}\right).$$
The source map $s$ projects onto the second factor. On a pair $(g:M\to G,\mathcal{ A}\in \Omega^1_{U}(M;\mathfrak{g}))$, the target map $t$ is given by $t(g,\mathcal{A})={\rm Ad}_g\mathcal{A}+g^{-1}d_pg$, i.e., the action by fiberwise gauge transformations. 

Now let us recall the adjunction $F\dashv i$, where 
\begin{equation}\label{pseudostraight}
i:\PSh(\FEmb_d;\mathscr{G}{\rm pd})\to \Psi\PSh(\FEmb_d,\mathscr{G}{\rm pd})
\end{equation}
denotes the inclusion of strict functors into pseudofunctors, denoted by $\Psi\PSh$ (see Power \cite{Power}).

\begin{remark}
The left adjoint $F$ formally adds pullbacks to a pseudofunctor $X$. More precisely, if $X\in \Psi \PSh(\FEmb_d,\mathscr{G}{\rm pd})$, then $(FX)(M\to U)$ is the groupoid whose objects are pairs $((f,g),x)$, with $(f,g):(M\to U)\to (N\to V)$ a morphism in $\FEmb_d$ and $x\in X(N\to V)$. The morphisms $((f,g),x)\to ((f',g'),x')$ are morphisms $X(f,g)(x)\to  X(f',g')(x')$ in $X(M\to U)$. The unit $\eta$ of the adjunction $F\dashv i$ has $M\to U$ component given by the morphism of groupoids
$$X(M\to U)\to FX(M\to U),$$
sending an object $x$ to the pair $((\id,\id),x)$ and a morphism $x\to x'$ to itself, viewed as a morphism $X(\id,\id)(x)=x\to x'=X(\id,\id)(x')$.
\end{remark}

In the special case where $G=\GL(d)$, we define the stack $\mathbf{B}_{\nabla}\GL(d)$ differently, analogous to the way $\mathbf{B}\GL(d)$ is defined. 

\begin{definition}\label{strictconn}
Let ${\rm Vect}^{\times}_{\nabla}(d)$ be the pseudofunctor on $\FEmb_d^{\rm op}$ that sends a submersion $M\to U$ to the groupoid of rank $d$ vector bundles with connection on $M$, with isomorphisms between them.
We define the strict presheaf of groupoids 
$\mathbf{B}_{\nabla}{\rm GL}(d)=F{\rm Vect}_{\nabla}^{\times}(d).$
We have a canonical pseudonatural transformation
$$\eta:{\rm Vect}^{\times}(d)\to \mathbf{B} \GL(d),$$
given by the unit of the adjunction \eqref{pseudostraight}. In particular, given a $d$-dimensional vector bundle with connection $(V,\nabla)$, we have a corresponding object $\eta(V,\nabla)\in \mathbf{B} \GL(d)(M\to U)$, given by the pair $(\id,(V,\nabla))\in F{\rm Vect}^{\times}_{\nabla}(d):=\mathbf{B} \GL(d)$. 

We will often abuse notation and write $(V,\nabla)\in \mathbf{B}_{\nabla} \GL(d)$ for the object corresponding to a vector bundle with connection $(V,\nabla)\in {\rm Vect}^{\times}_{\nabla}(d)(M\to U)$ under the map $\eta$.
\end{definition}

Observe that the fundamental theorem of Riemannian geometry implies that we have a canonical map
$$\nabla: \Riem\to \mathbf{B}_{\nabla}\GL(d),$$
defined by sending $g$ to the corresponding fiberwise Levi-Civita connection $\nabla_g$. Given a representation $\rho_d:H_d\to {\rm O}(d)$, we also have a canonical map
$$\mathbf{B}_{\nabla}H_d\overset{\mathbf{B}(\rho_d)}{\longrightarrow} \mathbf{B}_{\nabla}{\rm O}(d)\to \mathbf{B}_{\nabla}\GL(d),$$
where the second map sends the trivial ${\rm O}(d)$-bundle with connection to the corresponding vector bundle with connection associated to the trivial principal ${\rm O}(d)$-bundles with connection.

\begin{definition}\label{diffstr}
Fix $d\geq 0$.  Let $(H_d,\rho_d)$ be a symmetry type in the sense of Definition \ref{symmtype}. We define a corresponding simplicial presheaf on $\FEmb_d$, which we denote by $\mathcal{H}^{\nabla}_d$, as the simplicial presheaf given by the homotopy pullback (in the local projective model structure):
\begin{equation}
\mathcal{H}^{\nabla}_d=\Riem\times^h_{\nabla,\mathbf{B}_{\nabla}\GL(d),\rho}\mathbf{B}_{\nabla}H_d. \label{pulldiffstr}
\end{equation}
\end{definition}

We observe that a vertex in the simplicial set $\mathcal{H}_d^{\nabla}(M\to U)$ recovers the notion of a differential $H_d$-structure in the sense of Freed--Hopkins \cite{FreedHopkins}, but adapted to suit families.

\begin{proposition}\label{diffhddata}
Fix $d\geq 0$ and let $p:M\to U\in \FEmb_d$. A map $p\to \mathcal{H}_d^{\nabla}$ determines the following data:
\begin{enumerate}
\item A fiberwise Riemannian metric $g$ on $\tau_pM$.
\item A principal $H_d$-bundle $P\to M$ with fiberwise connection $\nabla$.
\item An isomorphism of vector bundles with fiberwise connection $\phi:(P\times_{\rho_d}\RR^d,\nabla)\overset{\cong}{\to} (\tau_pM,\nabla_g)$. 
\end{enumerate} 
\end{proposition}
\proof
Let $H_d\text{-}{\Psi{\rm Bun}}_{\nabla}$ denote the pseudofunctor that sends a submersion $p:M\to U$ with $d$-dimensional fibers to the groupoid of principal $H_d$-bundles with fiberwise connection.  Let $H_d\text{-}{\rm Bun}_{\nabla}$  denote the composition of $H_d\text{-}{\Psi{\rm Bun}}_{\nabla}$ with the strictification functor $F$ in \eqref{pseudostraight} and the nerve $N:\mathscr{G}{\rm pd}\to \sset$. Using the fact that principal $H_d$-bundles are locally trivial, it is straightforward to verify that $H_d\text{-}{\rm Bun}_{\nabla}$ satisfies descent, hence is a fibrant object in the local projective model structure. Since every principal $H_d$-bundle is trivializable over an open ball centered at the origin in $\RR^d$, it follows at once that the map
$$\mathbf{B}_{\nabla}H_d\to H_d\text{-}{\rm Bun}_{\nabla}$$
is a stalkwise weak equivalence, hence a local weak equivalence by Proposition \ref{stalks}. This implies that  $H_d\text{-}{\rm Bun}_{\nabla}$ is a fibrant replacement of $\mathbf{B}_{\nabla}H_d$ in the local projective model structure. 

With this observation, we compute the homotopy fiber by the standard resolution
$$
\xymatrix{
\mathcal{H}^{\nabla}_d\ar[r]\ar[d] & H_d\text{-}{\rm Bun}_{\nabla}\times_{\mathbf{B}_{\nabla}\GL(d)}(\mathbf{B}_{\nabla}\GL(d))^{\Delta^1}\ar[d]
\\
{\rm Riem}\ar[r]^-{\nabla} & \mathbf{B}_{\nabla}\GL(d)\;.
}
$$
Mapping into the above pullback and using the universal property gives precisely the data (1)--(3).

\endofproof

Next, we define a simplicial presheaf that encodes fiberwise flat $H_d$-structure. Let $\RR^d\in \Struct_d$ be the representable presheaf associated to the projection $\RR^d\to \RR^0$. The group of smooth maps $C^{\infty}(U,H_d)$ acts on the set of fiberwise embeddings $\FEmb_d(p:M\to U, \RR^d\to \RR^0)$ as follows. Let $f:U\to H_d$ be a smooth map and let $\phi:M\to \RR^d$ be a fiberwise open embedding. Define the fiberwise embedding $f\cdot \phi$ at $x\in p^{-1}(u)$ by
\begin{equation}\label{action}
(f\cdot \phi)(x)=\rho_d(f(u))(\phi(x)).
\end{equation}
 This action commutes with pullback along fiberwise embeddings. Indeed, for a fiberwise embedding $(g:N\to M,h:V\to U)$, we have  
$$
(f\cdot \phi)(g(y))=\rho_d(f(h(v)))(\phi(g(y)))=\rho_d((f\circ h)(v))(\phi\circ g)(y)=((f\circ h)\cdot (\phi\circ g))(y).
$$
Therefore, the smooth group $H_d$ acts on $\RR^d$ in the sense of a group action internal to $\Struct_d$. 
\begin{definition}\label{flatstr}
Fix $d\geq 0$. Let $(H_d,\rho_d)$ be a symmetry type in the sense of Definition \ref{symmtype}. We define a corresponding simplicial presheaf $\mathcal{H}_d^{\rm fl}$ as the homotopy quotient 
$$\mathcal{H}_d^{\rm fl}:=\RR^d/\!/_{\rho_d} H_d,$$
 taken in the local projective model structure on $\Struct_d$. The action of $H_d$ on $\RR^d$ is defined by \eqref{action}. 
\end{definition}

\begin{remark}
The homotopy quotient $\mathcal{H}^{\rm fl}_d$ encodes fiberwise flat structures in the following sense. Let $M\to \RR^0$ be a representable with trivial base space, so that the family direction is trivial. A derived map $M\to \mathcal{H}_d$ is equivalently a choice of atlas $\{U_{a}\}_{a\in A}$ of $M$, an embedding $\phi_a:U_{a}\into \RR^d$ of each $U_a$, and an element $h_{ab}\in H_d$, for each intersection $U_a\cap U_b$, such that the diagram
$$
\xymatrix{
\RR^d\ar[rr]^{\rho_d(h_{ab})} && \RR^d
\\
&U_a\cap U_b\ar[ur]_-{\phi_{ab}}\ar[ul]^-{\phi_{ba}}&
},
$$
commutes. In particular, this implies that the fiberwise tangent bundle $\tau_p:M\to \mathbf{B}{\rm GL}(d)$ canonically lifts to a map $\ell_p:M\to \mathbf{B}H_d$, but even more is true. 

If $t:\RR^d\to \mathbf{B}H_d$ is the trivial $H_d$-structure on $\RR^d$, then the lift $\ell_p$ is defined on each $U_a$ by taking $\ell_{U_a}:=\phi_a^*(t)$. The transition maps of the $H_d$-bundle are constant on $U_a\cap U_b$, equal to $h_{ab}\in H_d$. 
This means that the $H_d$-structure determined by the lift $\ell$ is \emph{locally flat} or \emph{integrable}. Such geometries are called \emph{rigid geometries} in Stolz--Teichner \cite{StolzTeichner}. This explains our terminology \emph{fiberwise flat} $H_d$-structures, for geometric structures encoded by $\mathcal{H}^{\rm fl}_d$. 
\end{remark}

\section{Reflection positivity and the geometric bordism category}\label{bordrefl}

\subsection{The geometric bordism category}
Let $\Delta$ denote the usual simplex category and let $\Gamma$ denote Segal's gamma category. As usual, we let $\PSh_{(\infty,1)}(\cartsp\times \Gamma\times \Delta^{\times d})$ denote the relative category of simplicial presheaves on $\cartsp\times \Gamma\times \Delta^{\times d}$, with weak equivalences given by objectwise weak equivalences. In Grady--Pavlov \cite{GradyPavlov}, we defined a left Bousfield localization of the injective model structure at morphisms that implement the following conditions:
\begin{enumerate}\label{localcond}
\item Segal's special $\Delta$-condition ;
\item Completeness;
\item Globularity;
\item Segal's special $\Gamma$-condition;
\item Homotopy descent. 
\end{enumerate}
We let $S$ denote the set of maps at which we localize. We defer the reader to \cite[Section 2]{GradyPavlov} for details of this localization. 

\begin{definition}
Following \cite{GradyPavlov}, we define the relative category $\smcat_{\infty,d}$ by adjoining $S$-local weak equivalences to $\PSh_{(\infty,1)}(\cartsp\times \Gamma\times \Delta^{\times d})$. These weak equivalences implement the conditions (1) -- (5) above.
\end{definition}

We also recall from \cite{GradyPavlov} the smooth symmetric monoidal $(\infty,d)$-category of bordisms, 
$$\Bord^{\mathcal{S}}_d\in \smcat_{\infty,d}$$
where $\mathcal{S}$ is an $\infty$-sheaf on the site $\FEmb_d$. Evaluating $\Bord_d^{\mathcal{S}}$ on an object $(U,\langle \ell\rangle,{\bf m})\in \cartsp\times \Gamma\times \Delta^{\times d}$ gives a simplicial set. Ignoring the $\Gamma$-direction, which encodes symmetric monoidal structure, one should think of the resulting simplicial set as the space of $U$-parametrized families of ${\bf m}$-composable chains of $d$-morphisms. In terms of bordisms, a vertex in this simplicial set can be pictured as a $U$-parametrized grid of ``cut manifolds'' \cite{GradyPavlov} embedded fiberwise (over $U$) in an ambient manifold $M$. Composition corresponds to applying a face map to the multisimplex ${\bf m}$, which has the effect of ``removing a cut manifold''. 

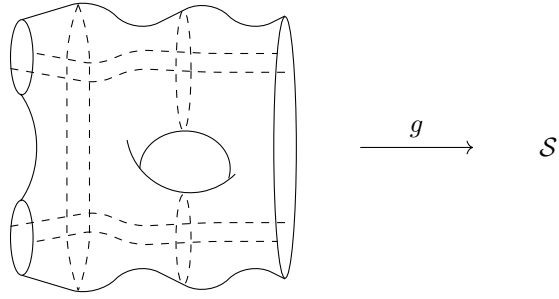
\begin{figure}[h]\label{bordism}
\begin{center}
\begin{tikzpicture}[scale=.5]
\draw[rounded corners=25pt](0,1.38)--(1,0)--(0,-1.38);
\draw (0.3,2.4) arc (0:360:0.3 and 1);
\draw  (0.3,-2.4) arc (0:360:0.3 and 1);
\draw[dashed,rounded corners=10pt] (1.5,3.8)--(1.8,3)--(1.8,-3)--(1.5,-3.8);
\draw[dashed,rounded corners=10pt] (1.5,3.8)--(1.2,3)--(1.2,-3)--(1.5,-3.8);
\draw[dashed,rounded corners=5pt] (0.4,2.5)--(2,2.2)--(3,2.6)--(5,2.5)--(6.8,2.5);
\draw[dashed,rounded corners=5pt] (-.3,2.1)--(2,1.7)--(3,2.1)--(5,2)--(7.2,2);
\draw[dashed,rounded corners=5pt] (0.4,-2.5)--(2,-2.2)--(3,-2.6)--(5,-2.5)--(6.8,-2.5);
\draw[dashed,rounded corners=5pt] (-.3,-2.1)--(2,-1.7)--(3,-2.1)--(5,-2)--(7.2,-2);
\draw[dashed] (4.5,2.05) arc (0:360:0.2 and 1.6);
\draw[dashed] (4.5,-2.45) arc (0:360:0.2 and 1.2);
\draw[rounded corners=10pt](0,3.4)--(2,4)--(3,3)--(5,4)--(6,3)--(7,3.5);
\draw[rounded corners=10pt](0,-3.4)--(2,-4)--(3,-3)--(5,-4)--(6,-3)--(7,-3.5);
\draw (7.3,0) arc (0:360:0.3 and 3.48);
\draw (2.8,0.2) arc (190:315:1.7cm and 1.7cm);
\draw (5.5,-0.82) arc (-15:180:1.2cm and 1cm);

\draw[->] (9,0) -- (12,0);

\node at (14,0) {$\mathcal{S}$};
\node at (10.5,.5) {$g$};
\end{tikzpicture}
\end{center}
\caption{A bordism (vertex) in the simplicial set $\Bord^{\mathcal{S}}_d(\RR^0,\langle 1\rangle,([1],[3]))$. The core is embedded in an ambient manifold $M$, which extends outside the core. A geometric structure is given by  a morphism $g:(M\to \RR^0)\to \mathcal{S}$ in $\sh_{(\infty,1)}(\FEmb_d)$ around the core.}
\end{figure}

Figure \ref{bordism} depicts a vertex in the simplicial set $\Bord^{\mathcal{S}}_d(\RR^0,\langle 1\rangle,([1],[3]))$. An $l$-simplex is given by a $\Delta^l$-parametrized family of cuts on the manifold along with the germ of a morphism $g:M\to \mathcal{S}$  around the union of cores (indexed by $\Delta^l$). In general, for arbitrary $U\in \cartsp$, an $l$-simplex is a $U\times \Delta^l$-parametrized family of cuts along with the germ of a morphism $g:M\times U\to \mathcal{S}$. 

In Grady--Pavlov \cite{GradyPavlov2}, we proved a powerful result about field theories that arise as symmetric monoidal functors on $\mathfrak{Bord}_d^{\mathcal{S}}$. Namely, such field theories can be classified using a \emph{geometric cobordism hypothesis}, similar to the way that topological field theories are classified by Baez--Dolan--Lurie's cobordism hypothesis \cite{Lurie.cob}. The following example is Freed--Hopkins \cite[Example 3.3]{FreedHopkins}, adapted to our setting. In this example, we start with a Lagrangian and derive a functorial QFT using the geometric cobordism hypothesis established in \cite{GradyPavlov2}. We observe that the reflection structure is seen even at the prequantum level.

\begin{example}\label{particleoncircle}
Let $d=1$. Let us construct a prequantum functorial field theory which describes a free particle on a circle. We regard $S^1$ as a geometric structure for field theories, as in Example \ref{targetman}. Let ${\rm Riem}_1$ be the sheaf of fiberwise Riemannian metrics (Definition \ref{riem}). The relevant geometric structure for the theory is the product ${\rm Riem}_1\times S^1$. 

We will interpret the Lagrangian as a natural transformation 
$$
L:{\rm Riem}_1\times S^1\to \Omega ^1
$$
as follows. For a submersion $p:M\to U$ with $1$-dimensional fibers, a morphism $p\to {\rm Riem}\times S^1$ is equivalently a choice of fiberwise Riemannian metric $g$ and a smooth map $f:M\to S^1$. The metric $g$ induces a unit speed parametrization $\gamma_g:\RR\to M$. Let $x_{(g,f)}=f\circ \gamma_g$. The Lagrangian sends a pair $(g,f)$ to the differential 1-form
$$L=\tfrac{1}{2}\dot{x}_{(g,f)}^2 dt,$$
where $dt$ is the standard volume form on $\RR$. 

 We have a further map $\Omega^1 \to {\rm Vect}_{\nabla}$ that sends a 1-form to a connection on the trivial line bundle, $\mathcal{A}\mapsto d+\mathcal{A}$. In total, we have produced a morphism   
\begin{equation}\label{maplagrange}
L:{\rm Riem}_1\times S^1\to {\rm Vect}^{\times}_{\nabla}.
\end{equation}

Under the geometric cobordism hypothesis of Grady--Pavlov \cite{GradyPavlov2}, this map corresponds to a functorial field theory, which we now describe. The field theory takes values in ${\rm Vect}_{\CC}$. 
Since every connected open 1-dimensional Riemannian manifold is isometric to an open interval in $\RR$, we can focus attention on open intervals in $\RR$. An object in the bordism category is given by an interval with a distinguished point, which we can take to be of the form $(t-\epsilon,t+\epsilon)$ \footnote{Embeddings of open intervals that preserve the distinguished point are isomorphisms in the bordism category.} along with a smooth map $x:(t-\epsilon,t+\epsilon)\to S^1$. To this data we assign the complex vector space $\CC$, regarded as the fiber of the trivial bundle at $y=x(t)$. 

 A morphism in the bordism category is a path $x:(t_0-\epsilon,t_1+\epsilon)\to S^1$ along with a Riemannian metric, which we can take to be $edt^2$ for some smooth function $e$ taking values in the positive reals. Let us denote the action corresponding to the Lagrangian density by
$$S(x)=\int_{t_0}^{t_1}\tfrac{1}{2e}\dot{x}^2dt.$$ 
To such data, we associate the linear map given by multiplication
\begin{equation}\label{transmap}
e^{-S(x)}:\CC\to \CC.
\end{equation}
This data uniquely determines the remainder of the field theory. On the other hand, the above data also uniquely determines a connection on the trivial bundle $\CC\times S^1\to S^1$. The parallel transport map associated to the connection on the trivial bundle picked out by \eqref{maplagrange} is precisely \eqref{transmap}. 

To see the time reversal symmetry, it is useful to add a $\theta$-angle to this field theory. Fix $\theta\in \RR$ and a normalized volume form $\omega\in \Omega^1(S^1)$, i.e., $\int_{S^1}\omega=1$. Now define 
$$L_{\theta}:{\rm Riem}_1\times S^1\to \Omega^1,$$
by assigning a pair $(edt^2,x)$ to 
$$L_{\theta}=\tfrac{1}{2e}\dot{x}^2 dt-i\theta x^*(\omega).$$
The corresponding transport map on the field theory side is 
$$e^{-S_{\theta}(x)}=e^{-S(x)}e^{-i\theta \int_{t_0}^{t_1}x^*(\omega)}:\CC\to \CC.$$
For example, if $x:\RR \to S^1$ is given by $x(t)=e^{it}$ and $\omega=\frac{1}{2\pi i}d \log$, then $x^*(\omega)= dt$. The bordism given by an $\epsilon$-neighborhood of the interval $[0,2\pi]$ is sent to the linear map given by multiplication by $e^{-S(x)}e^{-i\theta}$. 

Now time reversal corresponds to change of orientation of $S^1$. An orientation reversing isometry $\phi$ induces a natural transformation of functors $\phi:L_{\theta}\to L_{-\theta}$. On the field theory side, $\phi$ induces a transformation of field theories that sends the field theory with transport operator $S_{\theta}$ to the field theory with transport operator $S_{-\theta}$. Again letting $\omega=\frac{1}{2\pi i}d\log$ and $x(t)=e^{it}$, we see that the field theory sends the bordism $((-\epsilon,2\pi+\epsilon),x,dt^2)$ to the linear map
$$e^{-(S_{-\theta}(x))}z=e^{-S}e^{i\theta}z=e^{-S(x)}e^{-i\theta}\cdot z=\overline{e^{-S_{\theta}(x)}}z,$$
where $\cdot$ denotes scalar multiplication in the conjugate vector space $\overline{\CC}$. Hence, changing the sign of $\theta$ corresponds to conjugation on the field theory side.
\end{example}

\subsection{Reflection structure}

We now want to define an involution on $\mathfrak{Bord}_d^{\mathcal{H}_d}$ analogous to Freed--Hopkins \cite[Section 4]{FreedHopkins}. There, extended symmetry groups $\widehat{H}_d$ of $H_d$ were considered, given by adjoining an involution. Concretely, the group $\widehat{H}_d$ is given by an extension \begin{equation}\label{jmap}
1\to H_d\overset{j_d}{\to} \widehat{H}_d\to \ZZ/2\to 1\;,
\end{equation}
which is (noncanonically) split by a choice of hyperplane reflection $\sigma\in {\rm O}(d)$ \cite[Proposition 3.13]{FreedHopkins}. The group $\widehat{H}_d$ fits into a pullback diagram 
$$
\xymatrix{
H_d\ar[r]\ar[d] & \widehat{H}_d\ar[d]
\\
{\rm O}(d)\ar[r] & {\rm O}(d)\times \ZZ/2\;.
}
$$
A canonical example is the following. 
\begin{example}
Let $(H_d,\rho_d)$ be a symmetry type. Suppose that the image of $\rho_d:H_d\to {\rm O}(d)$ is precisely ${\rm SO}(d)$. Let $K=\ker(\rho_d)$. Since ${\rm Pin}_+(d)$ is a double cover of ${\rm O}(d)$, it follows from \cite[Theorem 2.7]{FreedHopkins} that there is a central element $k_0\in K$ such that 
$$H_d\cong {\rm Spin}(d)\times K/\langle (-1,k_0)\rangle.$$
Set
$$\widehat{H}_d={\rm Pin}_+(d)\times K/\langle (-1,k_0)\rangle.$$
Then projection onto $\pi_0{\rm Pin}_+(d)$ gives rise to an exact sequence 
$$1\to {\rm Spin}_d\times K/\langle (-1,k_0)\rangle\to \widehat{H}_d\to \ZZ/2\to 1$$
\end{example}

Given a smooth $(H_d,\rho_d)$-structure, we would like to define the opposite smooth structure $(H_d',\rho_d')$ analogous to the construction of \cite[Definition 4.1]{FreedHopkins}. To this end, let us observe that the map $j_d:H_d\to \widehat{H}_d$ in \eqref{jmap} induces a map $j_d:\mathbf{B}H_d\to \mathbf{B}\widehat{H}_d$ that fits into the diagram 
$$
\xymatrix{
\mathbf{B}H_d\ar[r]^{j_d}\ar[d] & \mathbf{B}\widehat{H}_d\ar[d]
\\
\mathbf{B}{\rm O}(d)\ar[r]^-{i_+} & \mathbf{B}{\rm O}(d)\times \ZZ/2\;,
}
$$
where $i_+$ includes into the connected component indexed by $1\in \ZZ/2$. Now a splitting $\sigma$ of \eqref{jmap} gives rise to an automorphism $\phi_{\sigma}$ of $H_d$. At the level of classifying stacks, we have a commutative diagram 
$$
\xymatrix{
\mathbf{B}H_d\ar[r]^{\overline{\phi}_{\sigma}}\ar[d]^-{\rho_d} & \mathbf{B}H_d\ar[d]^-{\rho_d}
\\
\mathbf{B}{\rm O}(d)\ar[r]^-{\sigma} & \mathbf{B}{\rm O}(d)
},
$$
where $\overline{\phi}_\sigma=\mathbf{B}(\phi_\sigma)$. Post-composition by the automorphism $\sigma:\mathbf{B}\GL(d)\to \mathbf{B}\GL(d)$, given by the hyperplane reflection $\sigma$, fixes the fiberwise tangent bundle map $\tau_p:M\to \mathbf{B}\GL(d)$, but acts non trivially on automorphisms of $\tau_p$. We  have an induced morphism on the simplicial set of sections 
\begin{equation}\label{invbord}
\phi_{\sigma}:\mathcal{H}_d(M\to U)=\map_{\mathbf{B}\GL(d)}(\tau_p,\rho_d)\to \map_{\mathbf{B}\GL(d)}(\tau_p,\rho_d)=\mathcal{H}_d(M\to U).
\end{equation}

\begin{proposition}\label{betahd1}
The morphism \eqref{invbord} is involutive, for all submersions $M\to U$. It is also natural with respect to fiberwise open embeddings and therefore gives rise to an involution in $\sh_{(\infty,1)}(\FEmb_d)$:
$$\beta_{\sigma}:\mathcal{H}_d\to \mathcal{H}_d \qquad \beta_{\sigma}^2=1.$$
\end{proposition}
\proof
Since $\sigma^2=1$, it follows also that the automorphism $\phi_{\sigma}:H_d\to H_d$ is an involution. Since \eqref{invbord} is defined by composition with $\overline{\phi}_\sigma=\mathbf{B}(\phi_{\sigma})$, the claim follows.
\endofproof

We extend the involution $\beta_{\sigma}$ to \emph{differential} $(H_d,\rho_d)$-structures as follows. Let $\phi_{\sigma}$ be the automorphism of $H_d$ induced by a splitting of \eqref{jmap}. From \cite[Proposition 3.13]{FreedHopkins}, we recall that this automorphism $\phi_{\sigma}:H_d\to H_d$ is given by (twisted) conjugation \footnote{In the case where the image of $\rho_d$ is ${\rm SO}(d)$, the automorphism is conjugation by the unique lift $\tilde \sigma\in {\rm Pin}^+(d)$ of $\sigma$. When the image is all of ${\rm O}(d)$, the automorphism is twisted conjugation by this element, where the twist is multiplication by the nontrivial character ${\rm Pin}^+(d)\to \ZZ/2$.} by a lift $\tilde{\sigma}$ of $\sigma$. 

\begin{proposition}\label{invbord2}
We have an induced automorphism of differential 1-forms $\Omega^1(-;\mathfrak{h}_d)$, via 
\begin{equation}\label{mapgauge}
d\phi_{\sigma}\vert_{1}:\mathcal{A}\mapsto {\rm Ad}_{\tilde{\sigma}}(\mathcal{A}).
\end{equation}
The following hold.
\begin{enumerate}
\item The map \eqref{mapgauge} commutes with gauge transformations and induces a morphism of stacks
$$\beta_{\sigma}:\mathbf{B}_{\nabla}H_d\to \mathbf{B}_{\nabla}H_d.$$
\item conjugation by $\sigma$ induces a map on $\Omega^1(-;\mathfrak{o}(d))$, analogous to \eqref{mapgauge}, and induces a morphism of stacks 
$$\bar\beta_{\sigma}:\mathbf{B}_{\nabla}{\rm O}(d)\to \mathbf{B}_{\nabla}{\rm O}(d).$$\
\item We have a commutative diagram 
$$
\xymatrix{
\mathbf{B}_{\nabla}H_d\ar[r]^-{\beta_{\sigma}} \ar[d] & \mathbf{B}_{\nabla}H_d\ar[d]
\\
\mathbf{B}_{\nabla}{\rm O}(d)\ar[r]^{\bar\beta_{\sigma}} & \mathbf{B}_{\nabla}{\rm O}(d)
}
$$
\item We have an induced involution 
$$\beta_{\sigma}:\mathcal{H}_d^{\nabla}\to \mathcal{H}_d^{\nabla} \qquad \beta_{\sigma}^2=1.$$
\end{enumerate}
\end{proposition}
\proof
To see that ${\rm Ad}_{\sigma}$ commutes with gauge transformations, let $\mathcal{A}:(M\to U)\to \Omega^1(-;\mathfrak{h}_d)$ be a fiberwise differential 1-form on $M$ and let $g:M\to H_d$ be a smooth function. Observe that for all smooth maps $h:M\to H_d$, we have
$$
\phi_{\sigma}(g)\tilde\sigma h\tilde \sigma^{-1} \phi_{\sigma}(g)^{-1}= \tilde \sigma g \tilde\sigma^{-1}\tilde \sigma h\tilde \sigma^{-1}\tilde \sigma g\tilde \sigma^{-1}=\tilde \sigma ghg^{-1}\tilde \sigma^{-1}.
$$
Hence, ${\rm Ad}_{\phi_{\sigma}(g)}{\rm Ad}_{\tilde \sigma}={\rm Ad}_{\tilde \sigma}{\rm Ad}_{g}$. Moreover, $d(\phi_{\sigma}\circ g)={\rm Ad}_{\tilde \sigma}(dg)$ and $dL_{\tilde \sigma g\tilde \sigma^{-1}}{\rm Ad}_{\tilde \sigma}={\rm Ad}_{\tilde \sigma}dL_{g}$. So
\begin{eqnarray*}
{\rm Ad}_{\phi_{\sigma}(g)}({\rm Ad}_{\tilde \sigma}(\mathcal{A}))+dL_{\phi_{\sigma}(g)^{-1}}d\phi_{\sigma}(g) 
&=& {\rm Ad}_{\tilde \sigma}{\rm Ad}_g(\mathcal{A})+dL_{\phi_{\sigma}(g)^{-1}}{\rm Ad}_{\tilde \sigma}(dg)
\\
&=& {\rm Ad}_{\tilde \sigma}\left({\rm Ad}_g(\mathcal{A})+dL_{g^{-1}}dg\right)
\end{eqnarray*}
and ${\rm Ad}_{\tilde \sigma}$ indeed commutes with gauge transformations. We immediately deduce that ${\rm Ad}_{\tilde \sigma}$ induces a morphism on action groupoids $\Omega^1(-;\mathfrak{h}_d)/\!/H_d$ and hence gives rise to a well defined morphism of smooth stacks 
$$\beta_{\sigma}:\mathbf{B}_{\nabla}H_d\to \mathbf{B}_{\nabla}H_d.$$
The proof of (2) is completely analogous, with ${\rm Ad}_{\sigma}$ replacing ${\rm Ad}_{\tilde \sigma}$. Part (3) follows at once from the commutativity 
\begin{equation}\label{comprefl}
\xymatrix{
H_d\ar[rr]^{\phi_{\sigma}}\ar[d]^-{\rho_d} && H_d\ar[d]^-{\rho_d}
\\
{\rm O}(d)\ar[rr]^-{\sigma (-)\sigma^{-1}} && {\rm O}(d)\;,
}
\end{equation}
which implies commutativity $d\rho_d{\rm Ad}_{\tilde \sigma}={\rm Ad}_{\sigma}d\rho_d$ at the level of Lie algebras, and hence commutativity at the level of classifying stacks. Finally, for part (4), the maps $\bar{\beta}_{\sigma}$ and $\beta_{\sigma}$ from parts (1) and (2) induce an automorphism of the homotopy pullback  
$$
\beta_{\sigma}={\rm id}\times_{\bar\beta_{\sigma}}\beta_{\sigma}:\mathcal{H}_d^{\nabla}={\rm Riem}\times^h_{\mathbf{B}_{\nabla}{\rm O}(d)}\mathbf{B}_{\nabla}H_d\to {\rm Riem}\times^h_{\mathbf{B}_{\nabla}{\rm O}(d)}\mathbf{B}_{\nabla}H_d=\mathcal{H}_d^{\nabla}.
$$ 
This automorphism is involutive, since $\beta_{\sigma}$ and $\bar \beta_{\sigma}$ are involutive. 
\endofproof

For fiberwise flat structures, we define the involution as follows. By the commutativity of \eqref{comprefl}, the hyperplane reflection $\sigma:\RR^d\to \RR^d$ commutes with $\phi_{\sigma}$. That is, for all $h:U\to H_d$, we have
$$\sigma(\rho_d(h))=\rho_d(\phi_{\sigma}(h))\sigma.$$
Therefore, we have an induced involution on action groupoids
\begin{equation}\label{flatinv}
(\sigma,\phi_d):\RR^d/\!/H_d\to \RR^d/\!/H_d.
\end{equation}

\begin{definition}\label{betahd2}
We denote the involution \eqref{flatinv} on $\mathcal{H}_d^{\rm fl}$ again by 
$$
\beta_{\sigma}:\mathcal{H}_d^{\rm fl}\to \mathcal{H}_d^{\rm fl} \qquad \beta_{\sigma}^2=1.
$$
\end{definition}

\section{Deformations of structures}\label{structdeform}

In this section, we introduce several different ``deformation space constructions'', which will be used at various points throughout the paper. We also prove some basic results related to the construction. In particular, a convenient criteria for showing that a morphism of geometric structures induces a weak equivalence on deformation spaces is established in Theorem \ref{critwe}. 

Throughout the remainder of the section, we let ${\sf V}$ denote a model category. The main examples of ${\sf V}$ that we will consider will be the category of spaces $\mathscr{S}$, symmetric spectra, or a left Bousfield localization of simplicial presheaves. In the case where ${\sf V}$ is a symmetric monoidal model category, we let $F(-,-)$ denote the internal hom in ${\sf V}$. For example, when ${\sf V}=\Sp$, $F(-,-)$ denotes the function spectrum. In general, we let $\map(-,-)$ denote the simplicial mapping space, given by hammock localization. In the case where ${\sf V}$ is simplicial, then we can alternatively use the derived mapping simplicial set, which is weakly equivalent to the hammock localization.

Let $\PSh_{(\infty,1)}(\cartsp;{\sf V})$ be the relative category whose underlying ordinary category is the category of presheaves with values in ${\sf V}$, and whose weak equivalences are those natural transformations that are objectwise weak equivalences. 
In the case where ${\sf V}$ is symmetric monoidal and simplicial, we can canonically enrich $\cartsp$ over ${\sf V}$ by setting $\cartsp_{\sf V}(U,V)=C^{\infty}(U,V)\otimes \mathbf{1}$, where ${\bf 1}$ is the monoidal unit. Moreover, the category of simplicial presheaves $\PSh_{(\infty,1)}(\cartsp;{\sf V})$ can be canonically enriched in ${\sf V}$. For $X,Y\in \PSh_{(\infty,1)}(\cartsp;{\sf V})$, we define the enrichment via the end
\begin{equation}\label{enrichcomb}
F(X,Y):=\int_{U\in \cartsp}F(X(U),Y(U)).
\end{equation}

\subsection{The smooth deformation space construction} Much of the material in this subsection has already been established in Berwick-Evans -- Boavida de Brito -- Pavlov \cite{Pavlov} and Pavlov \cite{Pavlov2}. Here, we adapt these results to our setting.

 Let $t:\cartsp\to \ast$ be the canonical functor. By left and right Kan extension, we have a triple adjunction 
\begin{equation}\label{prestr}
\xymatrix{
\PSh_{(\infty,1)}(\cartsp;{\sf V})\ar@<.4cm>[rr]^-{t_!}\ar@<.25cm>@{}[r]|<<<<<<<<<<<<<<<<<<{\downv} \ar@<-.4cm>[rr]_-{t_*} \ar@<-.25cm>@{}[r]|<<<<<<<<<<<<<<<<<<{\downv} && \ar@<0cm>[ll]|-{t^*}{\sf V}\;.
}
\end{equation}
If the objects of ${\sf V}$ are called \emph{gadgets}, then we call the objects of $\PSh_{(\infty,1)}(\cartsp;{\sf V})$ \emph{smooth gadgets}. Using the above notation, if $X$ is a smooth gadget, we can define the \emph{pre deformation gadget of $X$} as the left Kan extension $t_!(X)$, defined in \ref{prestr}. The right adjoint functor $t_*$ is given by evaluating at the point $\RR^0$.

\def\smV{{\rm C}^{\infty}{\sf V}}
\begin{notation}
Let ${\sf V}$ be a model category. We denote the relative category of ${\sf V}$-valued $\infty$-sheaves (see Notation \ref{notation}) by
$$\smV:=\sh_{(\infty,1)}(\cartsp;{\sf V})$$
If the objects of ${\sf V}$ are called gadgets, we call the objects of $\smV$ \emph{smooth gadgets}.
\end{notation}

\begin{definition}\label{constglob}
We define the \emph{constant presheaf functor} as the functor $t^*$ in \eqref{prestr}. Henceforth, we denote this functor by $\delta^{pre}_{\sf V}:=t^*$. We define the \emph{global sections functor} as the functor $t_*$. Henceforth, we denote this functor by $\gamma_{\sf V}:=t_*$
\end{definition}

Note that the right adjoint $\delta^{pre}_{\sf V}=t^*$ preserves all weak equivalences, since these are objectwise. This is also true at the level of sheaves, however $\delta^{pre}_{\sf V}$ does not send fibrant objects to fibrant objects (in either the local projective or injective model structure). Hence, in the case of sheaves, $\delta^{pre}_{\sf V}$ must be appropriately derived. The following functor will give us a convenient presentation for the corresponding derived functor. We begin with a few observations.

In the following, we will restrict to the case when ${\sf V}$ is a \emph{simplicial} model category. Although this restriction is not really necessary, it allows us to avoid some pedantry that arises when defining certain models for derived functors. Furthermore, we lose no generality in restricting to simplicial case, since the relative category of all combinatorial model categories and left Quillen functors is Dwyer-Kan equivalent to the relative category of all simplicial left proper combinatorial model categories and simplicial left Quillen functors (Pavlov \cite[Theorem 1.2]{Pavlov.C}).

 
\begin{definition}\label{constsheaf}
Let ${\sf V}$ be a simplicial combinatorial model category. As described in the preceding paragraph, ${\sf V}$ is tensored and powered over $\mathscr{S}$. Let ${\rm sing}:\mathscr{T}{\rm op}\to \mathscr{S}$ denote the singular simplicial set functor. We define the \emph{constant sheaf functor} 
$$
\delta_{\sf V}:{\sf V}\to \smV
$$
by the formula 
$$\delta_{\sf V}(v)(U)=v^{{\rm sing}(U)}\;,$$
for all $v\in {\sf V}$, $U\in \cartsp$. Here the superscript on the right denotes the powering by the singular simplicial set of $U$. 
\end{definition}

Observe that there is a canonical comparison map $\delta^{pre}_{\sf V}\to \delta_{\sf V}$, defined as follows. For all $v\in {\sf V}$ and $U\in \cartsp$, the terminal map ${\rm sing}(U)\to \ast$ induces a map $\delta^{pre}_{\sf V}(v)(U)=v\to v^{{\rm sing}(U)}=\delta_{\sf V}(v)(U)$, which is natural in both $v$ and $U$.
\begin{lemma}\label{contractopens}
For all $v\in {\sf V}$, the canonical map $\delta^{pre}_{\sf V}(v)\to \delta_{\sf V}(v)$ evaluates to a weak equivalence on all contractible $U\in \cartsp$.
\end{lemma}
\proof
A choice of point $u\in U$ gives a canonical map $u:\ast\to {\rm sing}(U)$, which is an acyclic cofibration in the standard model structure on $\mathscr{S}$, since $U$ is contractible. By the pullback power axiom, the induced map $v^{{\rm sing}(U)}\to v$ is an acyclic fibration in ${\sf V}$. Hence $\delta_{\sf V}(v)(U)\simeq v$. On the other hand, the constant presheaf functor satisfies $\delta^{pre}_{\sf V}(v)(U)=v$ by definition. 
\endofproof

Now turning our attention to the left adjoint $t_!$, observe that this functor does not preserve weak equivalences even at the level of presheaves. Deriving the functor $t_!$ amounts to passing from pre deformation gadgets to \emph{deformation gadgets}. 

\begin{definition}\label{extsimp}
We define the \emph{extended $n$-simplex} as the smooth manifold
$$\Delta^n_e=\{(t_0,\hdots,t_n)\in \RR^{n+1}:\sum t_i=1\}\subset \RR^{n+1}.$$
We define the \emph{boundary of the extended $n$-simplex} $\partial \Delta^n_e$ as the (non-representable) presheaf of sets on $\cartsp$ given by the coequalizer (taken in presheaves on $\cartsp$):
\begin{equation}\label{boundarynsimp}
\xymatrix{
\coprod_{0\leq i<j\leq n}\Delta^{n-2}_e\ar@<.1cm>[r]\ar@<-.1cm>[r] & \coprod_{0\leq i \leq n}\Delta^{n-1}_e\to \partial \Delta^n_e\;,
}
\end{equation}
where on the coproduct summand indexed by $i<j$, the two horizontal maps are the two coface maps $d_i$ and $d_{j-1}$. There is a canonical inclusion $\partial \Delta^n_e\into \Delta^n_e$ induced by the universal property of the coequalizer.
\end{definition}

The extended $n$-simplex gives us a natural way to define smooth families of geometric objects, parametrized by simplices. In fact, one can naturally form an object in ${\sf V}$ from a smooth $(\infty,1)$-sheaf with values in ${\sf V}$, using extended $n$-simplices. When ${\sf V}=\mathscr{S}$, this construction can be thought of as an operation that takes in a smooth $(\infty,1)$-sheaf and outputs its smooth deformation space. The following is due to Berwick-Evans -- Boavida de Brito -- Pavlov \cite{Pavlov}. It provides a very convenient presentation for the derived functor of $t_!$.

\begin{definition}\label{struts}
Let $X\in \smV$ be a smooth gadget. Then we define the \emph{deformation gadget of $X$} as the homotopy colimit
$$
\csp_{\sf V}(X):=\hocolim_{[n]\in \Delta^{\rm op}}X(\Delta^n_e),
$$
where $\Delta^n_e$ is the extended $n$-simplex. The face and degeneracy maps in the above diagram are induced by functoriality of $X$ on $\cartsp$. This defines a functor:
$$\csp_{\sf V}:\smV\to {\sf V}.$$
We also define two related endofunctors (for general ${\sf V}$, the second lands only in presheaves) as follows
\begin{itemize}
\item  We define the \emph{internal deformation gadget functor} 
$$
\Icsp_{\sf V}:\smV\to \smV,
$$
by setting $\Icsp_{\sf V}:=\delta_{\sf V}\csp_{\sf V}$. 

\item We define the \emph{concordification functor}
$$
\Con_{\sf V}:\smV\to \PSh_{(\infty,1)}(\cartsp;{\sf V}),
$$
by the formula:
$$\Con_{\sf V}(X)(U):=\hocolim_{[n]\in \Delta^{\rm op}}X(\Delta^n_e\times U).$$
\end{itemize}
 When ${\sf V}$ is a left proper, simplicial model variety \cite[Definition 13.3]{Pavlov2}, the functor $\Con_{\sf V}$ sends {\v C}ech local objects to {\v C}ech local objects, so that $\Con_{\sf V}$ factors through sheaves and is equivalent to $\Icsp_{\sf V}$ \cite[Theorem 13.8]{Pavlov2}. 

\end{definition}

\begin{proposition}\label{combhleft}
Suppose ${\sf V}$ is combinatorial. Then the functor $\csp_{\sf V}$ in Definition \ref{struts} is homotopy left adjoint to the constant sheaf functor $\delta_{\sf V}$ in Definition \ref{constsheaf}. 
\end{proposition}
\proof
For ${\sf V}=\mathscr{S}$, this is \cite[Proposition 1.3]{Pavlov}. By Dugger's theorem \cite[Theorem 1.1]{Dugger}, it suffices to treat the case where ${\sf V}$ is a left Bousfield localization of simplicial presheaves on some small category $K$. In this case ${\sf V}$ is simplicial, hence also $\PSh_{(\infty,1)}(\cartsp;{\sf V})$ is a module over $\PSh_{(\infty,1)}(\cartsp)$, where both model structures are the injective model structure and the tensoring is defined pointwise. 

Now $\PSh_{(\infty,1)}(\cartsp;{\sf V})$ is generated under homotopy colimits by objects of the form $U\otimes k$, where $U\in \cartsp$ and $k\in K$, viewed as representable objects in $\PSh_{(\infty,1)}(\cartsp)$ and $\PSh_{(\infty,1)}(K)$, respectively. Moreover, since every $U\in \cartsp$ admits a good covering by contractible open subsets, the corresponding category of sheaves $\sh_{(\infty,1)}(\cartsp;{\sf V})$ is generated under homotopy colimits by objects of the form $U\otimes k$, where $U$ is contractible. Since $\csp$ is manifestly homotopy cocontinuous, it suffices to prove the claim on objects of the form $U\otimes k$, with $U$ contractible. We have 
$$
\csp(U\otimes k)=\hocolim_{[n]\in \Delta^{\rm op}}\left(C^{\infty}(\Delta^n_e,U)\otimes k\right)\simeq \left(\hocolim_{[n]\in \Delta^{\rm op}}C^{\infty}(\Delta^n_e,U)\right)\otimes k\simeq k\;,
$$
where the last equivalence follows form the well known fact that the inclusion of the smooth singular nerve into the continuous singular nerve is a weak equivalence. It follows that for all $v\in {\sf V}$, we have an isomorphism in $\Ho \mathscr{S}$:
$$\map(\csp(U\otimes k),v)\simeq v(k)\;.$$
On the other hand, by Lemma \ref{contractopens} and Yoneda, we have an isomorphism in $\Ho \mathscr{S}$: 
$$\map(U\otimes k, \delta(v))\simeq \map_{\sf V}(k,\delta(v)(U))\simeq \map_{\sf V}(k,v^{\sing(U)})\simeq v(k).$$
Hence we have verified the homotopy adjunction on objects of the form $U\otimes k$, completing the proof.

\endofproof


\def\fraksmcat{\mathfrak{Cat}_{\infty,d}^{\otimes}}

\begin{proposition}\label{descent}
For ${\sf V}=\mathscr{S}$, ${\sf V}=\mathscr{S}^{{\rm O}(d)}$, and ${\sf V}=\mathscr{P}{\rm Sh}_{(\infty,1)}(\Delta^{\times d}\times \Gamma)$,  the functors $\Con_{\sf V}$ sends homotopy sheaves to homotopy sheaves and is equivalent to $\Icsp_{\sf V}=\delta_{\sf V}\csp_{\sf V}$. For ${\sf V}=\Sp$ and an $(\infty,1)$-sheaf $X\in \smV$, the presheaf $\Con_{\Sp}(X)$ satisfies homotopy descent along finite covers. 
\end{proposition}
\proof
The first claim follows from the main result of \cite{Pavlov}. In fact, $\Con_{\mathscr{S}}$ even satisfies descent on all smooth manifolds. The second and third claim follows from the first, since in both cases, $\mathcal{C}_{\sf V}$ is simply given by the prolongation of $\mathcal{C}_{\sf V}$ to presheaves and hence factors through sheaves precisely when $\mathcal{C}_{\mathscr{S}}$ does. The third claim is \cite[Proposition 7.6]{BNV}.
\endofproof

\subsection{The fiberwise deformation space construction}

Our next goal is to recall two other homotopy cocontinuous functors that were introduced in Grady--Pavlov \cite{GradyPavlov} \footnote{The composite of the two functors was called $\mathcal{C}_d$ in \cite{GradyPavlov}}. The first functor, which we denote as $\R$, takes a presheaf on the \emph{enriched} site of fiberwise open embeddings (Definition \ref{enrichedfemb}) and evaluates on representables of the form $p_U:\RR^d\times U\to U$. The resulting object is a presheaf on $\cartsp$ with values in equivariant spaces $\mathscr{S}^{{\rm O}(d)}$, where the action of ${\rm O}(d)$ is provided by the action of ${\rm O}(d)$ on $\RR^d$, by rotations. The second functor $\E$ promotes a presheaf on the non-enriched site to the enriched site by formally adding isotopies to geometric structure. 

\begin{figure}[h] 
\begin{center}
\begin{tikzpicture}
\node at (0,0) {$\left\{\parbox{3.15cm}{\itshape geometric structures}\right\}$};
\node at (2.5,2) {$\left\{\parbox{3.15cm}{\itshape geometric structures with isotopies}\right\}$};
\node at (8,2) {$\left\{\parbox{2.9cm}{\itshape smooth equivariant spaces}\right\}$};
\node at (10,0) {$\left\{\parbox{2.8cm}{\itshape equivariant spaces}\right\}$};

\draw[->] (-.3,.3) -- (1.3,1.5);
\draw[->] (1.8,1.5) -- (.2,.3);

\draw[->] (4.5,2.2) -- (6.1,2.2);
\draw[->] (6.1,1.9) -- (4.5,1.9);

\draw[->] (8,1.5) -- (9.5,.3);
\draw[->] (10,.3) -- (8.5,1.5);

\draw[->] (2,.15) -- (8,.15);
\draw[->] (8,-.15) -- (2,-.15);

\node at (5.5,3) {$\parbox{2cm}{scannning adjunction}$};
\node at (-1,1.3) {$\parbox{2cm}{isotopification adjunction}$};
\node at (11,1.3) {$\parbox{2cm}{deformation adjunction}$};
\node at (5.5,-1) {$\parbox{4cm}{deformation of structures adjunction}$};

\end{tikzpicture}
\end{center}
\caption{Adjunctions between various $(\infty,1)$-categories of sheaves.}\label{adjuncts}
\end{figure}
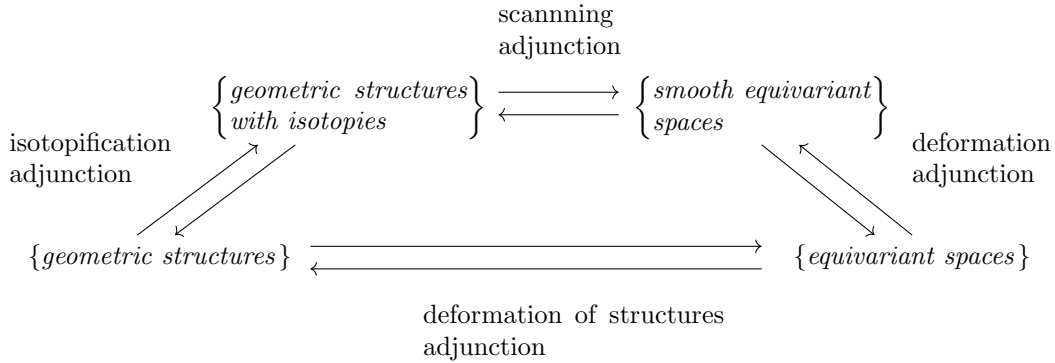

Since there are several different categories of $(\infty,1)$-sheaves and adjunctions which are introduced in this sections, it may be helpful to have a schematic picture of how these $(\infty,1)$-sheaves are related. We depict these relationships in Figure \ref{adjuncts}. In the figure, adjunctions are labeled by a phrase which we feel best describes the essence of the adjunction. For example, the left adjoint in the \emph{isotopification adjunction} formally adds isotopies to geometric structures. The adjunction labeled \emph{scanning adjunction} is very interesting. The left adjoint is the restriction functor $\R_d$ described in the preceeding paragraph, while the right adjoint forms a bundle associated to the fiberwise frame bundle of a submersion $p:M\to U$. The unit of the adjunction is a fiberwise version of Segal's scanning map. Composing the top three adjunctions in the figure together yields the bottom adjunction. Thus, the unit of the bottom adjunction can be thought of as a variant of the Segal's scanning map. Morally, taking deformations of geometric structures on bordisms amounts passing to a corresponding topological structure via the scanning map. In the next section we relate deformations of the bordism category with some fixed geometric structure to bordisms equipped with the \emph{topological} structure given by deformations of that geometric structure.

In the previous section, we defined geometric structures, via the $(\infty,1)$-category of sheaves $\sh_{(\infty,1)}(\FEmb_d)$. We also implicitly defined smooth equivariant spaces, via the presheaf category $\sh_{(\infty,1)}(\cartsp;{\sf V})$ with ${\sf V}=\mathscr{S}^{{\rm O}(d)}$. It therefore remains to define the $(\infty,1)$-category of geometric structures \emph{with isotopies}.

\def\rloc{\mathbb{R}\text{-}{\rm loc}}
\def\smset{{\rm C}^{\infty}\mathscr{S}{\rm et}}
\def\smsset{{\rm C}^{\infty}\mathscr{S}}

\def\frakFEmb{\mathfrak{FEmb}}

The next definition is an adaptation of Grady--Pavlov \cite[Definition 4.1.5]{GradyPavlov}.
\begin{definition}\label{enrichedfemb}
We define the enriched site $\frakFEmb_d$, with enrichment in $\mathscr{S}$, as follows.
The objects are the same objects as in $\FEmb_d$.
Given two objects $W→U$ and $W'→U'$,
the corresponding hom-object between them is simplicial set
whose $n$-simplices are given by
\begin{itemize}
\item morphisms in $\FEmb_d$
of the form $(f:W\times \Delta^n_e→W',u:U\times \Delta^n_e→U')$ with source $(W\times \Delta^n_e→U\times \Delta^n_e)=(W→U)\times \id_{\Delta^n_e}$ and target $W'→U'$.
We require that the map $u$ is the composition of the projection $U\times \Delta^n_e\to U$ and a map $U\to U'$. 
\end{itemize}
The Grothendieck topology of $\frakFEmb_d$ is generated by the same coverage as for $\FEmb_d$.
\end{definition}

\begin{remark}
In Grady--Pavlov \cite{GradyPavlov}, we used a slightly different enrichment of $\frakFEmb_d$. The enrichment was in the $\RR$-local model structure on $\smset$. The $(\infty,1)$-categories presented by $\smset$ and $\smsset$ (with the $\RR$-local model structure) are both equivalent to the the $(\infty,1)$-category of spaces. In fact, the functor 
$${\sf B}_{\mathscr{S}}:\smset\into \smsset\to \mathscr{S}$$
is a Dwyer-Kan equivalence of corresponding simplicial categories. This follows immediately from existence of the transferred model structure on $\smset$ (Pavlov \cite[Theorem 1.2]{Pavlov2}). Application of $\csp_{\mathscr{S}}$ to the enriched homs in Grady--Pavlov \cite[Definition 4.1.5]{GradyPavlov} yield precisely the simplicial enrichment in Definition \ref{enrichedfemb}, hence we are free to use either enriched model for $\frakFEmb_d$. We choose the model in Definition \ref{enrichedfemb}, since it simplifies some of the proofs. 
\end{remark}

\def\frakStruct{\mathfrak{Struct}}
\begin{notation} Given $d\geq 0$, we denote by 
$$\frakStruct_{d}=\PSh(\frakFEmb_d,\mathscr{S})_{{\rm \check{C}ech}}$$
the relative category of enriched functors with weak equivalences given by {\v C}ech local weak equivalences. This relative category admits an $\mathscr{S}$-enriched model category presentation given by taking the enriched left Bousfield localization of the injective (or projective) model structure at {\v C}ech covers. 
The existence of the enriched left Bousfield localization exists by Grady--Pavlov \cite[Proposition 2.2.2]{GradyPavlov}.
\end{notation}

Now we can regard the non-enriched site ${\sf FEmb}_d$ (Definition \ref{femb}) as enriched in $\mathscr{S}$ by taking hom objects to be the constant presheaf on the set ${\sf FEmb}_d(M\to U, N\to V)$. Observe that there is a canonical enriched functor 
\begin{equation}\label{jfunctor}
j:{\sf FEmb}_d\to \frakFEmb_d
\end{equation}
which on homs is given by the inclusion as constant functions
$${\sf FEmb}_d(M\to U, N\to V)\into \frakFEmb_d(M\to U, N\to V)$$
The enriched left Kan extension along $j$ gives rise to a functor 
\begin{equation}\label{jfun}
j_!:\Struct_d\to \frakStruct_d
\end{equation}
Since $j_!$ sends representables to corresponding representables on the enriched site, it also sends {\v C}ech covers to {\v C}ech covers. Therefore, it defines a left Quillen functor in the local injective model structure. 

\begin{definition}\label{fiberwiseconc}
Fix $d\geq 0$. We define the \emph{isotopification functor} 
$$
\E_d:\Struct_d\to \frakStruct_d
$$
as the left derived functor $\E_d:=\mathbb{L} j_!$ of the functor $j_!$ in \eqref{jfun}, in the local injective model structure. 
We also define the endofunctor
$$\Ei_d:\Struct_d\to \Struct_d$$
as the composition $\Ei_d:=\mathbb{R}j^*\circ \mathbb{L}j_!$, where $\RR j^*$ is the right derived functor of the restriction functor $j^*$.
\end{definition}

\def\smequi{{\rm C}^{\infty}\mathscr{S}^{{\rm O}(d)}}
Our next goal is to define a functor 
$$
\R_d:\Struct_d\to \smequi
$$
We do this as follows.
\begin{itemize}
\item Given an object $X\in \frakStruct_d$, we can evaluate on projections of the form $p_U:\RR^d\times U\to U$. For each $U\in \cartsp$, this gives rise to a simplicial set $X(\RR^d\times U\to U)$ equipped with an action of the monoid of fiberwise open embeddings of $\RR^d\times U\to U$ into itself. More precisely, we have a map of simplicial sets
\begin{equation}\label{kister-maz}
\frakFEmb_d(\RR^d\times U\to U,\RR^d\times U\to U)\times X(\RR^d\times U\to U)\to X(\RR^d\times U\to U),
\end{equation}
which is obtained by the enriched presheaf structure maps of $X$. 

\item There is a canonical monomorphism of simplicial sets
\begin{equation}\label{ontofemb}
C^{\infty}(U,{\rm O}(d))\to \frakFEmb_d(\RR^d\times U\to U,\RR^d\times U\to U),
\end{equation}
that sends an $n$-simplex $\sigma:\Delta^n_e\times U\to {\rm O}(d)$ to the fiberwise open embeddings $\bar{\sigma}:\Delta^n_e\times \RR^d\times U\to \RR^d\times U$ defined by 
$$\bar\sigma(t,x,u)=(\sigma(t,u)(x),u), \quad t\in \Delta^n_e,x\in \RR^d,u\in U\;.$$
Restricting the map \eqref{kister-maz} along the map \eqref{ontofemb} equips $X(\RR^d\times U\to U)$ with an action of $C^{\infty}(U,{\rm O}(d))$. 

\item We have a canonical inclusion of simplicial sets ${\rm O}(d)\into C^{\infty}(U,{\rm O}(d))$, which picks out constant maps in the $U$-direction. The simplicial set $X(\RR^d\times U\to U)$ admits an action of ${\rm O}(d)$ by further restricting along this inclusion. 
\end{itemize}

\begin{definition}\label{kfunctor}
Fix $d\geq 0$. The construction defined in the previous paragraph gives rise to a functor, which we call the \emph{restriction functor}. We denote this functor by 
$$
\R_d:\frakStruct_d\to \smequi\;,
$$
Summarizing the above construction, this functor restricts a presheaf on the enriched site to the following subcategory:
\begin{itemize}
\item  Objects are given by projections $p_U:\RR^d\times U\to U$ 
\item For two objects $p_U:\RR^d\times U\to U$ and $p_{V}:\RR^d\times V\to V$, the simplicial set of morphisms has $n$-simplices given by pairs $(f\times g,g)$, where $g:U\to V$ and $f:\Delta^n_e\to {\rm O}(n)$ are smooth maps. The map $f$ acts as a $\Delta^n_e$-parametrized family of diffeomorphisms on $\RR^d$ by $f(x,t)=f(t)\cdot x$, for all $(x,t)\in \RR^d\times \Delta^n_e$.
\end{itemize}
This functor is a left Quillen functor in a Quillen equivalence, by \cite[Proposition 3.3.12]{GradyPavlov2}. 
\end{definition}
\def\Sec{{\sf Sec}}

\begin{remark}\label{scanning}
The functor $\R_d$ admits a homotopy right adjoint $\Sec_d$, which sends a presheaf $X$ on $\cartsp$ with values in $\mathscr{S}^{{\rm O}(d)}$ to the presheaf on $\frakFEmb_d$ that sends $p:M\to U$ to the simplicial set of sections 
$$\Gamma(M;{\rm Fr}(p)\times_{{\rm O}(d)}X(U)),$$
where ${\rm Fr}(p)$ is the bundle of orthonormal frames of the vertical tangent bundle associated to $p$ and ${\rm O}(d)$ acts diagonally on both factors. The $\mathcal{S}$-component of the unit of the adjunction is a map of the form 
$$\mathcal{S}(p:M\to U)\to \Gamma(M;{\rm Fr}(p)\times_{{\rm O}(d)}\mathcal{S}(\RR^d\times U\to U)).$$
It is families variant of Segal's scanning map. 
\end{remark}

The functor $\R_d$ is a key ingredient in the "deformation space" construction for geometric structures. For example, if the geometric structure is smooth families of flat $H_d$-bundles, application of $\R_d$ will just restrict to flat $H_d$-bundles on $\RR^d$, parametrized by open subspaces of cartesian spaces. It is easy to see that the space of all fiberwise flat $H_d$-bundles on $\RR^d\times U\to U$ is equivalent to $\deloop {\rm C}^{\infty}(U,H_d)$, since every fiberwise flat $H_d$-bundle on $\RR^d$ is trivializable by a connection preserving trivialization, and the automorphisms of the fiberwise trivial $H_d$-bundle on $\RR^d$ are smooth families of elements in $H_d$, or smooth maps $U\to H_d$.

\begin{remark}\label{enrichall}
In the following, we will always work with presheaves on the enriched site $\frakFEmb_d$. Given a presheaf $\mathcal{S}$ on the non-enriched site $\FEmb_d$, we will always implicitly convert it to a presheaf on the enriched site via the derived Kan extension of $\E_d$ in \eqref{fiberwiseconc}:
$$\E_d:\Struct_d\to \frakStruct_d.$$
We remind the reader that the definition of the bordism category in Grady--Pavlov \cite{GradyPavlov2} uses enriched sheaves. When $\mathcal{S}\in \Struct_d$, then $\mathfrak{Bord}_d^{\mathcal{S}}$ means that a family of bordisms $p:M\to U$ is equipped with a morphism in $\frakFEmb_d$, $(p:M\to U)\to \E_d(\mathcal{S})$. This has the effect of adding isotopies to geometric structure. We defer the reader to \cite[Section 4]{GradyPavlov} for more details.
\end{remark}

We now define a functor that constructs a tangential structure, in the sense of Galatius--Madsen--Tillmann--Weiss \cite{GMTW}, from a geometric structure. We remind the reader that the functor $\csp_{\mathscr{S}^{{\rm O}(d)}}$ from Definition \ref{struts}, with ${\sf V}=\mathscr{S}^{{\rm O}(d)}$, sends a presheaf on $\cartsp$ with values in ${\rm O}(d)$-equivariant spaces, to ${\rm O}(d)$-equivariant space. 
\begin{definition}\label{dfunctor}
Fix $d\geq 0$. We define the \emph{Deformation of structures functor} 
$$
\D_d:\frakStruct_d\to \mathscr{S}^{{\rm O}(d)}
$$ 
as the composition $\D_d:=\csp_{\mathscr{S}^{{\rm O}(d)}}\R_d$. Given a geometric structure $\mathcal{S}\in \frakStruct_d$, we define its \emph{deformation space} by $\D_d(\mathcal{S})$.
\end{definition}

\def\lSec{l{\sf Sec}}

Since both functors $\R_d$ and $\csp_{\mathscr{S}^{{\rm O}(d)}}$ are homotopy cocontinuous,  the composite functor $\D_d=\csp_{\mathscr{S}^{{\rm O}(d)}}\R_d$ admits a homotopy right adjoint, which we denote by  
$$\lSec_d:\mathscr{S}^{{\rm O}(d)}\to \frakStruct_d.$$
We now give an explicit presentation of this functor. Let $\Sec_d$ be the right derived Quillen adjoint of $\R_d$ in \eqref{kfunctor}, which exists since $\R_d$ is a left Quillen functor (see Remark \ref{scanning} for a brief description). Let $\delta_{\mathscr{S}^{{\rm O}(d)}}:\mathscr{S}^{{\rm O}(d)}\to \smequi$ denote the constant sheaf functor Definition \ref{constsheaf}. We have the following.

\begin{proposition}\label{righadjbe}
Fix $d\geq 0$ and let $\delta_{\mathscr{S}^{{\rm O}(d)}}$ and $\Sec_d$ be as in the preceding paragraph. Then $\lSec_d=\Sec_d\circ \delta_{\mathscr{S}^{{\rm O}(d)}}$ is homotopy right adjoint to the functor $\D_d=\csp_{\mathscr{S}^{{\rm O}(d)}}\R_d$. 
\end{proposition}
\proof
Since $\delta_{\mathscr{S}^{{\rm O}(d)}}$ is homotopy right adjoint to $\csp_{\mathscr{S}^{{\rm O}(d)}}$, by Proposition \ref{combhleft}, this follows at once from the fact that $\Sec_d$ is homotopy right adjoint to $\R_d$. 

\endofproof

The following example provides some motivation for the equivalence between (1) and (4) in Theorem \ref{freedhopconj}.

\begin{example}\label{flatbunexample}
Consider the simplicial presheaf $H_d\text{-}{\rm Bun}_{{\rm fl}\text{-}\nabla}$ on $\FEmb_d$, which assigns a submersion $p:M\to U$ with $d$-dimensional fibers to the groupoid of principal $H_d$ bundles on $M$ with fiberwise \emph{flat} connection on $p$. The left kan extension $\E_d$ applied to $H_d\text{-}{\rm Bun}_{{\rm fl}\text{-}\nabla}$ can be described explicitly via the formula for the left Kan extension as follows. 
\begin{itemize}
\item For $p:M\to U\in \frakFEmb_d$, the simplicial set $\E_d(H_d\text{-}{\rm Bun}_{{\rm fl}\text{-}\nabla})(p:M\to U)$ has $n$-simplices given by equivalence classes of pairs $(g,(P,\nabla))$ where $g:\Delta^n_e\times M\to M$ is a fiberwise embedding (over $\Delta^n_e$), and $(P,\nabla)$ is a fiberwise flat $H_d$-bundle on $M$. 
\end{itemize}
Now a flat $H_d$-bundle on $\RR^d$ is trivializable by a connection preserving trivialization. Moreover, automorphisms of the trivial bundle with the trivial connection can be identified with elements of $H^d$. Taking into account smooth $U$-families of such objects, we identify the resulting simplicial set (when $U$ is contractible) as 
$$\E_d\left(H_d\text{-}{\rm Bun}_{{\rm fl}\text{-}\nabla}\right)(\RR^d\times U\to U)\simeq {\rm O}(d)\times \deloop {\rm C}^{\infty}(U,H_d),$$
The second factor on the right is the constant smooth simplicial set on $\deloop {\rm C}^{\infty}(U,H_d)$. The first factor on the right is the constant smooth simplicial set on the simplicial group ${\rm O}(d)$. This identification comes from identifying the space of maps $g$ in the above bullet point via the Kister-Mazur theorem 
\begin{align*}
\frakFEmb_d(\RR^d\times U\to U,\RR^d\times U\to U) & \simeq \map(U,\frakFEmb_d(\RR^d,\RR^d))
\\
&\simeq \frakFEmb_d(\RR^d,\RR^d)
\\
& \simeq {\rm O}(d).
\end{align*}
The first equivalence is just an unraveling of the definitions. The second equivalence holds when $U$ is contractible. The last equivalence is a consequence of the Kister--Mazur theorem. In the last line ${\rm O}(d)$ means the singular simplicial set of ${\rm O}(d)$. 
  
Since $\R_d$ simply evaluates on projections $p_U:\RR^d\times U\to U$, and then applies $\csp$, we have a local equivalence of presheaves of ${\rm O}(d)$-equivariant spaces
$$\R_d\E_d\left(H_d\text{-}{\rm Bun}_{{\rm fl}\text{-}\nabla}\right)\simeq  {\rm O}(d) \times \mathbf{B} H_d,$$
where $\mathbf{B}H_d$ is the moduli stack of trivial principal $H_d$-bundles. In other words, we have the free ${\rm O}(d)$-equivariant space on the stack of principal $H_d$-bundles. 

The takeaway from this example is the appearance of the moduli space of \emph{all} principal $H_d$-bundles on the base space of the parametrizing family, starting from the geometric structure given by fiberwise \emph{flat} $H_d$-bundles. After passing to smooth deformation spaces, this gives a further identification 
$$\D_d\E_d\left(H_d\text{-}{\rm Bun}_{{\rm fl}\text{-}\nabla}\right)\simeq  {\rm O}(d) \times \deloop H_d,$$
where on the right we see the whole homotopy type of $H_d$, not just the underlying discrete group. 
\end{example}

The following table serves as a reference for the various functors we have established so far. 
\begin{center}
{\renewcommand{\arraystretch}{2}
\begin{tabular}{c|c|p{6.5cm}}
{\bf Symbol} & {\bf Name} & {\bf Description}
\\
\hline
$\gamma_{\sf V}:\smV\to {\sf V}$ & global sections functor & Evaluates at the point $\RR^0\in \cartsp$
\\
\hline
$\delta_{\sf V}:{\sf V}\to \smV$ & constant sheaf functor & Sends an object $v\in {\sf V}$ to the sheaf whose value at $U$ is the powering $v^{{\rm sing}(U)}$.
\\
\hline
$\csp_{\sf V}: \smV\to {\sf V}$ & deformation functor & Evaluates a sheaf on extended $n$-simplces and glues together an object in ${\sf V}$ from this data. 
\\
\hline
 $\mathcal{B}_{\sf V}:\smV\to \smV$ & internal deformation functor & The composition $\delta_{\sf V}\csp_{\sf V}$
\\
\hline 
 $\R_d:\frakStruct_d\to \smequi$ & restriction functor & Sends a simplicial presheaf to its value on representables of the form $\RR^d$. 
\\
\hline
$\D_d:\frakStruct\to \mathscr{S}^{{\rm O}(d)}$ & deformation of structures functor & The composition $\csp_{\sf \mathscr{S}^{{\rm O}(d)}} \R_d$, 
\\
\hline
 $\Sec_d: \smequi\to \frakStruct_d$ & fiberwise bundle functor & Sends a smooth equivariant space to the sections of the bundle of smooth spaces associated to the fiberwise frame bundle. 
\\
\hline 
$\lSec_d:\mathscr{S}^{{\rm O}(d)}\to \smequi$ & constant bundle functor & The composition $\Sec_d\delta_{\mathscr{S}^{{\rm O}(d)}}$
\\
\hline
$\E_d:\Struct_d\to \frakStruct_d$ & isotopification functor & Adds isotopies to geometric structures. 
\end{tabular}}
\end{center}

\medskip
The following proposition establishes relationships between all the adjoint functors we have established. 

\begin{proposition}\label{adjointcalc}
Let ${\sf V}$ be a combinatorial and simplicial model category. Recall the functors  $\gamma_{\sf V}$, $\delta_{\sf V}$, $\csp_{\sf V}$, $\Icsp_{\sf V}$, $\R_d$, $\Sec_d$, $\D_d$, $\lSec_d$ from Definitions \ref{constglob}, \ref{constsheaf}, \ref{struts}, \ref{fiberwiseconc}, \ref{kfunctor} and  Proposition \ref{righadjbe} (respectively). We have the following homotopy adjunctions:
\begin{itemize}
\item $\delta_{\sf V}\dashv \gamma_{\sf V}$
\item $\csp_{\sf V}\dashv \delta_{\sf V}$
\item $\R_d\dashv \Sec_d$
\item $\D_d\dashv \lSec_d$.
\item $\Icsp_{\sf V}=\delta_{\sf V}\csp_{\sf V}\dashv \delta_{\sf V}\gamma_{\sf V}$
\end{itemize}
Moreover, the functors obey the following calculus:
\begin{enumerate}
\item $\D_d=\csp_{\mathscr{S}^{{\rm O}(d)}}\R_d$
\item $\lSec_d=\Sec_d\delta_{\mathscr{S}^{{\rm O}(d)}}$
\item $\id\overset{\simeq}{\to}\gamma_{\sf V}\delta_{\sf V}$
\item $\id\overset{\simeq}{\to}\csp_{\sf V}\delta_{\sf V}$
\item $\R_d \Sec_d\overset{\simeq}{\to} \id$ and $\Sec_d\R_d\overset{\simeq}{\leftarrow} \id$.
\item $\Icsp_{\sf V}=\delta_{\sf V}\csp_{\sf V}\dashv \delta_{\sf V}\gamma_{\sf V}$
\end{enumerate}
where $\simeq$ are natural weak equivalences in the corresponding relative category.
\end{proposition}
\proof
The homotopy adjunction $\csp_{\sf V}\dashv \delta_{\sf V}$ is Proposition \ref{combhleft}. The homotopy adjunction $\R_d\dashv \Sec_d$ follows from the definition of $\R_d$. This is part of a Quillen equivalence by \cite[Proposition 3.3.12]{GradyPavlov2}, which implies the two relations in (5), where the equivalences are given by the derived unit and counit of the Quillen adjunction. The homotopy adjunction $\delta_{\sf V}\dashv \gamma_{\sf V}$ follows from the homotopy adjunction $\delta^{pre}_{\sf V}\dashv \gamma_{\sf V}$ in Definition \ref{constglob}, along with the fact that the canonical map
$\delta^{pre}_{\sf V}\to \delta_{\sf V},$
described in the paragraph preceding Lemma \ref{contractopens}, is a local equivalence, by Lemma \ref{contractopens}. The homotopy adjunction $\D_d\dashv \lSec_d$ follows by composing adjunctions. Relations (1) and (2) are tautological. The remaining relations (3) and (4) are proved as follows. Since $\gamma_{\sf V}$ sends local weak equivalences to weak equivalences, the local weak equivalence $\delta^{pre}_{\sf V}\to \delta_{\sf V}$ gives rise to a weak equivalence 
$$\id=\gamma_{\sf V}\delta^{pre}_{\sf V}\to \gamma_{\sf V}\delta_{\sf V}.$$
The homotopy adjunction $\Icsp_{\sf V}=\delta_{\sf V}\csp_{\sf V}\dashv \delta_{\sf V}\gamma_{\sf V}$ follows by composing adjoint functors. Similarly, for all $X\in \mathscr{S}$, the composition $\csp\delta^{pre}_{\sf V}(X)$ is the homotopy colimit over the constant diagram on $X$, indexed by $\Delta^{\rm op}$. Since $\Delta$ has a terminal object, its nerve is contractible and the homotopy colimit is just $X$. Therefore, composing the local weak equivalence $\delta^{pre}_{\sf V}\to \delta_{\sf V}$ with $\csp_{\sf V}$ gives rise to a local weak equivalence $\id \overset{\simeq}{\to} \csp_{\sf V}\delta_{\sf V}$.
\endofproof

\subsection{Detecting deformation space equivalences}

We now turn out attention to establishing a convenient criteria for determining when a morphism of presheaves $f:X\to Y$, $X,Y\in \Struct_d$ induces an equivalence after passing to deformation spaces, via the functor $\D_d\E_d$. 

\begin{definition}
Fix $n\in \NN$ and $i\in [n]$. Let $\Delta^n_s$ denote the simplicial $n$-simplex and let $\Delta^n_{st}$ denote the standard $n$-simplex in $\RR^{n+1}$. For $\epsilon>0$, we also let $\Delta^n_{st,\epsilon}$ denote the $\epsilon$-neighborhood of the standard $n$-simplex in the extended $n$-simplex $\Delta^n_e$ (Definition \ref{extsimp}). We define various face and coface maps as follows.
\begin{itemize}
\item We let $\dg^i:\Delta^{n-1}_e\to \Delta^n_e$ denote the inclusion of the $i$th face into the extended $n$-simplex, i.e. the $i$th geometric coface map. We denote the corresponding restriction map by $\dg_i$, i.e. the $i$th face map.
\item We let $d^i:\Delta^{n-1}_s\to \Delta^{n}_s$ denote the inclusion of the $i$th face into the simplicial $n$ simplex, i.e. the $i$th coface map. We denote the induced restriction maps by lower indices $d_i$, i.e. the $i$th geometric face map. 
\item We define the inclusion of the $\epsilon$-neighborhood of the $i$th face of the standard $n$-simplex as $\dg^i_{\epsilon}:\Delta^{n-1}_{st,\epsilon}\into \Delta^n_e$. Again, lower indices indicate the corresponding restriction. 
\end{itemize}
\end{definition}

As is familiar from motivic homotopy theory \cite{Voevodsky}, there are two sorts of homotopies between morphisms in $\sm\mathscr{S}$. One is the usual simplicial homotopy. The other is a \emph{geometric} homotopy, which is defined using the extended $1$-simplex $\Delta^1_e$. Homotopies using the geometric version of homotopy are called \emph{concordances}. We recall the definition.

\begin{definition}\label{concordance}
Let $X,Y\in \smsset$. Two morphisms $f,g:X\to Y$ are said to be \emph{concordant} is there exists a morphism $H:X\times \Delta_e^1\to  Y$ such that the restrictions $\dg_0^*H=f$ and $\dg_1^*H=g$. A map $f:X\to Y$ is called a \emph{concordance equivalence} if there is a map $g:Y\to X$ such that $fg$ and $gf$ are both concordant to the identity.
\end{definition}

We can also define concordant equivalent sections of a presheaf as follows.
\begin{definition}
Let $X\in \smsset$ be an objectwise Kan complex satisfying homotopy descent. Let $U\in \cartsp$. Let $s,t\in X_0(U)$ and let $A\subset U$ be a closed subset. 

\begin{itemize}
\item We say that $s$ and $t$ are \emph{concordant} if there exists $\epsilon>0$ and $H\in X_1(\Delta^1_e\times U)$ such that $\dg_0d_0H=p_U^*s$ and $\dg_1d_1H=p_{U}^*t$, where $p_U^*$ is restriction along the projection $U\times \Delta^0_{st,\epsilon}\to U$. 

\item We say that $s$ and $t$ are \emph{concordant relative $A$} if there exists a concordance $H\in X_1(\Delta^1_e\times U)$ connecting $t$ and $s$ that restricts to a constant concordance along the inclusion $\Delta^1_e\times V\into \Delta^1_e\times U$, where $A\subset V$ is an open neighborhood. 
\end{itemize}
One can use the the $(\infty,1)$-sheaf condition and Kan fibrancy to prove that concordance equivalence defines an equivalence relation on sections \cite{Pavlov}. 

We denote the set of concordance classes by $X[U]$. We denote the set of relative concordance classes as $X[U;A]$. 
\end{definition}

\begin{notation}
Recall the functors $\csp_{\mathscr{S}}$, $\gamma_{\mathscr{S}}$, $\delta_{\mathscr{S}}$, $\mathcal{C}_{\mathscr{S}}$ defined in Definition \ref{constglob}, \ref{constsheaf}, and \ref{struts} for ${\sf V}=\mathscr{S}$. Throughout this section, we will use the simplified notation $\csp:=\csp_{\mathscr{S}}$, $\gamma:=\gamma_{\mathscr{S}}$, $\delta:=\delta_{\mathscr{S}}$, $\mathcal{C}:=\mathcal{C}_{\mathscr{S}}$.
\end{notation}
In the following, we will need to use all three types of simplices $\Delta^n_s$, $\Delta^n_{st,\epsilon}$ and $\Delta^n_e$.

\begin{lemma}\label{equivssimp}
Let $n\geq 0$ and let $\Delta^n_s$ denote the simplicial $n$-simplex. 
The map $\Delta^n_s\to \csp(\Delta^n_e)$ that picks out the identity on $\Delta^n_e$ is a weak equivalence. Similarly, the canonical map $\partial \Delta^n_s\to \csp(\partial \Delta^n_e)$ is an equivalence. 
\end{lemma}
\proof
As a model for the homotopy colimit $\csp(\Delta^n_e)$ in Definition \ref{struts}, we can take the smooth singular simplicial set of $\Delta^n_e$. It is well known that the inclusion of the smooth singular simplicial set into the continuous singular simplicial set is a weak homotopy equivalence. Since the extended $n$-simplex $\Delta^n_e$ deformation retracts to the standard $n$-simplex $\Delta^n_{st}\subset \RR^{n+1}$, it follows that $\csp(\Delta^n_e)$ is weak equivalent to the usual singular simplicial set of the standard $n$-simplex. Since the geometric realization of $\Delta^n_s$ is equivalent to the standard (topological) $n$-simplex, it follows from the preceding discussion that if we apply the geometric realization functor to map $\Delta^n_s\to \csp(\Delta^n_e)$, the result is weakly equivalent to the unit of the adjunction $\eta_{\Delta^n_{st}}:\Delta^n_{st} \to \vert \mathcal{N}(\Delta^n_{st}) \vert$, where $\mathcal{N}$ is the singular simplicial nerve functor. But $\eta_{\Delta^n_{st}}$ is a weak equivalence, so we are done.

The proof for $\partial \Delta^n_s$ proceeds in a similar fashion.  The geometric realization of $\partial \Delta^n_s$ is weakly equivalent to the standard topological boundary, so the only nontrivial step is to see that $\csp(\partial \Delta^n)$ is equivalent to the continuous singular simplicial set of the topological boundary $\partial \Delta^n_{st}$. For this, we observe that the strict colimit diagram \eqref{boundarynsimp} is also a homotopy colimit diagram, since the horizontal maps are both injective cofibrations. Using the fact that homotopy colimits commute and the fact that $\csp(\Delta^n_e)$ is weak equivalent to $\Delta^n_{s}$, as shown in  the preceding paragraph, we deduce that $\csp(\partial \Delta^n_e)$ is weak equivalent to $\partial \Delta^n_s$, as claimed. 
\endofproof

Let $\epsilon\ll 1$ and let $\partial \Delta^n_{st,\epsilon}$ denote the $\epsilon$-neighborhood of the boundary of the standard $n$-simplex in $\RR^{n+1}$. Fix a map 
$$r^i_{n-1}:\Delta^{n-1}_e\to \Delta^{n-1}_{st}\overset{d^i}{\into} \partial \Delta^n_{st,\epsilon}$$
that retracts the extended $(n-1)$-simplex to the standard (compact) $(n-1)$-simplex and then includes it into the boundary. By the universal property of the colimit \eqref{boundarynsimp}, the $r^i_{n-1}$'s induce a map 
\begin{equation}\label{rmap}
r:\partial \Delta^n_{e}\to \partial \Delta_{st,\epsilon}^n.
\end{equation}
\begin{lemma}\label{siboundary}
For $\epsilon\ll 1$, the map $r$ in \eqref{rmap} induces a weak equivalence 
$$\csp(r):\csp(\partial \Delta^n_{st,\epsilon})\overset{\simeq}{\to} \csp(\partial \Delta^n_e).$$
\end{lemma}
\proof
For all $\epsilon\ll 1$, we can construct a cover $\{U_{i}\}_{i=0}^{n}$ of $\partial \Delta^n_{st,\epsilon}$ that satisfies the following:
\begin{itemize}
\item Each $U_i$ contains the $i$th face of $\partial \Delta^n_{st}$ 
\item Each $U_i$ deformation retracts smoothly onto the $i$-th face.
\item Each intersection $U_i\cap U_j$, $i<j$ smooth deformation retracts to the common face $\Delta^{n-2}_{st}$ indexed by $i<j$. \end{itemize}
For $n=2$, such a cover is depicted below. 
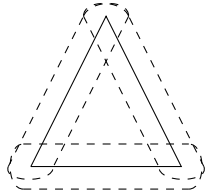
\begin{center}
\begin{tikzpicture}
\draw (0,0) -- (1,2) -- (2,0) -- (0,0);
\draw[dashed] (-.3,0) -- (.7,2);
\draw[dashed] (1.3,2) -- (.3,0);
\draw[dashed] (.7,2) to[out=70, in=90] (1.3,2);
\draw[dashed] (-.3,0) to[out=250, in=270] (.3,0);
\begin{scope}[yscale=1,xscale=-1,xshift=-2cm]
\draw[dashed] (-.3,0) -- (.7,2);
\draw[dashed] (1.3,2) -- (.3,0);
\draw[dashed] (.7,2) to[out=70, in=90] (1.3,2);
\draw[dashed] (-.3,0) to[out=250, in=270] (.3,0);
\end{scope}
\draw[dashed] (-.1,.3) -- (2.1,.3);
\draw[dashed] (-.1,-.3) -- (2.1,-.3);
\draw[dashed] (-.1,-.3) to[out=180, in=180] (-.1,.3);
\draw[dashed] (2.1,-.3) to[out=0, in=0] (2.1,.3);

\end{tikzpicture}
\end{center}
Let $C(\{U_i\})$ denote the {\v C}ech nerve of this cover. Then we have a local equivalence $C(\{U_i\})\to \partial \Delta^n_{st,\epsilon}$. Observe that by construction, each map $r^i_{n-1}:\Delta^{n-1}_e\to \partial \Delta^n_{st,\epsilon}$ inducing the map  \eqref{rmap} factors through $U_i$, and the restriction to the common face $\Delta^{n-2}_e$ factors through $U_i\cap U_j$. Hence, the induced map $r:\partial \Delta^n_e\to \partial \Delta^n_{st,\epsilon}$ factors through the monomorphism $i:C(\{U_i\})\into \partial \Delta^n_{st,\epsilon}$. Since each $U_i$ smooth deformation retracts to the $i$-th face of $\partial \Delta^n_{st}$, the factorizing 	 map $j:\partial \Delta^n_e\to C(\{U_i\})$, with $ij=r$, is a concordance equivalence. Since $\csp$ preserves both local equivalences and concordance equivalences, it follows that both $\csp(i)$ and $\csp(j)$ are equivalences. Hence $\csp(r)$ is also a weak equivalence.  
\endofproof

Fix some small $0<\epsilon\ll 1$ and let $n\geq 1$. Fix a global point $x:\ast\to X$ of an $(\infty,1)$-sheaf $X\in \smsset$. Let 
$$\overline{x}:\ast\to X(\partial \Delta^n_{st,\epsilon})$$
be the map that picks out the vertex given by the composition
$$
\partial \Delta^n_{st,\epsilon} \to \ast \overset{x}{\to} X\;.
$$
We define the space $X(\Delta^n_e;x,\epsilon)$ as the homotopy fiber in simplicial sets 
\begin{equation}\label{fiberhomsheav}
\xymatrix{
X(\Delta^n_e;x,\epsilon) \ar[rr]\ar[d] && \ast\ar[d]^-{\overline{x}}
\\
X(\Delta^n_e) \ar[rr]^-{r^*}  && X(\partial \Delta^n_{st,\epsilon}),
}
\end{equation}
where the map $r^*$ is induced by restriction along the map \eqref{rmap}. We remove dependence on $\epsilon$ by taking the germ:
$$X(\Delta^n_e;x):=\colim_{\epsilon\to 0}X(\Delta^n;x,\epsilon),$$
where the colimit is taken over filtered diagram of neighborhood inclusions, whose transition maps are $\partial \Delta^n_{st,\epsilon}\into \partial \Delta^n_{st,\epsilon'}$, for $\epsilon<\epsilon'$. 

Recall the concordification functor $\Con$ and internal deformation space functor $\Icsp$ in Definition \ref{struts}. Recall also the global sections functor $\gamma$ in Definition \ref{constglob}. Notice that $\gamma\Con=\csp$. Following Berwick-Evans -- Boavita de Brtio--Pavlov \cite{Pavlov}, we define a comparison map 
$$\overline{c}:\pi_0(\Con(X)(\Delta^n_e;x))\to \pi_n(\csp(X),x),$$
as follows. There is a natural transformation $\Con\to \Icsp$ of the form
\begin{equation}
ev_{X}:\Con(X)\to \delta\gamma \Con (X)= \Icsp(X),
\end{equation}
called the \emph{evaluation map} (see \cite[Proposition 2.17]{Pavlov}). The evaluation map gives rise to a further map
\begin{equation}\label{concmap}
\Con(X)(\Delta^n_e)\to \map(\Delta^n_e,\Icsp(X))\overset{\simeq}{\to}\map(\csp(\Delta^n_e),\csp(X))\overset{\simeq}{\to} \map(\Delta^n_s,\csp(X)),
\end{equation}
where the middle map is the weak equivalence obtained by first applying the enriched functor 
$$\csp:\map(\Delta^n_e,\Icsp(X))\to \map(\csp(\Delta^n_e),\csp\delta\csp(X))$$ 
and then composing with the the map obtained by applying the functor $\map(\csp(\Delta^n_e),-)$ to the map
$$
\csp\delta(Y)=\hocolim_{[n]\in \Delta^{\rm op}}Y^{{\rm sing}(\Delta^n_e)}\to \hocolim_{[n]\in \Delta^{\rm op}}Y^{\Delta^n_s}=Y,
$$
with $Y=\csp(X)$. This map is a weak equivalence, as it witnesses the homotopy adjunction $\csp\dashv \delta$. The second map in \eqref{concmap} is obtained by precomposing with the map $\Delta^n_s\to \csp(\Delta^n_e)$ in Proposition \ref{equivssimp}. An analogous map can be constructed with $\partial \Delta^n_s$ replacing $\Delta^n_s$ and $\partial \Delta^n_{st,\epsilon}$ replacing $\Delta^n_e$. 

One checks that we have a commutative diagram: 
\begin{equation}\label{commdiag}
\xymatrix{
\Con(X)(\Delta^n_e) \ar[rr]\ar[d] && \map(\Delta^n_s,\csp(X))\ar[d]
\\
\Con(X)(\partial \Delta^n_{st,\epsilon}) \ar[rr] && \map(\partial \Delta^n_s,\csp(X)),
}
\end{equation}
where the horizontal maps are the ones described in the preceding paragraph and the vertical maps are induced by the inclusion of the boundary. The diagram induces a map on homotopy fibers:
\begin{equation}\label{isosim}
c_{\epsilon}:\Con (X)(\Delta^n_e;x,\epsilon)\to \map_{\csp(x)}(\Delta^n_s/\partial \Delta^n_s,X),
\end{equation}
where the right side is the space of simplicial maps $\Delta^n_s\to \csp(X)$ that restrict to the map $\csp(x):\ast\to \csp(X)$ on the boundary. Passing to the colimit as $\epsilon\to 0$, the $c_{\epsilon}$ maps induce a map 
\begin{equation}\label{isosimp}
c:\Con (X)(\Delta^n_e;x)\to \map_{\csp(x)}(\Delta^n_s/\partial \Delta^n_s,X),
\end{equation}
which yields the desired map on $\pi_0$.
\begin{proposition}\label{hgroupiso}
The map $c$ in \eqref{isosimp} is an equivalence. In particular, $c$ induces an isomorphism on relative concordance classes
$$X[\Delta^n_e;x]\cong \pi_0(\Con(X)(\Delta^n_e;x))\to \pi_n(\csp(X),x).$$
\end{proposition}
\proof
It suffices to show that for sufficiently small $0<\epsilon\ll 1$, the map $c_{\epsilon}$ induces an equivalence. Indeed, taking homotopy fibers commutes with taking filtered colimits. Moreover, weak equivalences are stable under filtered colimits. Hence, if $c_{\epsilon}$ is a weak equivalence for small $\epsilon$, then the $c_{\epsilon}$ will assemble to the weak equivalence $c$.  

To show that $c_{\epsilon}$ is an equivalence, observe that $c_{\epsilon}$ is the map homotopy fibers that is induced by the two horizontal maps in \eqref{commdiag}. Therefore, it suffices to prove that these two horizontal maps are weak equivalences. By the main theorem of Berwick-Evans -- Boavita de Brito -- Pavlov \cite{Pavlov}, the evaluation map $ev:\Con\to \Icsp$ is a natural weak equivalence. By Lemma \ref{equivssimp}, the composite map \eqref{concmap} is an equivalence. By Lemma \ref{equivssimp} and Lemma \ref{siboundary}, the map obtained by replacing $\Delta^n_s$ by $\partial \Delta^n_s$ and $\Delta^n_e$ by $\partial \Delta^n_{st,\epsilon}$ in \eqref{concmap} is a weak equivalence. These two maps are precisely the two horizontal maps in \eqref{commdiag}, so we are done. 
\endofproof

\begin{corollary}\label{concequiv}
Let $X,Y\in \smsset$ and let $f:X\to Y$ be a morphism. Fix a point $x:\Delta^0_e\to X$. Suppose $f$ induces an isomorphism on concordance classes $f_*:X[\Delta^0_e]\to Y[\Delta^0_e]$ and an isomorphism on relative concordance classes 
$$f_*:X[\Delta^n_e;x]\to Y[\Delta^n_e;f(x)],$$
for all $n\geq 1$. Then $f$ induces a weak equivalence of spaces
$$\csp(f):\csp(X)\to \csp(Y).$$
\end{corollary}
\proof
By Proposition \ref{hgroupiso}, it follows at once that $\csp(f)$ induces an isomorphism on $\pi_n$ for all $n\geq 1$. For $n=0$, observe that $\csp(X)$ can be computed by taking the diagonal of the bisimplicial set $([n],[m])\mapsto X_m(\Delta^n_e)$. There is a canonical surjective map
$$
\phi:X[\Delta^0_e]\onto \pi_0(\csp(X))
$$
which sends concordances classes, where the concordances are constant on a germ around $0$ and $1$ in the extended 1-simplex, to corresponding concordance classes with no restrictions on the concordance. We claim that $\phi$ is also injective. Let $s,t\in X(\Delta^0_e)$ be such that $[s]=[t]\in \pi_0(\csp(X))$. Then there is a concordance $H\in X_1(\Delta^1_e)$ joining $s$ and $t$. Let $0<\delta'<\delta \ll 1$ be a small real numbers and let $\psi:\RR\to \RR$ be a smooth function that satisfies the following conditions:
\begin{itemize}
\item $\psi'\geq 0$
\item $\psi(x)=x$ on $\RR\setminus B_{\delta}(0)\cup B_{\delta}(1)$.
\item $\psi(x)=1$ on $B_{\delta'}(1)$
\item $\psi(x)=0$ on $B_{\delta'}(0)$
\end{itemize}

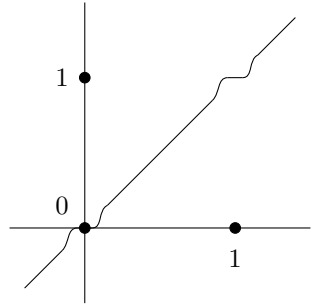
\begin{center}
\begin{tikzpicture}
\draw (-2,-1) -- (2,-1);
\draw (-1,-2) -- (-1,2);
\filldraw (-1,-1) circle (2pt);

\draw (-.7,-.7) -- (.7,.7);
\draw (-.7,-.7) to[out=205, in =0] (-.9,-1);
\draw (.7,.7) to[out=45, in =180] (.9,1);
\draw (-.9,-1) -- (-1.1,-1);
\draw (.9,1) -- (1.1,1);
\draw (-1.1,-1) to[out=180, in =45] (-1.3,-1.3);
\draw (1.1,1) to[out=0, in =205] (1.3,1.3);
\draw (-1.3,-1.3) -- (-1.8,-1.8);
\draw (1.3,1.3) -- (1.8,1.8);

\filldraw (-1,1) circle (2pt);
\filldraw (1,-1) circle (2pt);

\node at (-1.3,1) {$1$};
\node at (1,-1.4) {$1$};
\node at (-1.3,-.7) {$0$};
\end{tikzpicture}
\end{center}

Then pulling back the concordance $H\in X_1(\Delta^1_e)$ by $\psi:\Delta^1_e\to \Delta^1_e$ yields a concordance joining $s$ and $t$ that is constant on the $\delta'$ interval around $0$ and $1$. Hence $[s]=[t]\in X[\Delta^0_e]$.
\endofproof

In practice, it is often convenient to just verify the following stronger condition. In the case where the sheaves are ordinary presheaves (as opposed to simplicial presheaves), this strategy for proving that one has an equivalence on corresponding smooth deformation spaces was used in the main proof of \cite{GMTW}.  

\begin{corollary}\label{relcongiso}
Let $X,Y\in \smsset$ and let $f:X\to Y$ be a morphism. Suppose that for all closed subsets $A\subset U$, $f$ induces an isomorphism on relative concordance classes 
$$f_*:X[U;A]\to Y[U;A].$$
Then $f$ induces a weak equivalence of spaces
$$\csp(f):\csp(X)\to \csp(Y).$$
\end{corollary}
\proof
Consider the two cases $A=\emptyset$, $U=\Delta^0_e$ and $A=\partial \Delta^n_{st}$, $U=\Delta^n_e$ and apply Corollary \ref{concequiv}.
\endofproof

Let $\iota_d:\cartsp\into \FEmb_d$ be the faithful functor that sends a cartesian space $U$ to the projection $U\times \RR^d\to U$ and sends a smooth map $f:U\to V$ to the fiberwise open embedding $f\times {\rm id}_{\RR^d}:U\times \RR^d\to V\times \RR^d$. Then restriction along $\iota_d$ gives a functor 
\begin{equation}\label{iotafunctor}
\iota_d^*:\Struct_d\to \smsset\;,
\end{equation}
which simply evaluates a simplicial presheaf on objects of the form $p_U:U\times \RR^d\to U$ and morphisms that restrict to identities on the fibers. 

\begin{lemma}\label{etaequiv}
Fix $d\geq 0$. Let $\Ei_d$ be as Definition \ref{fiberwiseconc} and let $X\in \Struct_d$. Applying $\csp$ to the map obtained by evaluating $\iota^*_d$ on the derived unit $\eta_X:X\to \Ei_dX$ gives rise to a weak equivalence of spaces
$$\csp\iota_d^*(\eta_X):\csp\iota^*_d(X)\to \csp\iota_d^*\Ei_d(X)\;.$$
\end{lemma}
\proof
First, observe that $\iota_d^*$ preserves all weak equivalences of simplicial presheaves, since these are defined objectwise. Moreover, $\iota^*_d$ manifestly preserves all monomorphisms. Since $\iota_d^*$ admits a right adjoint, by Kan extension, it follows that $\iota_d^*$ is a left Quillen functor in the injective model structure. Since $\iota_d^*$ is already derived in the injective model structure, $\iota^*_d$ is homotopy cocontinuous. The functor $\csp$ is also homotopy cocontinuous by definition. Hence it suffices to prove the claim in the case where $X=M\to V$ is representable. 

By definition of $\Ei_d$, the map $\eta_{M\to V}$ evaluates to a natural map (natural on $\cartsp$)
\begin{equation}\label{mapfembs}
\iota_d^*(\eta_{M\to V})(U):\FEmb_d(\RR^d\times U\to U,M\to V)\to \frakFEmb_d(\RR^d\times U\to U,M\to V)\;.
\end{equation}
Observe that both sides satisfy homotopy descent on $\cartsp$, since the property of being a smooth family of embeddings is local with respect to the parametrizing family direction \footnote{It is clearly not a sheaf in the fiber direction, as being an embedding is not a local property.}. Therefore, we are in a position to apply Corollary \ref{concequiv}, so that it suffices to show that the map \eqref{mapfembs} induces an isomorphism on relative concordance classes when $U=\Delta^n_e$ and the germ of the boundary is sent to a fixed point. Since both sides have the same vertices, the map is surjective on concordance classes. To prove injectivity, we must show that two fiberwise embeddings $f,g:\RR^d\times \Delta^n_e\to M$ are concordant by a map $H:\RR^d\times \Delta^n_e\times \Delta^1_e\to M$ covering a map $h:\Delta^n_e\times \Delta^1_e\to V$:
$$
\xymatrix{
\RR^d\times \Delta^n_e \times \Delta^1_e\ar[r]^-{H}\ar[d]^-{p_{2,3}} & M\ar[d]
\\
\Delta^n_e\times \Delta^1_e\ar[r]^-{h} & V\;,
}
$$
such that both $H$ and $h$ are constant on an $\epsilon$-neighborhood of the boundary $\partial \Delta_{st,\epsilon}$, if they are concordant by a map $H':\RR^d\times \Delta^n_e\times \Delta^1_e\times \Delta^1_e\to M$ covering a map $h':U\times \Delta^1_e\to V$:
$$
\xymatrix{
\RR^d\times \Delta^n_e\times \Delta^1_e\times \Delta^1_e\ar[r]^-{H'}\ar[d]^-{p_{2,3}} & M\ar[d]
\\
\Delta^n_e\times \Delta^1_e\ar[r]^-{h'} & V\;,
}
$$
that is constant on some $\epsilon'$-neighborhood of the boundary. To prove this, set $\epsilon=\epsilon'$. Let $d:\Delta^1_e\to \Delta^1_e\times \Delta^1_e$ denote the diagonal map. Then we have a commutative diagram 
$$
\xymatrix{
\RR^d\times \Delta^n_e\times \Delta^1_e\ar[r]^-{\id\times d}\ar[d]^-{p_{2,3}} & \RR^d\times \Delta^n_e\times \Delta^1_e\times \Delta^1_e\ar[r]^-{H'}\ar[d]^-{p_{2,3}} & M\ar[d]
\\
\Delta^n_e\times \Delta_e^1\ar[r]^-{\id} & \Delta^n_e\times \Delta^1_e\ar[r]^-{h'} & V\;,
}
$$
and setting $H=(\id\times d)^*H'$ and $h'=h$ yields the desired concordance. 
\endofproof

The next corollary gives us a convenient criteria for showing that a morphism of presheaves on $\FEmb_d$ induces a weak equivalence on deformation spaces. 

\begin{corollary}\label{critwe}
Fix $d\geq 0$. Recall the restriction functor $\iota_d^*$ from \eqref{iotafunctor}. Recall also the isotopification functor $\E_d$ in Definition \ref{fiberwiseconc} and the deformation space functor $\D_d$ in Definition \ref{dfunctor}. Let $f:X\to Y$ be a morphism in $\Struct_d$ such that for all closed subsets $A\subset U$, the  map 
$$\iota_d^*(f):\iota_d^*(X)\to \iota^*_d(Y)\in \smsset$$
induces an isomorphism on relative concordance classes
$$\iota_d^*(f)_*:\iota^*_d(X)[U;A]\overset{\cong}{\to} \iota^*_d(Y)[U;A]\;.$$ 
 Then the map
$$\D_d\E_d (f):\D_d\E_d (X)\to \D_d\E_d (Y)$$
is a weak equivalence in $\mathscr{S}^{{\rm O}(d)}$. 
\end{corollary}
\proof
Let us recall that $\mathscr{S}^{{\rm O}(d)}=\PSh_{(\infty,1)}({\rm B}{\rm O}(d))$, by definition. Here ${\rm B}{\rm O}(d)$ is the simplicial category with a single object and ${\rm sing}({\rm O}(d))$ as the simplicial set of morphisms. Let $\ast\to {\rm B}{\rm O}(d)$ be the canonical map and let $\gamma:\mathscr{S}^{{\rm O}(d)}=\PSh_{(\infty,1)}({\rm B}{\rm O}(d))\to \mathscr{S}$ denote the induced restriction functor. Consider the corresponding diagram 
$$
\xymatrix{
\smequi\ar[r]^-{\csp_{\mathscr{S}^{{\rm O}(d)}}}\ar[d]^{\gamma} & \mathscr{S}^{{\rm O}(d)}\ar[d]^-{\gamma}
\\
\smsset\ar[r]^-{\csp_{\mathscr{S}}} & \mathscr{S}\;.
}
$$
Since $\gamma$ is homotopy cocontinuous (it admits a right homotopy adjoint by the derived Kan extension), it follows at once that the diagram is commutative. Explicitly, we have
$$
\csp_{\mathscr{S}}(\gamma(X))=\hocolim_{[n]\in \Delta}\gamma(X(\Delta^n_e))\simeq \gamma\left(\hocolim_{[n]\in \Delta^{\rm op}}X(\Delta^n_e)\right)=\gamma(\csp_{\mathscr{S}^{{\rm O}(d)}}(X))\;.
$$
By definition, a morphism $f:X\to Y$ of equivariant spaces $\mathscr{S}^{{\rm O}(d)}$ is a weak equivalence (i.e., an objectwise weak equivalence) if and only if the corresponding map $\gamma(f):\gamma(X)\to \gamma(Y)$ is a weak equivalence. Therefore, to prove that $\D_d \E_d(f):\csp \R_d\E_d(X)=\D_d \E_d(X)\to \D_d\E_d (Y)=\csp\R_d\E_d(Y)$ is a weak equivalence, it suffices to show that the map $\gamma \D_d (f)$ (forgetting the ${\rm O}(d)$-action) is a weak equivalence in $\mathscr{S}$. 

Now the functor $\iota^*_d$ evaluates on just the fiberwise identity maps $\RR^d\times U\to \RR^d\times U$. By definition of $\E_d$ and $\Ei_d$ (Definition \ref{fiberwiseconc}), the composition $\iota_d^*\Ei_d$ is precisely $\gamma \R_d\E_d$. Since $\csp$ commutes with $\gamma$, we have a commutative diagram 
$$
\xymatrix{
\csp_{\mathscr{S}}(\iota_d^*(X))\ar[rr]^-{\csp_{\mathscr{S}}(\iota_d^*(\eta_X))}\ar[d]_{\csp_{\mathscr{S}}(\iota_d^*(f))} && \csp_{\mathscr{S}}(\iota_d^*(\Ei_d(X))) \ar[r]^{\simeq}\ar[d]_{\csp_{\mathscr{S}}(\iota_d^*\Ei_d(f))}  &  \gamma \csp_{\mathscr{S}^{{\rm O}(d)}}\R_d\E_d(X)\ar[d]_{\gamma\D_d\E_d(f)}
\\
\csp_{\mathscr{S}}(\iota_d^*(Y))\ar[rr]_-{\csp_{\mathscr{S}}}(\iota_d^*(\eta_Y)) && \csp_{\mathscr{S}}(\iota_d^*(\Ei_d(Y))) \ar[r]^{\simeq} &  \gamma \csp_{\mathscr{S}^{{\rm O}(d)}}\R_d\E_d(Y)
}
$$
By Lemma \ref{etaequiv}, the horizontal maps are weak equivalences. By Corollary \ref{relcongiso}, the left vertical map is a weak equivalence. By 2-out-of-3, so is the right. 
\endofproof

\section{Non-topological invertible field theories and smooth stable homotopy theory}\label{stablehom}

\subsection{Brown–Comenetz and Anderson duality for sheaves of spectra}\label{invertible}

Recall the Brown--Comenetz dual \cite{BrownComenetz} of the sphere spectrum, defined as the spectrum $I_{\QQ/\ZZ}$ which represents the cohomology theory 
$$X\mapsto \hom(\pi_{-*}^{st}(X),\QQ/\ZZ).$$
That this indeed defines a cohomology theory follows from divisibility of $\QQ/\ZZ$.
Similarly, one has a spectrum $I_{\QQ}\simeq H\QQ$ representing cohomology with $\QQ$-coefficients $X\mapsto \hom(\pi_{-*}^{st}(X),\QQ)\cong H^*(X;\QQ)$. 
The \emph{Anderson dual} $I_{\ZZ}$ is defined as the homotopy fiber of the map $I_{\QQ}\simeq H\QQ\to I_{\QQ/\ZZ}$ (see Anderson \cite{Anderson} for the original construction, Hopkins--Singer \cite{HopkinsSinger} for the construction used here). 

Our first goal is to establish an analogue of Brown--Comenetz duality for sheaves of spectra on the site of cartesian spaces. In Definition \ref{smanddu}, we use the following variant of the Brown representability theorem. The proof is a trivial adaptation of Patchkoria--Pstragowski \cite[Proposition 2.15]{PatchkoriaPstragowski} to the setting of model categories. 
\begin{proposition}
\label{Brown}
Suppose $\cC$ is a combinatorial stable model category and let $E:\cC^{\rm op}\to \Ab$ be a homological functor that takes arbitrary direct sums in $\cC$ to products. Then $E$ is representable in the homotopy category $\Ho({\sf C})$. 

\end{proposition}
\proof
The proof proceeds exactly as in \cite[Proposition 2.15]{PatchkoriaPstragowski}. The crucial observation is that when $\cC$ is combinatorial and stable, it is a left Bousfield localization of the projective model structure on simplicial presheaves with values in spectra on some small category $K$ (Dugger \cite{Dugger2}). The homotopy category of $\PSh(K;\Sp)_{\rm proj}$ is compactly generated and triangulated. Therefore, we can apply the results of \cite{Neeman} to the composition 
$$\PSh(K;\Sp)_{\rm proj}^{\op}\to \cC^{\rm op}\overset{E}{\to} \Ab,$$
where the first functor is the homotopy localization functor (fibrant replacement in the local model structure). This gives an object $X$ in $\PSh(K;\Sp)$ representing $E$ in the homotopy category. But this object is manifestly a local object, since $E$ factors through $\cC^{\rm op}$. Hence $E$ is representable by an object in $\cC$. 
\endofproof

\def\smSp{{\rm C}^{\infty}\Sp}

As an immediate consequence of the Brown representability theorem Proposition \ref{Brown}, the following definition makes sense.

\begin{definition}\label{smanddu}
Let $\CC^{\times}$ be the the units of $\CC$ with its smooth structure. 
We define the \emph{smooth Brown--Comenetz dual} (of the sphere spectrum) with values in $\CC^{\times}$ as the sheaf of spectra $I_{\CC^{\times}}$ associated to the cohomology theory:
$$X\mapsto {\rm hom}_{\sh(\cartsp;\Ab)}(\tilde{\pi}_*(X), \CC^{\times}).$$
Here $X\in \smSp$, $\tilde{\pi}_*$ are its sheaves of stable homotopy groups, and $\CC^{\times}$ is embedded as a smooth sheaf of abelian groups in the usual way. 
The hom on the right is taken in sheaves of abelian groups
\end{definition}

\begin{notation}
Recall the functors $\csp_{\Sp}$, $\gamma_{\Sp}$, $\delta_{\Sp}$, $\mathcal{C}_{\Sp}$ defined in Definition \ref{constglob}, \ref{constsheaf}, and \ref{struts} for ${\sf V}=\Sp$. Throughout this section, we will use the simplified notation $\csp:=\csp_{\Sp}$, $\gamma:=\gamma_{\Sp}$, $\delta:=\delta_{\Sp}$, $\mathcal{C}:=\mathcal{C}_{\Sp}$.
\end{notation}

\begin{proposition}\label{shapeglobal}
The functors 
$$\gamma:\smSp\to \Sp \qquad \csp:\smSp\to \Sp \qquad \delta:\Sp\to\smSp $$
 induce a triple homotopy adjunction $\csp\dashv \delta \dashv \gamma$.
\end{proposition}

Suppose we take a divisible abelian group $A$. We could first form the locally constant sheaf $\delta(A)\in \sh(\cartsp;\Ab)$, and then apply the construction of Definition \ref{smanddu}, with $\delta(A)$ replacing $\CC^{\times}$, to obtain a sheaf of spectra $I_{\delta(A)}$. On the other hand, we could first take the usual Brown--Comenetz dual $I_A$ and then form the constant sheaf of spectra $\delta(I_A)$. The next proposition shows that the two construction agree.

\begin{proposition}
Let $A$ be a discrete divisible group. Let $\delta(I_A)$ be the locally constant sheaf of spectra associated to the Brown--Comenetz dual $I_A$. Let $I_{\delta(A)}$ be the sheaf of spectra associated to the locally constant sheaf of groups $\delta(A)$.
Then, we have an equivalence of sheaves of spectra 
$$\iota:\delta(I_A)\to I_{\delta(A)}$$
\end{proposition}
\proof
We first define a map $\iota:\delta(I_A)\to I_{\delta(A)}$. From the universal property of $I_{\delta(A)}$, we have an isomorphism
\begin{equation}\label{mapofspec}
\pi_0\map(\delta(I_A),I_{\delta(A)})\cong \hom_{\sh(\cartsp;\Ab)}(\tilde \pi_0(\delta(I_A)),\delta(A)).
\end{equation}
Since $\delta(I_A)$ is locally constant, we have an isomorphism of sheaves of abelian groups $\tilde \pi_0(\delta(I_A))\cong \delta(\pi_0(I_A))$. By the adjunction $\delta\dashv \gamma$ on sheaves of abelian groups, and the relation $\gamma\delta=\id$, the right side of \eqref{mapofspec} is isomorphic to  
\begin{equation}\label{mapofspec2}
\hom_{\sh(\cartsp;\Ab)}(\delta(\pi_0(I_A)),\delta(A))\cong \hom_{\Ab}(\pi_0(I_A),\gamma\delta A)= \hom_{\Ab}(\pi_0(I_A), A).
\end{equation}
By the universal property of $I_A$, we have an isomorphism 
\begin{equation}\label{mapofspec3}
\hom_{\Ab}(\pi_0(I_A), A)\cong \pi_0\map(I_A,I_A).
\end{equation}
Then passing through the isomorphisms \eqref{mapofspec3},\eqref{mapofspec2} and \eqref{mapofspec}, we have a homotopy class of maps $\iota:\delta(I_A)\to I_{\delta(A)}$ that corresponds to the class of the identity map $I_A\to I_A$. We claim that any representative $\iota$ of this class is a local equivalence. To this end, it suffices to prove that for contractible $U\in \cartsp$, the $U$ component of $\iota$ induces an isomorphism on stable homotopy groups. 

Now $\pi_*(\iota_{U})$ can be identified with the map
\begin{eqnarray}\label{nattr}
\pi_0 \map(U\otimes \mathbf{S}, \Sigma^*\delta(I_A))\cong \pi_*(\delta(I_A)(U))&\overset{\pi_*(\iota_{U})}{\longrightarrow}& \pi_*(I_{\delta(A)}(U)) \nonumber
\\
&\cong & \pi_0\map(U\otimes \mathbf{S},\Sigma^*I_{\delta(A)}) \nonumber
\\
&\cong &\hom_{\sh(\cartsp;\Ab)}(\tilde{\pi}_*(U), \delta(A))\;,
\end{eqnarray}
where the last isomorphism is obtained via the universal property of the Anderson dual. We first prove that both the source an target of $\pi_*(\iota_U)$ are isomorphic to $\hom_{\Ab}(\pi_*(\mathbf{S}),A)$. Since $U$ is contractible, the homotopy adjunction $\csp\dashv \delta$, along with the universal property of $I_A$, implies also that we have an isomorphism 
\begin{equation}\label{iso1}\pi_0\map(U\otimes \mathbf{S},\Sigma^*\delta(I_A))\cong \pi_0\map(\csp(U)\otimes \mathbf{S},\Sigma^*I_{A}) \cong \hom_{\Ab}( \pi_*(\csp(U)),A)\cong \hom_{\Ab}( \pi_*(\mathbf{S}),A).
\end{equation}
 
Since $U$ is a presheaf of sets, we have $\tilde\pi_*(U)\cong \tilde \ZZ(U)\otimes \pi_*(\mathbf{S})$, where $\tilde \ZZ$ denotes the functor that takes the sheafification of the free abelian group. Computing the right side of \eqref{nattr}, we obtain
\begin{eqnarray}\label{iso2}
\hom_{\sh(\cartsp;\Ab)}(\tilde{\pi}_*(U), \delta(A)) &\cong & \hom_{\sh(\cartsp;\Ab)}(\tilde \ZZ(U)\otimes \pi_*(\mathbf{S}), \delta(A)) \nonumber
\\
&\cong &\hom_{\sh(\cartsp;\Ab)}(\pi_*(\mathbf{S}), \hom_{\sh(\cartsp;\Ab)}(\tilde \ZZ(U),\delta(A))) \nonumber
\\
&\cong & \hom_{\sh(\cartsp;\Ab)}(\pi_*(\mathbf{S}), \hom_{\sh(\cartsp)}(U,\delta(A))) \nonumber
\\
&\cong & \hom_{\Ab}(\pi_*(\mathbf{S}), A).
\end{eqnarray}
By 2-put-of-3 for isomorphisms, it suffices to prove that pre and post composing $\iota_{U}$ with the isomorphisms \eqref{iso1} and \eqref{iso2} gives an isomorphism $\hom_{\Ab}(\pi_*(\mathbf{S}),A)\overset{\cong}{\to} \hom_{\Ab}(\pi_*(\mathbf{S}),A)$. But this map is precisely the map induced the morphism $I_A\to I_A$ (via the universal property) that corresponds to $\iota$ under the isomorphisms \eqref{mapofspec},\eqref{mapofspec2}, and \eqref{mapofspec3}. Since this corresponding map is identity, by definition, the map is an isomorphism, as claimed.
\endofproof

For a sheaf of spectra $X\in \smSp$, we denote the powering with an object $Y\in\smsset$ as $X^Y$. This is the sheaf of spectra whose levels are the internal hom $\hom(Y,X_n)$ in $\smsset$. We denote the corresponding tensoring of $X\in \smSp$ by $Y\in \smsset$ by $X\wedge Y_+$, which levelwise smashes with $Y_+$. 

\begin{proposition}\label{anderson.value}
The value of the sheaf of spectra $I_{\CC^{\times}}$ on a cartesian space $U$ is a Brown--Comenetz dual of the sphere spectrum, 
with values in the group $C^{\infty}(U;\CC^{\times})$. 
That is, $I_{\CC^{\times}}(U)$ represents the cohomology theory (on spaces)
$$X\mapsto {\rm hom}_{\Ab}( \pi_*(X),C^{\infty}(U,\CC^{\times})).$$
\end{proposition}

\proof
Observe that from the enriched Yoneda lemma, it follows at once that for a representable $U\in \smsset$ and a presheaf of spectra $X\in \smSp$, we have $
\gamma(X^U)=X^U(\RR^0)\cong X(U)$.

Let $\tilde{\ZZ}(U)$ denote the sheafification of the presheaf given by taking the free abelian group of the representable $U$ objectwise. By the universal property of $I_{\CC^{\times}}$, the adjunction $\delta\dashv \gamma$, and the commutativity $\tilde{\pi}_*(\delta(X))\cong \delta(\pi_*(X))$, we have the string of isomorphisms
\begin{eqnarray*}
\pi_0\map(X, \Sigma^nI_{\CC^{\times}}(U)) &\cong & \pi_0\map(X, \gamma(\Sigma^nI_{\CC^{\times}}^{U}))\\
&\cong & \pi_0\map(\delta(X)\wedge U_+, \Sigma^nI_{\CC^{\times}})\\
&\cong & {\rm hom}_{\sh(\cartsp;\Ab)}(\tilde{\pi}_n(\delta(X))\otimes \widetilde{\ZZ}(U),\CC^{\times})\\
&\cong & {\rm hom}_{\sh(\cartsp;\Ab)}(\delta(\pi_n(X)), \mathcal{H}{\rm om}(\widetilde{\ZZ}(U),\CC^{\times}))\\
&\cong & {\rm hom}_{\Ab}(\pi_n(X), \gamma\hom(\widetilde{\ZZ}(U),\CC^{\times}))\\
&\cong & {\rm hom}_{\Ab}(\pi_n(X),C^{\infty}(U,\CC^{\times})).
\end{eqnarray*}
Here $\mathcal{H}{\rm om}$ denotes the internal hom in sheaves of abelian groups.
\endofproof

\begin{remark}
Observe that since $\CC$ is also a divisible group, we can define a smooth variant $I_{\CC}$ completely analogously to $I_{\CC^{\times}}$. Moreover, Proposition \ref{anderson.value} still holds with $\CC$ replacing $\CC^{\times}$. 
\end{remark}

\begin{remark}
  We will use the notation $\ZZ(1):=2\pi i \ZZ$ in defining the smooth Anderson dual. This, of course, is motivated by the exact sequence
    $$\ZZ(1)\to \CC \overset{\rm exp}{\longrightarrow} \CC^{\times}.$$
    The next proposition holds equally well with the spectrum $I_{\ZZ(1)}$, defined as the homotopy fiber of the exponential map $I_{(\CC)^{\delta}}\to I_{(\CC^{\times})^{\delta}}$, replacing $I_{\ZZ}$ (the proof is verbatim the same with $\CC^{\delta}$ replacing $\QQ$ and $(\CC^{\times})^{\delta}$ replacing $\QQ/\ZZ$). Here the superscript $\delta$ indicates that we are taking the discrete topology. The spectrum $I_{\ZZ(1)}$ is used by Freed–Hopkins in \cite{FreedHopkins} and is equivalent to the Anderson dual $I_{\ZZ}$ (defined in the first paragraph of this section).
   Following the next proposition, we will always use either $I_{\ZZ(1)}$ or $I_{\delta(\ZZ(1))}$ for the Anderson dual or smooth Anderson dual in order to stay consistent with the notation of Freed–Hopkins \cite{FreedHopkins}.
\end{remark}

\begin{definition}
  We define the \emph{smooth Anderson dual of the sphere} $I_{\delta(\ZZ(1))}$ as the homotopy fiber of the map ${\rm exp}:I_{\CC}\to I_{\CC^{\times}}$, induced by the exponential map ${\rm exp}:\CC\to \CC^{\times}$.
  \end{definition}

\begin{proposition}\label{anderson}
The sheaf of spectra $I_{\delta(\ZZ(1))}$ is equivalent to $\delta(I_{\ZZ})$, where $I_{\ZZ}$ is the usual Anderson dual of the sphere. 
\end{proposition} 
\proof 
Consider the homotopy pullback diagram of sheaves of spectra
\begin{equation}\label{cartsq}
\xymatrix{
P\ar[r] \ar[d] &  I_{\delta(\QQ/\ZZ)}\ar[d]\\
 I_{\CC}\ar[r] & I_{\CC^{\times}},
}
\end{equation}
where the vertical map on the right is induced by the inclusion $i:\QQ/\ZZ\into \CC^{\times}$. 
From the universal property of the homotopy pullback, there is a map $j:I_{\delta(\QQ)}\to P$, which is unique up to a contractible choice, induced by the obvious cone over the diagram \eqref{cartsq}. 
We will show that $j$ is a local equivalence.
Let $U\in \cartsp$ be a contractible open subset of $\RR^n$, for example an open ball. The long exact sequence on cohomology induced by the homotopy pullback \eqref{cartsq}, along with Proposition \ref{anderson.value}, implies that the $P$ cohomology of $U$ sits in an exact sequence:
\begin{center}
\adjustbox{scale=.85}{
\begin{tikzpicture}[descr/.style={fill=white,inner sep=1.5pt}]
        \matrix (m) [
            matrix of math nodes,
            row sep=1em,
            column sep=2.5em,
            text height=1.5ex, text depth=0.25ex
        ]
        {  
        && \hom_{\Ab}(\pi_{*-1}\mathbf{S},C^{\infty}(U;\CC^{\times})) &&\\
        & P^*(U) & \hom_{\Ab}(\pi_*\mathbf{S},\QQ/\ZZ)\oplus \hom_{\Ab}(\pi_*\mathbf{S},C^{\infty}(U;\CC)) &&  \hom_{\Ab}(\pi_*\mathbf{S},C^{\infty}(U;\CC^{\times})) \\
         &  & P^{*+1}(U)  \;.&&\\
        };

        \path[overlay,->, font=\scriptsize,>=latex]
        (m-1-3) edge[out=355,in=175] node[descr,yshift=0.3ex] {$\beta$} (m-2-2)
        (m-2-2) edge (m-2-3)
        (m-2-3) edge node[above] {${\rm exp}-i$} (m-2-5)
        (m-2-5) edge[out=355,in=175] node[descr,yshift=0.3ex] {$\beta$} (m-3-3);
\end{tikzpicture}
}
\end{center}
Now we claim the map 
$${\rm exp}:\hom_{\Ab}(\pi_{*}\mathbf{S},C^{\infty}(U;\CC))\to \hom_{\Ab}(\pi_{*}\mathbf{S},C^{\infty}(U;\CC^{\times}))$$ 
is surjective. For $\phi:\pi_*\mathbf{S}\to C^{\infty}(U;\CC^\times)$, we define $\phi'\in \pi_{*}\mathbf{S}\to C^{\infty}(U;\CC)$ by 
$$\phi'(\alpha)=\log_{U}(\phi(\alpha)),$$
where $\log_U$ is a branch of the logarithm that is globally defined on the image of $\phi(\alpha):U\to \CC^{\times}$. This is possible since $U$ is an open ball in $\RR^n$. Then ${\rm exp}(\phi')=\phi$, proving surjectivity. This immediately implies that the map ${\rm exp}-i$ is also surjective.  

Since the kernel of $\beta$ is the image of this map, it follows that both $\beta$'s are the zero map.  This leads to a short exact sequence 
$$0\to P^*(U)\to \hom_{\Ab}(\pi_*\mathbf{S},\QQ/\ZZ)\oplus \hom_{\Ab}(\pi_*\mathbf{S},C^{\infty}(U;\CC)) \to \hom_{\Ab}(\pi_*\mathbf{S},C^{\infty}(U;\CC^{\times}))\to 0\;.$$
Then the five lemma applied to the natural transformation between the long exact sequence associated to map $j:I_{\delta(\QQ)}\to P$ and the similar exact sequence with $\hom_{\Ab}(\pi_*(\mathbf{S}),\QQ)$ replacing $P^*(U)$ implies that the map $j$ induces an isomorphism $P^*(U)\cong \hom_{\Ab}(\pi^*(\mathbf{S}),\QQ)$.

Now in any stable model category, a square is homotopy cartesian precisely when the induced map on homotopy fibers is an equivalence. The homotopy fiber of the top map in \eqref{cartsq} is $\delta(I_{\ZZ})$, since $\delta(I_{\QQ})\simeq I_{\delta(\QQ)}\simeq P$, $\delta(I_{\QQ/\ZZ})\simeq I_{\delta(\QQ/\ZZ)}$ and $\delta$ is homotopy continuous. The fiber of the bottom is 
$I_{\delta(\ZZ(1))}$ by definition. Thus, the induced map on fibers yields an equivalence $\delta(I_{\ZZ})\simeq I_{\delta(\ZZ(1))}$.

\endofproof

Recall the notion of a concordance between two morphisms of presheaves on $\cartsp$, Definition \ref{concordance}. Observe that we can define a concordance of presheaves of spectra the same way, using the tensoring over $\smsset$. 
More precisely, two morphisms of presheaves of spectra $f,g:X\to Y$ are said to be concordant if there is a morphism of presheaves of spectra $H:X\wedge (\Delta_e^1)_+\to Y$ making the diagram 
$$\xymatrix{
X \ar[dr]_{f}\ar[r] & X\wedge (\Delta^1_e)_+\ar[d]^-{H} &\ar[l] X\ar[dl]^{g}\\
&Y& 
}$$
where the two horizontal maps are the obvious analogous of the face inclusion, used in Definition \ref{concordance}. 
A map $f:X\to Y$ is said to be a \emph{concordance equivalence} if there is a map $g:Y\to X$ such that both $fg$ and $gf$ are concordant to the identity.

\begin{proposition}\label{shapeanderson}
The sheaf of spectra $\Icsp(I_{\CC})$ is contractible, and we have an equivalence of spectra 
$$\Icsp(I_{\CC^{\times}})\simeq \Sigma I_{\delta(\ZZ(1))}.$$
\end{proposition}
\proof First, observe that the sheaf $\CC$ is concordant to the trivial group through group homomorphisms (just take a straight-line homotopy). 
From the universal property of the smooth Brown--Comenetz dual, this concordance induces a concordance equivalence of sheaves of spectra $I_{\CC}\simeq_{con} I_0\simeq \ast$. 
By \cite[Corollary 2.4]{Pavlov}, $\csp$ sends concordance equivalences to homotopy equivalences and we have that $\Icsp = \delta\csp $ sends concordance equivalences to weak equivalences (since $\delta$ preserves all weak equivalences). We therefore have an induced equivalence of sheaves of spectra $\Icsp ( I_{\CC} )\simeq \ast$.

As is true in any stable model category, homotopy fiber sequences of spectra are also homotopy cofiber sequences. 
Since $\Icsp$ is a homotopy cocontinuous functor, it preserves these. Since $\Icsp (I_{\CC})$ is contractible, the connecting map associated to the long fiber sequences associated to $I_{\delta(\ZZ(1))}\to I_{\CC}\to  I_{\CC^{\times}}$ yields the desired equivalence. 
\endofproof

\subsection{The Madsen--Tillmann spectrum as a sheaf of spectra}

%

Fix $d\geq 0$. Let 
$$c:\PSh_{(\infty,1)}(\Gamma\times \cartsp)\to \sh_{(\infty,1)}(\Delta^{\times d}\times \Gamma\times \cartsp),$$ 
the functor obtained by restriction along the projection $\Delta^{\times d}\times \Gamma\times \cartsp\to \Gamma\times \cartsp$. The functor $c$ admits a further left adjoint $L$, by left Kan extension. The functor $L$ formally inverts all morphisms in an $(\infty,d)$-category. In the case where $d=1$, the functor $L$ is the usual geometric realization of a simplicial object $\PSh_{(\infty,1)}(\Gamma\times \cartsp)$. 

The functor $c$ preserves all objectwise weak equivalences. However, the functor $L$ must be appropriately derived in order to give rise to a morphism of relative categories. In the following, we write 
$\PSh_{(\infty,1)}(\Gamma)_{\rm loc},$
for simplicial presheaves on $\Gamma$ with weak equivalences given by local weak equivalences in the left Bousfield localization at the morphisms
$$
\langle \ell\rangle\coprod_{\langle 0 \rangle} \langle \kappa\rangle\to \langle \ell+\kappa \rangle, \qquad \emptyset \to \langle 0 \rangle, \qquad \langle \ell\rangle,\langle \kappa\rangle\in \Gamma,
$$
which implement Segal's special $\Gamma$-condition (see  Grady--Pavlov \cite[Section 2]{GradyPavlov} for an exposition). Connective spectra can be obtained from $\Gamma$-spaces by the work of Bousfield--Friedlander \cite{BousfieldFriedlander}, as follows. 

Let $j:\Delta\to \Gamma$ denote the functor that sends $[n]\in \Delta$ to $[n]\sqcup\{\ast\}$ and a morphism $f:[m]\to [n]$ to the morphism of finite pointed sets $g:[n]_+\to [m]_+$ defined by $g(+)=+$ and $g(l)=i$ for all $f(i-1)<l\leq f(i)$. Restricting a $\Gamma$-space along the functor $j$ gives rise to a simplicial space. The homotopy colimit, or realization, of this simplicial diagram in spaces gives the zero-th space of a connective spectrum. We denote the realization by $B$. Hence, we have a composite functor $Bj^*:\PSh_{(\infty,1)}(\Gamma)\to \mathscr{S}.$
This functor factors as a composition 
$$\PSh_{(\infty,1)}(\Gamma)\overset{S}{\to} \Sp^{\mathbb{N}}\to {\rm Fun}(\mathbb{N},\mathscr{S}_{/\ast})\overset{ev_0}{\to} \mathscr{S}_{/\ast}\into \mathscr{S}.$$ 
The composition of the first two functors in this factorization is given by taking the left Kan extension along the inclusion $\Gamma\into \mathscr{S}_{/\ast}$ and restricting to spheres $S^n\in \mathscr{S}_{/\ast}$. It is easy to show that this composition factors through $\mathbf{S}$-modules, hence through $\Sp^{\NN}$.  
By \cite[Section 5.7]{BousfieldFriedlander}, the functor $S$ is left Quillen in the stable model structure on $\PSh_{(\infty,1)}(\Gamma)$. The functor $S$ is homotopically fully-faithful and its (homotopical) essential image is the connective spectra $\Sp^{\mathbb{N}}_{\geq 0}$. Following Notation \ref{notation}, we write $\Sp_{\geq 0}$ for the relative category of $\Gamma$-spaces, with weak equivalences given by the stable local weak equivalences.

\begin{proposition}\label{formalinv}
The functors $c$ and $L$ induce morphisms of relative categories 
\begin{align}
& c:\sh_{(\infty,1)}(\cartsp;\Sp_{\geq 0})\to \smcat_{\infty,d}
\\
& \vert \cdot \vert:\smcat_{\infty,d}\to \sh_{(\infty,1)}(\cartsp;\Sp_{\geq 0})
\end{align}
where the functor $\vert \cdot \vert:=\mathbb{L}L$ is the left derived functor of the left adjoint $L\dashv c$ in the local injective model structure on $\PSh_{(\infty,1)}(\Delta^{\times d}\times \Gamma\times \cartsp)$, defined in \cite[Section 2]{GradyPavlov} (with the {\v C}ech morphisms left out). We have an induced homotopy adjunction $\vert \cdot \vert\dashv c$.  

We remark that $\vert \cdot \vert$ is not derived in the \emph{{\v C}ech local} model structure, and we prove that it automatically descends to sheaves. 
\end{proposition}
\proof
By definition of the model structure on $\smcat_{\infty,d}$, the homotopy adjunction holds at the level of presheaves. We will need to refer to the set of morphisms in Grady--Pavlov \cite[Section 2]{GradyPavlov} given by 
\begin{equation}\label{cechmaps2}
C(\{U_a\})\otimes j(\langle \ell\rangle,{\bf m})\to U\otimes j(\langle \ell\rangle,{\bf m}),
\end{equation}
where $(\langle \ell\rangle,{\bf m})\in \Gamma\times \Delta^{\times d}$, the functor $j$ is the embedding 
$$j:\PSh_{(\infty,1)}(\Gamma\times \Delta^{\times d})\to \PSh_{(\infty,1)}(\cartsp\times \Gamma\times \Delta^{\times d}),$$
induced by the corresponding projection,
and the tensoring is that of $\PSh_{(\infty,1)}(\cartsp\times \Gamma\times \Delta^{\times d})$ over $\PSh_{(\infty,1)}(\cartsp)$. Left Bousfield localizing at the morphisms \eqref{cechmaps} is responsible for implementing sheaves.   
%

The left derived functor $\vert \cdot \vert$ commutes with the tensoring over $\PSh_{(\infty,1)}(\cartsp)$, which can be seen by the equivalence of categories
$$\PSh_{(\infty,1)}(\cartsp\times \Gamma\times \Delta^{\times d})\cong \PSh_{(\infty,1)}(\cartsp,\PSh_{(\infty,1)}(\Gamma\times \Delta^{\times d})).$$
Since multisimplices are contractible, it follows that $\vert \cdot \vert$ sends {\v C}ech local maps
$
 C(\{U_a\})\otimes (\langle \ell\rangle,{\bf m})\to U\otimes (\langle \ell\rangle,{\bf m})
$
to morphisms that are weakly equivalent to  
$
 C(\{U_a\})\otimes \langle \ell\rangle\to U\otimes \langle \ell\rangle.
$
Hence, $c$ sends local objects to local objects. Therefore the homotopy left adjoint $\vert \cdot \vert$ sends all {\v C}ech local weak equivalences to weak equivalences and the adjunction $\vert \cdot \vert\dashv c$ descends to sheaves.

\endofproof

\begin{definition}\label{madsensheaf}

Fix $d\geq 0$. 
We define the \emph{Madsen--Tillmann} spectrum ${\sf MT}$ of a geometric structure via the composite functor 
$${\sf MT}:\Struct_d\xrightarrow{\Bord^{(-)}_d} \smcat_{\infty,d}\overset{\vert \cdot \vert}{\longrightarrow} \sh_{(\infty,1)}(\cartsp,\Sp_{\geq 0}).$$ 
\end{definition}

\section{Invertible functorial field theories}\label{freedhopsec}

In this section, we apply our results to invertible field theories, proving a conjecture by Freed and Hopkins \cite{FreedHopkins}.
We begin with some preliminary set up.

Let $H_d$ be a compact Lie group and $\rho_d:H_d\to {\rm O}(d)$ be a homomorphism.
We can regard $H_d$ as a tangential structure via the induced map $\deloop H_d\to \deloop {\rm O}(d)$.
Such tangential structures can be converted to a sheaf on $\FEmb_d$ by the construction in Section \ref{geometric.struct}.
According to Freed–Hopkins \cite[Ansatz 5.26]{FreedHopkins}, a \emph{continuous} invertible $d$-dimensional extended topological field theory should be represented by a morphism of spectra
$$\Sigma^dMTH_d\to \Sigma^{d+1}I_{\ZZ(1)}.$$
This ansatz comes from the following two observations.
First, if $\mathscr{C}$ is a smooth symmetric monoidal $(\infty,d)$-category with duals, then any field theory 
$Z:\Bord_d^{\mathcal{H}_d}\to \mathscr{C}$ 
that factors through the subobject ${\rm Pic}(\mathscr{C})\into \mathscr{C}$, containing only invertible objects and arrows, induces a map out of the realization of $\Bord_d^{\mathcal{H}_d}$, i.e., the space obtained by formally inverting all objects and arrows. 

We let $\vert \cdot \vert$ be as in Definition \ref{madsensheaf}. The above observation shows that we have a commutative diagram 
$$
\xymatrix{
\mathfrak{Bord}_d^{\mathcal{H}_d}\ar[r]^-{Z}\ar[d] & \mathscr{C}
\\
\vert \mathfrak{Bord}_d^{\mathcal{H}_d}\vert \ar@{-->}[r]^-{\overline{Z}} & {\rm Pic}(\mathscr{C})\;.\ar@{^{(}->}[u]
}
$$
This already gives us a description of invertible field theories in terms of morphisms of (sheaves of) spectra. In the notation of Definition \ref{madsensheaf}, the map $\overline{Z}$ in the above diagram is equivalently a map of the form
$$
Z':\MT(\mathcal{H}_d)\to {\rm Pic}(\mathscr{C}).
$$

 The second key observation is that the Brown–Comenetz dual $I_{\CC^{\times}}$ is a natural target for field theories. Indeed, the universal property of the Brown–Comenetz dual imply that we have an isomorphism
\begin{equation} \label{partitionfun}
\pi_0\map(\vert \Bord_d^{\mathcal{H}_d}\vert,\Sigma^dI_{\CC^{\times}})\cong {\rm hom}_{\sh(\cartsp;\Ab)}(\tilde{\pi}_d\vert \Bord_d^{\mathcal{H}_d}\vert,\CC^{\times}).
\end{equation}
 
Thus, homotopy classes of maps between the sheaf of spectra $\MT(\mathcal{H}_d):=\vert \Bord_d^{\mathcal{H}_d}\vert$ and $\Sigma^dI_{\CC^{\times}}$ are precisely fiberwise partition functions parametrized over cartesian spaces. The appearance of the Anderson dual then arises by observing that taking continuous deformation classes of such field theories should amount to taking smooth deformations of both the source and target in the left side of \ref{partitionfun}. By \ref{shapeanderson}, this would give the identification
\begin{equation}\label{ansatzfh}
\pi_0\map(\csp (\MT(\mathcal{H}_d)), \csp (\Sigma^dI_{\CC^{\times}}))\cong \pi_0\map(\Sigma^dMTH_d,\Sigma^{d+1}I_{\ZZ(1)}).
\end{equation}
It appears to us that the major hurdle in proving the above comes from the fact that a smooth variant of the fully-extended bordism category did not exist in the literature. Moreover, the smooth counterpart to the Brown--Comenetz dual was not yet developed.
For these reasons, Freed–Hopkins \cite{FreedHopkins} restrict themselves to \emph{discrete} field theories.
In this case, one takes the topological spectrum $I_{(\CC^{\times})^{\delta}}$, where $\CC^{\times}$ is equipped with the discrete topology.
The main theorem of Freed--Hopkins \cite[Theorem 5.23]{FreedHopkins} identifies the torsion subgroup of the right side of \ref{ansatzfh} with deformation classes of \emph{discrete} invertible field theories.

The following lemma is immediate from the definitions.

\begin{lemma}\label{involutions}
For a fixed hyperplane reflection $\sigma$, the corresponding involutions in Proposition \ref{betahd1}, Proposition \ref{invbord2} and Definition \ref{betahd2} induce corresponding involutions on the bordism category. We denote all three involutions by $\beta$.

\end{lemma}

\subsection{Comparing the topological and geometric bordism categories}

We would like to be able to compare deformation classes of geometric field theories and topological field theories with structure given by taking deformations of the geometric structure. The following proposition allows us to ``commute" smooth deformations inside the bordism category as deformations of the geometric structure.

\begin{proposition}\label{w.e.bords}
Fix $d\geq 0$. Recall the functors $\Icsp_{\mathscr{S}}$ (Definition \ref{struts}), $\D_d$ (Definition \ref{dfunctor}) and $\lSec_d$ (Definition \ref{righadjbe}). For all $\mathcal{S}\in \Struct_d$, there is a zig-zag of equivalences
$$\Icsp_{\mathscr{S}}(\mathfrak{Bord}_d^{\mathcal{S}})\simeq \mathfrak{Bord}_d^{\lSec_d\D_d(\mathcal{S})}$$
in $\smcatdual_{(\infty,d)}$.
\end{proposition}
\proof
Recall the homotopy left adjoint $\E_d$ in Remark \ref{enrichall}. By definition, 
$\mathfrak{Bord}_d^{\mathcal{S}}=\mathfrak{Bord}_d^{\E_d(\mathcal{S})}$. To reduce complicated notation, we use $\mathcal{S}=\E_s(\mathcal{S})$ below. For a smooth equivariant space $X\in \smequi$, let $X_n$ denote the smooth equivariant space given by 
$$X_n(U)=X(\Delta^n_e\times U),$$
 By definition of the concordification functor $\mathcal{C}_{\mathscr{S}^{{\rm O}(d)}}$, we have 
$$\mathcal{C}_{\mathscr{S}^{{\rm O}(d)}}(X)=\hocolim_{[n]\in \Delta^{\op}}X_n.$$
Moreover, $\mathcal{C}_{\mathscr{S}^{{\rm O}(d)}}$ computes $\Icsp_{\mathscr{S}^{{\rm O}(d)}}(\mathcal{S})$ by Proposition \ref{descent}.

 By \cite[Theorem 6.0.2]{GradyPavlov}, the functor 
$$\mathfrak{Bord}_d:\frakStruct_d\to \fraksmcat$$ 
preserves all weak equivalences. By the Quillen equivalence $\R_d\dashv \Sec_d$, the derived unit of the adjunction yields a weak equivalence 
$$\mathfrak{Bord}_d^{\mathcal{S}}\overset{\simeq}{\to} \mathfrak{Bord}_d^{\Sec_d\R_d\mathcal{S}}.$$
Since $\Icsp$ preserves all weak equivalences, it suffices to compute $\Icsp$ applied to the bordism category on the right above. 

Unraveling the definition of the bordism category in Grady--Pavlov \cite[Definition 4.5.2]{GradyPavlov}, we have 
$$
\mathfrak{Bord}_d^{\Sec_dX}(\Delta^n_e\times U)=\mathfrak{Bord}_d^{\Sec_d(X_n)}(U).
$$
Fix an arbitrary $U\in \cartsp$. By Proposition \ref{descent}, we can compute $\Icsp$ via the concordification functor $\mathcal{C}_{\sf V}$, with ${\sf V}=\mathscr{P}{\rm Sh}_{(\infty,1)}(\Delta^{\times d}\times \Gamma)$. Using the definition of $\mathcal{C}_{\sf V}$ and the discussion in the preceeding paragraph, we have a chain of weak equivalences
\begin{align*}
\Icsp\left(\mathfrak{Bord}_d^{\Sec_d\R_d\mathcal{S}}\right)(U) &=\hocolim_{[n]\in \Delta^{\op}}\mathfrak{Bord}_d^{\Sec_d\R_d\mathcal{S}}(\Delta^n_{e}\times U) 
\\
&\simeq \hocolim_{[n]\in \Delta^{\op}}\mathfrak{Bord}_d^{\Sec_d(\R_d\mathcal{S})_n}(U)
\\
&\simeq \mathfrak{Bord}_d^{\hocolim_{[n]\in \Delta^{\op}} \Sec_d(\R_d\mathcal{S})_n}(U)
\\
&\simeq \mathfrak{Bord}_d^{\Sec_d \Icsp_{\mathscr{S}^{{\rm O}(d)}} \R_d\mathcal{S}}(U)
\\
&\simeq \mathfrak{Bord}_d^{\lSec_d \D_d(\mathcal{S})}(U).
\end{align*}
Going from the second to the third line, we have used that $\mathfrak{Bord}_d$ is homotopy cocontinuous in its argument, by the main theorem of Grady--Pavlov \cite{GradyPavlov}. Going from the third to the fourth line, we have used that $\Sec_d$, being the right derived functor in a Quillen equivalence, commutes with homotopy colimits, along with the observations made in the first paragraph and the relation $\Icsp_{\mathscr{S}^{{\rm O}(d)}}=\delta_{\mathscr{S}^{{\rm O}(d)}}\csp_{\mathscr{S}^{{\rm O}(d)}}$, $\lSec_d=\Sec_d\delta_{\mathscr{S}^{{\rm O}(d)}}$, and $\D_d=\csp_{\mathscr{S}^{{\rm O}(d)}}\R_d$ (Proposition \ref{adjointcalc}). Since $U$ was arbitrary, this completes the proof.
\endofproof


We now want to compare the deformation space of the two geometric structures $\mathcal{H}^{\rm fl}_d$ and $\mathcal{H}_d^{\nabla}$, given by fiberwise flat $H_d$-structures and fiberwise differential $H_d$-structures. As a motivation for this comparison, we recall Example \ref{flatbunexample}. 

There is a canonical map
$$u:\mathcal{H}_d^{\rm fl}\to \mathcal{H}_d^{\nabla},$$ 
which sends a fiberwise flat $H_d$-structure to the underlying $H_d$-structure and flat connection. More precisely, we have a canonical map $(\RR^d\to \RR^0)\to \mathcal{H}_d^{\nabla}$ that picks out the trivial $H_d$-bundle, the standard metric on $\RR^d$, and the trivial connection $\nabla=d$. This map induces a map out of the homotopy quotient 
$$u:\mathcal{H}^{\rm fl}_d=\RR^d/\!/H_d\to \mathcal{H}_d^{\nabla},$$
since $H_d$ acts on $\RR^d$ by rotations, which induce metric and connection preserving automorphisms of the trivial bundle on $\RR^d$ by differentiation. Similarly, we have a canonical map 
$$u':\mathcal{H}^{\rm fl}_d\to \mathcal{H}_d\;,$$
which is obtained by further post-composing $u$ with the map $\mathcal{H}_d^{\nabla}\to \mathcal{H}_d$ that forgets connections. 

For a submersion with $d$-dimensional fibers $p:M\to U$, we denote a vertex in $\mathcal{H}^{\nabla}_d(p:M\to U)$ by a quadruple $(P,g,\nabla,\phi)$, where
\begin{itemize}
\item $P\to M$ is an $H_d$-principal bundle.
\item $g$ is a fiberwise metric on $p:M\to U$.
\item $\nabla$ is a fiberwise connection on the vector bundle $E\to M$, associated to the representation $\rho_d:H_d\to \GL(d)$. 
\item $\phi:(E,\nabla)\to (\tau_pM,\nabla_g)$ is an isomorphism of vector bundles with connection. 
\end{itemize}
This notation is justified by Proposition \ref{diffhddata}.

\begin{proposition}\label{shapestr}
Fix $d\geq 0$ and recall the functor $\D_d$ in Definition \ref{dfunctor} and the isotopification functor $\E_d$ in \ref{fiberwiseconc}. The maps $\D_d\E_d(u)$ and $\D_d\E_d(u')$ are weak equivalences.
\end{proposition}
\proof
The idea of the proof is to show that both maps induce a bijection on relative concordance classes. For the comparison between fiberwise differential $H_d$-structures and fiberwise flat $H_d$-structures (the map $\D_d(u)$), the idea is to use the straight-line deformation retraction of $\RR^d$ to the origin in order to produce a concordance to a flat $H_d$-structure.

We first prove the claim for the map $u$. By Corollary \ref{critwe}, it suffices to prove that $\iota_d^*(u)$ induces an isomorphism on relative concordance classes. We will construct an inverse map 
$$v:\mathcal{H}^{\nabla}_d[\RR^d\times U;A]\to \mathcal{H}^{\rm fl}_d[\RR^d\times U;A]$$

Let $V$ be an open neighborhood of $A$ in $U$ and let $tr_V$ denote the trivial fiberwise differential $H_d$-structure on the projection $p_V:\RR^d\times V\to V$. This means that the underlying $H_d$-bundle is the trivial bundle, the connection is constant on $V$ and is the trivial connection on each fiber of the projection. The metric constant on $V$ and is the standard metric on $\RR^d$ in each fiber. 


Let $r:\Delta_e^1\times \RR^d\to \RR^d$ be a straight-line deformation retraction of $\RR^d$ to the origin, defined by $r(t,x)=tx$. Let $s$ denote the fiberwise constant metric on the projection $\RR^d\times U\to U$, equal to the standard metric on $\RR^d$ in each fiber. The map $r$ induces a function 
$$r^*:\mathcal{H}^{\nabla}_d(\RR^d\times U\to U)_0\to \mathcal{H}^{\nabla}_d(\RR^d\times \Delta^1_e\times U\to \Delta^1_e\times U)_0,$$
on the set of vertices, which we now describe. We set
$$r^*(P,\nabla,g)=((r\times {\rm id}_{U})^*P, \widetilde{\nabla}_{s+(r\times {\rm id}_U)^*(g-s)}, s +(r\times {\rm id}_{U})^*(g-s),(r\times \id_U)^*\phi),$$
where the components are defined as follows.
\begin{itemize}
\item The bundle $(r\times \id_{U})^*P$ is the usual pullback bundle, considered as a fiberwise principal $H_d$-bundle with respect to the projection $\RR^d\times \Delta^1_e\times U\to \Delta^1_e\times U$ . 
\item The fiberwise metric $s +(r\times {\rm id}_{U})^*(g-s)$ is given on vertical tangent vectors $X,Y\in \RR^d\subset (T_{(x,t,m)}\RR^d\times \Delta^1_{e}\times U)$ by 
\begin{eqnarray*}
(s +(r\times {\rm id}_{U})^*(g-s))(X,Y) &=& s(X,Y)+g_{(x,m)}(tX,tY)-s(tX,tY)
\\
&=& s(X,Y)+t^2g_{(x,m)}(X,Y)-t^2s(X,Y)\;.
\end{eqnarray*}
This is a metric, since for $0<t\leq 1$, 
$$s(X,X)+t^2g_{(x,m)}(X,X)-t^2s(X,X)\geq t^2g_{(x,m)}(X,X)\geq 0$$
\item Let $E$ be the vector bundle associated to $(r\times \id_{U})^*P$ by the representation $\rho:H_d\to \GL(d)$. The fiberwise connection $\widetilde{\nabla}_{s+(r\times {\rm id}_U)^*(g-s)}$ at each $(t,m)\in \Delta^1_e\times U$ is given by the pullback of the Levi-Civita connection 
$
\nabla_{s+t(g-s)}
$, associated to the metric $s+t(g-s)$, on the pullback bundle 
$$
\xymatrix{
E_{(t,m)}\ar[r]\ar[d] & E_m\ar[d]\ar[r]^-{\phi} & (\tau_pM)_m\ar[d]
\\
\RR^d\ar[r]^{t} & \RR^d\ar[r] & \RR^d\;,
}
$$
where the map $t$ is multiplication by $t\in \Delta^1_e\subset \RR$. Note that if $t=0$ then $E_{(t,m)}=E_0\times \RR^d$, $\nabla_{s+0(g-s)}=\nabla_s=d$, and pullback connection is the trivial connection. 
\item At each $(t,m)\in \Delta_e^1\times U$, the bundle isomorphism $(r\times \id_{U})^*\phi_{(t,m)}$ is given by the pulling back the isomorphism $\phi_m:E_m\to (\tau_pM)_m$ along the map $t:\RR^d\to \RR^d$, which multiplies by $t$. 
\end{itemize}
From the definition, it is clear that $r^*$ sends concordant structures to concordant structures (the concordance itself can be pulled back by $r^*$ using the same definition, replacing the 1-parameter $t\in \Delta^1_e$ by a pair of parameters $(t,s)\in \Delta^1_e\times \Delta^1_e$). We claim that $r^*(P,\nabla,g)$ is also constant on $V$, so that the map sends relative concordances to relative concordances. On $V$, the metric is given by $g=s$ so that $r^*(g)= s +(r\times {\rm id}_{U})^*(g-s)=s$ on $V$. When $P$ is the trivial bundle, we have that the bundle $(r\times \id_V)^*P_{(t,m)}$ is the pullback bundle
$$
\xymatrix{
(r\times \id_{V})^*P_{(t,m)}\ar[r]\ar[d] & H_d\times \RR^d\ar[d]
\\
\RR^d\ar[r]^{t} & \RR^d
}
$$
which is again the trivial bundle $P_{(t,m)}=H_d\times \RR^d$. The fiberwise connection $\widetilde{\nabla}_s$ on the associated bundle is the trivial connection, which for $t>0$ can be seen by the formula for the pullback connection
\begin{equation}\label{conntriv}
(t^*\phi^*\nabla_s)_{X}f=(t^*\phi^*\nabla_s)_{X}(t^*f)\circ (1/t)=(\phi^*\nabla_s)_{tX}f\circ (1/t) = t(1/t)df (X)=df(X)\;.
\end{equation}
For $t=0$, the pullback sections are precisely the constant sections $0^*f=f(0)$ and $0^*(\phi^*\nabla_s)$ vanishes on all such sections, by definition. This again implies that $0^*\phi^*\nabla_s=d$. 
Hence, $r^*(P,g,\nabla,\phi)=tr_{V}$ on $V$, as claimed.

We define the map $v$ on $[P,g,\nabla,\phi]\in \mathcal{H}^{\nabla}_{d}[\RR^d\times U;A]$ by the formula
$$v[P,g,\nabla,\phi]=[r_0^*(P,g,\nabla,\phi)]\in \mathcal{H}^{\rm fl}_d[\RR^d\times U;A]\;.$$
This is well defined by the previous observations. The class of $r_0^*(P,g,\nabla,\phi)$ is given by evaluating the above data at $t=0$, which is a fiberwise flat $H_d$-structure. To see this, observe that evaluation at $t=0$ gives the $U$-family of bundles on $\RR^d$, which for each $m\in U$, $(P_m)_0\times \RR^d\to \RR^d$, where $P_m$ is the fiber of the original bundle $P_m\to \RR^d$ at $0$. The metric is (fiberwise) the standard metric on $\RR^d$ and the connection is (fiberwise) the trivial connection. This $H_d$-structure is fiberwise flat, since locally on the base $U$, for each $m\in U$, a choice of point $p(m)\in (P_m)_0$ gives rise to a trivialization of the form
$$
\xymatrix{
H_d\times \RR^d\ar[rrr]^{p(m)\times \rho_d(p(m))}\ar[d] &&& (P_m)_0\times \RR^d\ar[d]
\\
\RR^d\ar[rrr]^{\rho_d(p(m))} &&& \RR^d\;,
}
$$
where $p(m)$ acts on $\RR^d$ linearly, via the representation $\rho_d:(P_m)_0\to \GL(d)$.

It follows immediately from the construction that $v$ is an inverse of $u$:

\begin{itemize}
\item $u\circ v={\rm id}$
\end{itemize}

For all quadruples $(P,g,\nabla,\phi)$ representing a class in $\mathcal{H}^{\nabla}_{d}[\RR^d\times U;A]$, the element $(u\circ v)(P,g,\nabla,\phi)$ is concordant by a relative concordance to $(P,g,\nabla,\phi)$ via $r^*(P,g,\nabla,\phi)$, constructed above. Indeed, the element $r_1^*(P,g,\nabla,\phi)$ is given by evaluating the above data at $t=1$, which is $(P,g,\nabla,\phi)$ by construction.

\begin{itemize}
\item $v\circ u={\rm id}$
\end{itemize}

This is a simple unwinding of the definitions. A fiberwise flat $H_d$-structure on $\RR^d\times U\to U$ is given by the pullback of the trivial bundle on $\RR^d$ along a fiberwise open embedding. Let $P$ be such a bundle, and suppose that this bundle is the trivial family of trivial bundles on an open neighborhood $V$ of $A$, i.e. $P$ represents a relative concordance class. By definition $v\circ u(P)$ is the bundle that over each fiber $m\in U$ is the trivial bundle $H_d\times \RR^d=(P_m)_0\times \RR^d\to \RR^d$. Hence, we have $(v\circ u)(P)=P$, for all $P$.

Now to prove the claim for the map $u'$, observe that the map $r^*$ can be modified to produce an analogous map 
$$r^*:\mathcal{H}_d(\RR^d\times U\to U)\to \mathcal{H}_d(\RR^d\times \Delta^1_e\times U\to \Delta^1_e\times U),$$
given on vertices by $r^*(P)=(r\times \id_{U})^*(P)$, described in the first bullet point above. Then the entire argument for $\mathcal{H}_d^{\nabla}$ carries through verbatim, simply by removing the data of the connection and metric.

\endofproof


We conclude this section by relating the smooth and topological Madsen--Tillmann spectra. 

\begin{theorem}\label{madsentilleq}
Fix $d\geq 0$. We have a zig-zag of weak equivalences of spectra
$$
\csp_{\Sp}({\sf MT}(\mathcal{H}^{\rm fl}_d))\simeq \csp_{\Sp}({\sf MT}(\mathcal{H}^{\nabla}_d))\simeq \csp_{\Sp}({\sf MT}(\mathcal{H}_d))\simeq \Sigma^d{\rm MT}H_d,
$$
where on the right we have the $d$-fold suspension of the usual Madsen--Tillmann spectrum.
\end{theorem}
\proof
To simplify notation, we will denote $\csp_{\Sp}$ by $\csp$ below. By Grady--Pavlov \cite[Proposition 5.0.12]{GradyPavlov2}, there is a weak equivalence $\csp({\sf MT}(\mathcal{H}^{\rm fl}_d))\simeq \Sigma^d{\rm MT}H_d$. Observe that the functor $\vert \cdot \vert$ in Definition \ref{madsensheaf} commutes with $\csp$, since $\csp$ is given by a homotopy colimit and the functor $\vert \cdot \vert$ is homotopy cocontinuous. We claim that $\vert \cdot \vert$ also sends local weak equivalences (in the duals model structure) in $\smcatdual_{\infty,d}$ to local weak equivalences. The homotopy right adjoint $c$ of $\vert \cdot \vert$ includes a spectrum as a Picard $\infty$-groupoid into all symmetric monoidal $(\infty,d)$-categories. Since Picard $\infty$-groupoids have all duals, the right adjoint sends local objects to local objects. It follows that $\vert \cdot \vert$ sends local equivalences in $\smcatdual_{\infty,d}$ to local equivalences in $\Sp_{\geq 0}$. 

Applying ${\sf MT}=\vert \mathfrak{Bord}_d^{(-)}\vert $ (Definition \ref{madsensheaf}) to the zig-zag of weak equivalence in Proposition \ref{w.e.bords}, we obtain a zig-zag of equivalences
\begin{equation}\label{mtequiv}
\Sigma^d{\rm MT}H_d\simeq \csp({\sf MT}(\mathcal{H}_d^{\rm fl}))\simeq 
\left|\Icsp\left(\mathfrak{Bord}_d^{\mathcal{H}^{\rm fl}_d}\right)\right|
\overset{\simeq}{\to}  \left|\mathfrak{Bord}_d^{\lSec_d\D_d\E_d(\mathcal{H}_d^{\rm fl})}\right|.
\end{equation}
By Proposition \ref{shapestr}, the map $\D_d\E_d(u):\D_d\E_d(\mathcal{H}_d^{\rm fl})\to \D_d\E_d(\mathcal{H}_d^{\nabla})$ is a weak equivalence in ${\rm O}(d)$-equivariant spaces. Since $\mathfrak{Bord}_d$ (by definition) sends all weak equivalences to weak equivalences and $\vert \cdot \vert$ sends all local weak equivalences in $\smcatdual_{\infty,d}$ to weak equivalences in $\Sp_{\geq 0}$, by Proposition \ref{formalinv}, it follows that the right hand side is equivalent to ${\sf MT}(\lSec_d\D_d\E_d(\mathcal{H}_d^{\nabla}))$. Passing through the equivalence \eqref{mtequiv} once again, with $\mathcal{H}^{\nabla}_d$ replacing $\mathcal{H}_d^{\rm fl}$, shows that there is a zig-zag of equivalences
$$\Sigma^d{\rm MT} H_d\simeq \csp({\sf MT}(\mathcal{H}_d^{\rm fl}))\simeq \csp({\sf MT}(\mathcal{H}_{d}^{\nabla})).$$
For the remaining equivalence we argue as above, using the equivalence $\D_d\E_d(u')$ from Proposition \ref{shapestr}. 
\endofproof

\subsection{Deformation classes of invertible field theories}

In this subsection, we prove the main theorem. The crucial ingredient has already been established in Theorem \ref{madsentilleq}, and the remaining portion of the proof is formal. We begin with a definition. 

\begin{definition}\label{deffuns}
Let $\mathscr{C},\mathscr{D}\in {\sf Cat}^{\otimes}_{\infty,d}=\PSh_{(\infty,1)}(\Gamma\times \Delta^{\times d})_{\rm loc}$ be symmetric monoidal $(\infty,d)$-categories. We define the \emph{space of symmetric monoidal functors} between $\mathscr{C}$ and $\mathscr{D}$ simply as 
$${\rm Fun}^{\otimes}(\mathscr{C},\mathscr{D}):=\map(\mathscr{C},\mathscr{D})$$
where on the right we have the mapping space in $\smcat_{\infty,d}$, provided by hammock localization.

For two \emph{smooth} symmetric monoidal $(\infty,d)$-categories $\mathscr{C},\mathscr{D}\in {\sf Cat}^{\otimes}_{\infty,d}=\PSh_{(\infty,1)}(\cartsp\times \Gamma\times \Delta^{\times d})_{\rm loc}$, we define the \emph{deformation space of symmetric monoidal functors} as the space 
$$
\mathcal{F}{\rm un}^{\otimes}(\mathscr{C},\mathscr{D}):={\rm Fun}^{\otimes}(\csp(\mathscr{C}),\csp(\mathscr{D})).
$$
Here the functor 
$\csp:\smcat\to {\sf Cat}_{\infty,d}^{\otimes}$ is defined in Definition \ref{struts}, with ${\sf V}={\sf Cat}_{\infty,d}^{\otimes}$. We again suppress the subscript on $\csp$ to simplify notation.
\end{definition}

\begin{remark}
There are several different ``deformation spaces'' one could associate to the moduli of symmetric monoidal functors between symmetric monoidal $(\infty,d)$-categories. One alternative would be to use the derived internal hom in $\smcat_{\infty,d}$, take the core and evaluate on $\langle 1\rangle\in \Gamma$, to obtain a simplicial presheaf on $\cartsp$. Then one could apply $\csp$ to the result and call this the deformation space. 

This option apprears to give something different in spirit then the deformations described in \cite{FreedHopkins}. When $\mathscr{C}=\mathfrak{Bord}_d^{\mathcal{H}_d}$ and $\mathscr{D}=\Sigma^dI_{\CC^{\times}}$, this option gives a space whose 1-simplices are morphisms in $\smcat$:
$$\mathfrak{Bord}_d^{\mathcal{H}_d}\otimes \Delta^1\to \Sigma^dI_{\CC^{\times}},$$
where the tensoring is that of $\smcat$ over $\PSh_{(\infty,1)}(\cartsp)$. Such a map gives the data of a $\Delta^1$-indexed family of natural assignments of $U$-parametrized families of bordisms to a $U$-parametrized family of values.

On the other hand, if we set $\mathscr{C}=\mathfrak{Bord}_d^{\mathcal{H}_d}$ and $\mathscr{D}=\Sigma^dI_{\CC^{\times}}$ in Definition \ref{deffuns}, a 1-simplex gives an assignment of $\Delta^1$-indexed families of bordisms to a $\Delta^1$-indexed family of values in $\Sigma^dI_{\CC^{\times}}$. This seems to be more in keeping with what a ``deformation'' of a field theory should be. 
\end{remark}

\begin{theorem}\label{hopfred}
Let $H_d$ be a compact Lie group and let $\rho_d:H_d\to {\rm O}(d)$ be a homomorphism. 
In particular, we can take $(H_d,\rho_d)$ to be a symmetry type in the sense of Freed–Hopkins \cite[Definition~2.4]{FreedHopkins}. Consider the geometric structures, $\mathcal{H}_d^{\nabla}$ and $\mathcal{H}_d^{\rm fl}$, of differential-$(H_d,\rho_d)$-structures and flat $(H_d,\rho_d)$-structures (resp.).
The following spaces are isomorphic in $\Ho \mathscr{S}$:
\begin{enumerate}
\item Smooth deformations of field theories with smooth $(H_d,\rho_d)$-structure:
$$\mathscr{I}_{d}(\mathcal{H}_d):= {\mathcal{F}}{\rm un}^{\otimes}(\Bord_d^{\mathcal{H}_d},\Sigma^dI_{\CC^{\times}})$$
\item Smooth deformations of field theories with differential $(H_d,\rho_d)$-structure:
$$\mathscr{I}_{d}(\mathcal{H}^{\nabla}_d):= {\mathcal{F}}{\rm un}^{\otimes}(\Bord_d^{\mathcal{H}^{\nabla}_d},\Sigma^dI_{\CC^{\times}})$$
\item Smooth deformations of field theories with flat $(H_d,\rho_d)$-structure:
$$ \mathscr{I}_{d}(\mathcal{H}^{\rm fl}_d):= {\mathcal{F}}{\rm un}^{\otimes}(\Bord_d^{\mathcal{H}^{\rm fl}_d},\Sigma^dI_{\CC^{\times}})$$
\item The space of morphisms of spectra:
$$\map(\Sigma^d MTH_d, \Sigma^{d+1}I_{\ZZ(1)}),$$
where $MTH_d$ denotes the Madsen--Tillmann spectrum \cite[Section~3]{GMTW} and $I_{\ZZ(1)}$ is the Anderson dual of the sphere spectum.
\end{enumerate}
\end{theorem}
\proof
We begin by proving that (1) and (4) are isomorphic in $\Ho \mathscr{S}$. Recall the homotopy adjunction $\vert \cdot \vert\dashv c$ in Proposition \ref{formalinv}. Combining Theorem \ref{madsentilleq} and Proposition \ref{shapeanderson}, we have equivalences
\begin{eqnarray*}
\mathcal{F}{\rm un}^{\otimes}(\mathfrak{Bord}_d^{\mathcal{H}_d},\Sigma^dI_{\CC^{\times}}) &=& {\rm Fun}^{\otimes}(\csp(\mathfrak{Bord}_d^{\mathcal{H}_d}),\csp(\Sigma^d I_{\CC^{\times}}))
\\
&\overset{\ref{shapeanderson}}{\simeq} & {\rm Fun}^{\otimes}(\csp(\mathfrak{Bord}_d^{\mathcal{H}_d}),\Sigma^{d+1} I_{\ZZ(1)}) 
\\
&\overset{\vert \cdot\vert\dashv c}{\simeq} & \map (\vert \csp(\mathfrak{Bord}_d^{\mathcal{H}_d})\vert ,\Sigma^{d+1} I_{\ZZ(1)}) 
\\
&\overset{\star}{\simeq} & \map (\csp(\vert \mathfrak{Bord}_d^{\mathcal{H}_d}\vert ) ,\Sigma^{d+1} I_{\ZZ(1)}) 
\\
&=& \map (\csp({\sf MT}(\mathcal{H}_d) ) ,\Sigma^{d+1} I_{\ZZ(1)}) 
\\
&\overset{\ref{madsentilleq}}{\simeq} & \map (\Sigma^d {\rm MT}H_d  ,\Sigma^{d+1} I_{\ZZ(1)}) ,
\end{eqnarray*}
where we have labelled each equivalences by the Theorem, Proposition, or adjunction, used to establish it. The equivalence labelled $\star$ follows since $\vert \cdot \vert$ is the homotopy cocontinuous and $\csp$ is defined as a homotopy colimit (Definition \ref{struts}).

We prove the remaining equivalences by showing that (2) is equivalent to (4) and (3) is equivalent to (4). In both cases, the argument is verbatim the same as the argument for $\mathcal{H}_d$, given by replacing $\mathcal{H}_d$ by $\mathcal{H}_d^{\nabla}$ and $\mathcal{H}_d$ by $\mathcal{H}_d^{\rm fl}$, respectively. The crucial point is that Proposition \ref{madsentilleq} holds for all three of these geometric structures, so that the equivalence 
$$\map(\csp ({\sf MT}(\mathcal{H}_d)) ,\Sigma^{d+1} I_{\ZZ(1)}))\simeq \map(\Sigma^dMTH_d ,\Sigma^{d+1} I_{\ZZ(1)}),$$ 
labelled by \ref{madsentilleq} above, holds for all three structures. 
\endofproof

%
%

\subsection{The Freed--Hopkins conjecture}
In this section, we show how to deduce the Freed--Hopkins conjecture \cite[Conjecture~8.37]{FreedHopkins} from Theorem \ref{hopfred}. Most of the equivariant stable homotopy theory needed to prove the conjecture has already been developed in \cite{FreedHopkins}. In order to keep this work succinct, we choose not to develop these aspects in the setting of smooth spectra and we will mostly defer the reader to \cite{FreedHopkins}. Here we will present only some rather trivial modification which are needed in order to make sense of the spectrum $(I_{\CC^{\times}})_{\rm pos}$ as a sheaf of spectra. 

The main thrust of the argument rests on \cite[Theorem 8.20]{FreedHopkins}, which identifies a certain homotopy pullback of mapping spaces with the space of maps of the form $MTH\to \Sigma^{d+1}I_{\ZZ(1)}$. What we have accomplished in the present work is a comparison of each corner in this homotopy pullback diagram with corresponding deformation spaces of field theories. The starting point for \cite{FreedHopkins} is that certain maps of spectra are in correspondence with certain functorial field theories, which is taken as an Ansatz. Our contribution can be regarded as turning this Ansatz into a theorem.

In Freed–Hopkins \cite[Section 6.3]{FreedHopkins}, real structures on the spectra $I_{\ZZ(1)}$ and $I_{\delta(\CC^\times)}$ are considered. These considerations extend to the smooth spectrum $I_{\CC^{\times}}$ in a straightforward way, which we now explain. A complex conjugation action on $I_{\CC^{\times}}$ is a map
\begin{equation}\label{conjugation}
\nu':{\rm B}\ZZ/2_+\to {\rm B}{\rm h}{\rm Aut}(I_{\CC^{\times}}),
\end{equation}
where ${\rm h}{\rm Aut}(I_{\CC^{\times}})$ is the monoid of self homotopy equivalences of $I_{\CC^{\times}}$. We require that such a map induces complex conjugation on $\pi_1$.

\begin{proposition}\label{inducedand}
A map of the form $\nu'$ in \eqref{conjugation} induces a corresponding map 
$$\nu: {\rm B}\ZZ/2_+\to {\rm B}{\rm h}{\rm Aut}(I_{\ZZ(1)})$$
such that $\nu$ induces the action on $\pi_1$ given by the sign representation on $\ZZ(1)$. 
\end{proposition}
\proof
By Proposition \ref{shapeanderson}, we have an equivalence $b:\Icsp(I_{\CC^{\times}})\overset{\simeq}{\to} \Sigma I_{\delta(\ZZ(1))}$, which is induced by applying $\Icsp$ to the connecting homomorphism $I_{\CC^{\times}}\to \Sigma I_{\delta(\ZZ(1))}$ in the long fiber sequence associated to the exponential sequence. 

Since the smooth deformation space functor ${\sf B}$ sends homotopy equivalences to homotopy equivalences, the functor ${\sf B}$ and the connecting map $b$ induce a map
\begin{equation}\label{bmap}
b_*:{\rm h}{\rm Aut}(I_{\CC^{\times}})\to {\rm h}{\rm Aut}({\sf B}(I_{\CC^{\times}}))\simeq {\rm h}{\rm Aut}(\Sigma I_{\ZZ(1)})\;,
\end{equation}
which is given by applying ${\sf B}$ to a homotopy equivalence and then transferring this equivalence to $\Sigma I_{\ZZ(1)}$ via the equivalence $b$ above.
For a conjugation action $\nu'$ on $I_{\CC^{\times}}$, we define a corresponding action $\nu$ on $I_{\ZZ(1)}$ by $b_*(\nu')=\nu$. 

The identification of the action on $I_{\ZZ(1)}$ can be obtained by observing that the canonical map
$$I_{\delta(\ZZ(1))}\to I_{\CC}=H\CC$$
in the long fiber sequence gives rise to a commutative diagram on sheaves of $\pi_0$:
$$
\xymatrix{
\ZZ(1)\ar[r]\ar[d] & \CC\ar[d]^-{c}
\\
\ZZ(1)\ar[r] & \CC
}
$$
where the map $c$ on the right is complex congugation. The two horizontal arrows are the canonical inclusion $2\pi i \ZZ\into \CC$. The conjugation action $c$ induces restricts to the action by the sign representation $\sigma$ on $2\pi i \ZZ=\ZZ(1)$.

\endofproof

The above observations show that we have a homotopy cofiber sequence of sheaves of $\ZZ/2$-equivariant spectra 
$$I_{\ZZ(1)}^{\nu}\to H\CC^{\nu'}\to I_{\CC^{\times}}^{\nu'}$$
induced by any conjugation action $\nu'$. Moreover, the association $b_*:\nu'\mapsto \nu$ is an equivalence on the space of all maps of the form \eqref{conjugation}, as observed in \cite[Eq. 6.28]{FreedHopkins} \footnote{The argument in \cite{FreedHopkins} is uses the spectrum $I_{\CC^{\times}}$, where $\CC^{\times}$ is equipped with the discrete topology. The argument holds verbative, replacing discrete versions of $\CC$ and $\CC^{\times}$ with their corresponding smooth variants.}.

Now the space of all conjugation actions on $I_{\CC^{\times}}$ has infinitely many connected components. As in \cite{FreedHopkins}, we fix a convenient choice of connected component.

\begin{definition}\label{realstr}
The space of \emph{real structures} on $I_{\ZZ(1)}$ is the path connected component of the space of all $\nu$'s in \eqref{conjugation} containing $\gamma=1-\sigma$, where 
$$1-\sigma:{\rm B}\ZZ/2_+\to {\rm B}{\rm O}\to {\rm B}{\rm h}{\rm Aut}(S^0)\to {\rm B}{\rm h}{\rm Aut}(I_{\ZZ(1)}),$$
 picks out the automorphism induced by the action of $\ZZ/2={\rm O}(1)$ on the virtual representation sphere  $S^{1-\sigma}=(\RR-\RR_{\sigma})^+$, via the $J$-homomorphism. 
\end{definition}

By definition, the space of real structures is connected. Therefore any real structure is (non-canonically) equivalent to the real structure $\gamma=1-\sigma$. The above choice of connected component is motivated by \cite[Example 6.29]{FreedHopkins}, where it is observed that on the category whose objects are finite dimensional vector spaces over $\CC$ equipped with a positive definite Hermitian inner product, and whose morphisms are unitary linear maps, the action given by composing the dual space and complex conjugation is naturally isomorphic to the identity. The $\ZZ/2$ action that corresponds to duality is $\sigma-1$, so that when you compose the action by duality with the conjugation action $1-\sigma$, you indeed get the trivial action.

According to Freed–Hopkins \cite[Ansatz 7.13]{FreedHopkins}, a \emph{continuous} $d$-dimensional extended topological field theory with symmetry group $H_d$ should be a $\ZZ/2$-equivariant map
$$\Sigma^dMTH_d^{\beta}\to \Sigma^{d+1}I_{\ZZ(1)}^{\gamma},$$
where $I_{\ZZ(1)}^{\gamma}=S^{1-\sigma}\wedge I_{\ZZ(1)}$.
Here $MTH_d^{\beta}$ is the (Borel) $\ZZ/2$-equivariant spectrum given by the $\ZZ/2$-action induced by the involution $\beta$ in Lemma \ref{involutions}.

We now discuss the positivity part of ``reflection positivity''. In \cite{FreedHopkins}, an equivariant spectrum $(I_{\ZZ(1)})_H$ is introduced. This spectrum arizes as the homotopy fixed points of the $\ZZ/2$-equivariant spectrum given by the smash product of $I_{\ZZ(1)}^{\gamma}$ and the representation sphere $S^{\sigma-1}$, where $\sigma$ is the real sign representation. That is, we have
$$(I_{\ZZ(1)})_H:=(I_{\ZZ(1)}^{\gamma}\wedge S^{\sigma-1})^{h\ZZ/2},$$
where the superscript $h\ZZ/2$ denotes the homotopy fixed points of the $\ZZ/2$-action. 
Recalling that smashing by $\sigma-1$ corresponds to equipping a spectrum with the action given by taking duals, the spectrum $(I_{\ZZ(1)})_{H}$ is the homotopy fixed points of an equivariant spectrum that has an action corresponding to taking the conjugate dual, which explains the terminology Hermitian.

By \cite[Eq 6.15]{FreedHopkins}, we have an identification 
$$ (I_{\ZZ(1)})_H\simeq I_{\ZZ(1)}^{{\rm B}\ZZ/2_+}.$$
Hence, pulling back along the map ${\rm B}\ZZ/2_+\to \ast$ produces a map $I_{\ZZ(1)}\to (I_{\ZZ(1)})_{H}$, whose image is a summand. Following \cite{FreedHopkins}, we denote this summand as $(I_{\ZZ(1)})_{\rm pos}$.

\begin{definition}\label{sheafyH}
Let $\gamma'$ be the $\ZZ/2$-action on $I_{\CC^{\times}}$ that corresponds to the action $\gamma$ on $I_{\ZZ(1)}$ in Proposition \ref{inducedand}. We define the sheaf of spectra $(I_{\CC^{\times}})_{H}$ as the $\ZZ/2$-equivariant sheaf of spectra given by the homotopy fixed points of the smash product 
 $$ (I_{\CC^{\times}})_{H}:=\left(I_{\CC^{\times}}^{\gamma'}\wedge \delta(S^{\sigma -1})\right)^{h\ZZ/2}.$$
\end{definition}

We define a sheaf of spectra serving as a target for positive field theories as follows. 

\begin{definition}
We define the sheaf of spectra $(I_{\CC^{\times}})_{\rm pos}$ as the homotopy pullback of sheaves of spectra
\begin{equation}\label{pullspec}
\xymatrix{
(I_{\CC^{\times}})_{\rm pos}\ar[r]\ar[d] & (\Sigma I_{\delta(\ZZ(1))})_{\rm pos}\ar[d]
\\
(I_{\CC^{\times}})_{H}\ar[r] & (\Sigma I_{\delta(\ZZ(1))})_{H}.
}
\end{equation}
\end{definition}

The following proposition shows that under passage to deformation spectra, the sheaves of spectra defined above agree with the  corresponding spectra defined via the Anderson dual. 

\begin{proposition}\label{shapeofpos}
We have an equivalences of spectra 
$$\csp((I_{\CC^{\times}})_{\rm pos})\simeq (\Sigma I_{\ZZ(1)})_{\rm pos}, \qquad \csp((I_{\CC^{\times}})_{H})\simeq (\Sigma I_{\ZZ(1)})_{H}, \qquad \csp((I_{\CC^{\times}}^{\gamma'})^{h\ZZ/2})\simeq (\Sigma I_{\ZZ(1)}^{\gamma})^{h\ZZ/2}.$$

\end{proposition}
\proof
Since homotopy pullback squares of sheaves of spectra are also homotopy pushout squares, applying the homotopy cocontinuous functor $\csp$ to the square \eqref{pullspec} gives rise to a homotopy pushout, and hence also homotopy pullback square of spectra:
\begin{equation}\label{pullspecpos}
\xymatrix{
\csp((I_{\CC^{\times}})_{\rm pos})\ar[r]\ar[d] & \csp(\Sigma I_{\delta(\ZZ(1))})_{\rm pos}\ar[d]
\\
\csp((I_{\CC^{\times}})_{H})\ar[r] & \csp((\Sigma I_{\delta(\ZZ(1))})_{H}).
}
\end{equation}
By Proposition \ref{shapeanderson}, the map $\Icsp(I_{\CC^{\times}})\to \Sigma I_{\delta(\ZZ(1))}$ is a weak equivalence. By the relation $\csp\Icsp\simeq  \csp\delta\csp=\csp$ (Proposition \ref{adjointcalc}), applying $\csp$ to this weak equivalence yields a weak equivalence
$$\csp(I_{\CC^{\times}})\to \csp(\Sigma I_{\delta(\ZZ(1))})\simeq \Sigma I_{\ZZ(1)}.$$
Now the correspondance between real structures on $I_{\CC^{\times}}$ and $\Sigma I_{\ZZ(1)}$ is given by restricting the map
$$b_*:{\rm h}{\rm Aut}(I_{\CC^{\times}})\to {\rm h}{\rm Aut}(\csp(I_{\CC^{\times}}))\simeq {\rm h}{\rm Aut}(\Sigma I_{\ZZ(1)})$$
in \eqref{bmap} to the connected component of real structures. Hence for a given real structure on $\nu'$ on $I_{\CC^{\times}}$, the equivalence $b:\Icsp(I_{\CC^{\times}})\to \Sigma I_{\delta(\ZZ(1))}$ is canonically equivariant with respect to the induced real structure on $I_{\delta(\ZZ(1))}$. In particular, for the real structures $\gamma=1-\sigma$ and $\sigma-1$, the map $b$ is an equivariant equivalence. By stability for sheaves of spectra, we can commute the finite homotopy limit indexed by $\deloop \ZZ/2$ with the homotopy colimit defining the functor $\csp$. Hence, passing to homotopy fixed points  implies the third equivalence. It also implies that the bottom horizontal map of \eqref{pullspecpos} is an equivalence, again by passing to homotopy fixed points. Since the diagram is a homotopy pullback, the top horizontal map in the diagram is also an equivalence.

\endofproof

\begin{corollary}\label{reflection}
Fix $d\geq 0$. Let $I_{\CC^{\times}}^{\gamma'}$ be the (Borel) sheaf $\ZZ/2$-equivariant spectra with $\ZZ/2$-action $\gamma'$ corresponding to $\gamma$ in Proposition \ref{inducedand}. The following are isomorphic in $\Ho \mathscr{S}$
\begin{enumerate}\label{freedhopreflstr}
\item[(1)] Deformations of $\ZZ/2$-equivariant field theories with smooth $\mathcal{H}_d$-structure:
$$({\mathcal{F}}{\rm un}^{\otimes})^{\ZZ/2}((\Bord_d^{\mathcal{H}_d})^{\beta},\Sigma^{d+1}I_{\CC^{\times}}^{\gamma'}).$$
\item[(2)] Deformations of $\ZZ/2$-equivariant field theories with differential $\mathcal{H}_d$-structure 
$$({\mathcal{F}}{\rm un}^{\otimes})^{\ZZ/2}((\Bord_d^{\mathcal{H}^{\nabla}_d})^{\beta},\Sigma^{d+1}I_{\CC^{\times}}^{\gamma'}).$$
\item[(3)] Deformations of $\ZZ/2$-equivariant field theories with flat $H_d$-structure:
$$({\mathcal{F}}{\rm un}^{\otimes})^{\ZZ/2}((\Bord_d^{\mathcal{H}^{\rm fl}_d})^{\beta},\Sigma^{d+1}I_{\CC^{\times}}^{\gamma'}).$$
\item[(4)] The space of $\ZZ/2$-equivariant maps 
$$\map^{\ZZ/2}(\Sigma^dMTH^{\beta}_d,\Sigma^{d+1}I^{\gamma}_{\ZZ(1)})$$
\end{enumerate}

\end{corollary}

\proof
The corollary follows from Theorem \ref{hopfred} by passing to homotopy fixed points. The only nontrivial observation is that the correspondence between $\ZZ/2$-actions $\gamma'\mapsto \gamma$ is induced by the connecting map 
$$\csp(I_{\CC^{\times}})\overset{\simeq}{\to} \Sigma I_{\ZZ(1)},$$
as explained in the proof of Proposition \ref{inducedand}.
\endofproof

In the following, we will need to use $k$-truncated field theories, where $0\leq k\leq d$, in the sense of \cite[Section 8.2.1]{FreedHopkins}. The target for these field theories is fixed, depending only on $d$. However, the bordism category is the $k$-dimensional bordism category. 


\begin{lemma}\label{realfixedpts}
Fix nonnegative integers $d$ and $k\leq d$. Let $\gamma'$ be the $\ZZ/2$-action on $I_{\CC^{\times}}$ corresponding to the action $\gamma$ under Proposition \ref{inducedand}. Let $I_{\CC^{\times}}^{\gamma'}$ be the corresponding sheaf of equivariant spectra. The following are isomorphic in $\Ho \mathscr{S}$:
\begin{enumerate}
\item[(1)] The deformation space of $k$-truncated positive field theories with smooth $\mathcal{H}_k$-structure:
$${\mathcal{F}}{\rm un}^{\otimes}(\Bord_k^{\mathcal{H}_k},(\Sigma^{d}I_{\CC^{\times}}^{\gamma'})^{h\ZZ/2}).$$
\item[(2)] The deformation space of $k$-truncated positive field theories with differential $\mathcal{G}^{\nabla}$-structure 
$${\mathcal{F}}{\rm un}^{\otimes}(\Bord_k^{\mathcal{H}_k^{\nabla}},(\Sigma^{d}I_{\CC^{\times}}^{\gamma'})^{h\ZZ/2}).$$
\item[(3)] The deformation space of $k$-truncated positive field theories with flat $\mathcal{G}^{\rm fl}$-structure:
$${\mathcal{F}}{\rm un}^{\otimes}(\Bord_k^{\mathcal{H}_k^{\rm fl}},(\Sigma^{d}I_{\CC^{\times}}^{\gamma'})^{h\ZZ/2}).$$
\item[(4)] The space of maps
$$\map(\Sigma^kMTH_k,(\Sigma^{d+1}I_{\ZZ(1)}^{\gamma})^{h\ZZ/2})$$
\end{enumerate}
\end{lemma}
\proof
The homotopy fiber/cofiber sequence of sheaves of spectra
$$I_{\ZZ(1)}\to I_{\CC}\to I_{\CC^{\times}}$$
enhances to a homotopy fiber/cofibrer sequence of sheaves of equivariant spectra, when $I_{\ZZ(1)}$ equipped with the action $\gamma$, $I_{\CC^{\times}}$ is equipped with the action $\gamma'$, and $I_{\CC}$ is equipped with the action by complex conjugation. Since finite homotopy limits of sheaves of spectra commute with homotopy colimits, we have an induced corresponding homotopy fiber/cofiber sequence on homotopy fixed points. 

Now use Proposition \ref{shapeofpos} and proceed verbatim as in the proof of Theorem \ref{hopfred}, with $(\Sigma^{d}I_{\CC^{\times}}^{\gamma'})^{h\ZZ/2}$ replacing $I_{\CC^{\times}}$ and $\mathfrak{Bord}_k^{\mathcal{H}_k}$ replacing $\mathfrak{Bord}_d^{\mathcal{H}_d}$.
\endofproof

To define the deformation space of reflection positive field theories, we again follow \cite{FreedHopkins}. Recall the extension $\widehat{H}_d$ defined be the split exact sequence \eqref{jmap}. Recall that a splitting $\sigma$ of \eqref{jmap} is given by a choice of hyperplane reflection. Let $s_d$ be the splitting given by reflection accross the hyperplane spanned by all but the first standard basis vector. Let $\phi_{s_d}$ be the corresponding automorphism of $H_d$. Then we have a commutative diagram
\begin{equation}\label{splittingagain}
\xymatrix{
H_{d-1}\ar[d]^{\id}\ar[r]^-{i_d} & H_d\ar[d]^-{\phi_{s_d}}
\\
 H_{d-1}\ar[r]^-{i_d} & H_d
}
\end{equation}
where the inclusion $i_d:H_{d-1}\into H_d$ is given by pullback of the representation $\rho_d:H_d\to {\rm O}(d)$ along the inclusion ${\rm O}(d-1)\into {\rm O}(d)$ given by
$$Q\mapsto \left(\begin{array}{cc}
1 & 0
\\
0 & Q
 \end{array}\right).$$
Let $\beta$ be the involution of the bordism category in Lemma \ref{involutions}, with $\sigma=s_d$. There is a canonical map between bordism categories 
\begin{equation}\label{mapbordsequi}
\mathfrak{Bord}_{d-1}^{\mathcal{H}_{d-1}^{\nabla}}\to (\mathfrak{Bord}_{d}^{\mathcal{H}_{d}^{\nabla}})^{\beta},
\end{equation}
defined as follows. Fix ${\bf m}\in \Delta^{\times d}$, $\langle \ell\rangle\in \Gamma$ and $U\in \cartsp$. The component of the transformation \eqref{mapbordsequi} indexed by this data sends a $(d-1)$-bordism $M$, equipped with a map $f:M\to \mathcal{H}^{\nabla}_{d-1}$, to the $d$-bordism obtained by replacing the ambient manifold $M$ by $M\times \RR$, equipped with the composite map 
$$
M\times \RR\to \mathcal{H}_{d-1}^{\nabla}\xrightarrow{i_d} \mathcal{H}_{d}^{\nabla}.
$$ 
The map $i_d$ is the map on corresponding geometric structures, given by delooping the map \eqref{splittingagain}. On connection forms $i_d$ is given by post-composing with the map on Lie algebras induced by $i_d$. The cut ${\bf m}$-tuple $D$ on $M\times \RR$ is obtained by taking a product of each cut in the ${\bf m}$-tuple $C$ on $M$  by $\RR$. The map to $\langle \ell\rangle$ is given by precomposing the map $P:M\to \langle \ell\rangle$ with the projection $M\times \RR\to M$. As a consequence of the commutativity of the diagram \eqref{splittingagain}, the map \eqref{mapbordsequi} is equivariant when bordism category on the left is equipped with the trivial involution. 


\begin{definition}
Fix $d\geq 0$. We define the deformation space of \emph{reflection positive field theories with differential symmetry type $(H_d,\rho)$} as the homotopy fiber
\begin{equation}\label{reflposthrs}
\xymatrix{
\mathscr{I}_d(\mathcal{H}^{\nabla}_d)_{\rm ref,pos} \ar[r]\ar[d] & \ast \ar[d]
\\
({\mathcal{F}}{\rm un}^{\otimes})^{\ZZ/2}((\Bord_d^{\mathcal{H}^{\nabla}_d})^{\beta},\Sigma^{d}I_{\CC^{\times}}^{\gamma'})\ar[r] & {\mathcal{F}}{\rm un}^{\otimes}(\Bord_{d-1}^{\mathcal{H}^{\nabla}_{d-1}},\Sigma^{d}(I_{\CC^{\times}}^{\gamma'})^{\ZZ/2})\;.
}
\end{equation}
where $\Sigma^{d}(I_{\CC^{\times}}^{\gamma'})^{\ZZ/2}$ denotes the homotopy fixed point spectrum and the bottom map is obtained by pre-composition with \eqref{mapbordsequi}.
\end{definition}

The previous definition does not quite agree with \cite[Definition 8.8]{FreedHopkins}. Instead, we use the alternative definition given by \cite[Equation 8.14]{FreedHopkins}.

\begin{corollary}\label{reflectionpos}
Fix $d\geq 0$. The following spaces are isomorphic in $\Ho \mathscr{S}$:
\begin{enumerate}
\item The space of smooth deformations of reflection positive field theories with smooth $(H_d,\rho_d)$-structure:
$$\mathscr{I}_{d}(\mathcal{H}_d)_{\rm ref,pos}$$
\item The space of smooth deformations of reflection positive field theories with differential $(H_d,\rho_d)$-structure:
$$\mathscr{I}_d(\mathcal{H}^{\nabla}_d)_{\rm ref,pos}$$
\item The space of smooth deformations of reflection positive field theories with flat $(H_d,\rho_d)$-structure:
$$ \mathscr{I}_d(\mathcal{H}^{\rm fl}_d)_{\rm ref,pos}$$
\item The space of maps of spectra:
$$\map(MTH, \Sigma^{d+1}I_{\ZZ(1)}).$$

\end{enumerate}
\end{corollary}
\proof
We prove that (2) and (4) are isomorphic in $\Ho\mathscr{S}$. The claims $(1)\simeq (4)$ and $(3)\simeq (4)$ are  verbatim the same, given by replacing $\mathcal{H}^{\nabla}_d$ with the appropriate geometric structure. 

Use Corollary \ref{reflection} and Lemma \ref{realfixedpts} to identify the two bottom corners in the above diagram with their corresponding counterparts in stable homotopy theory. Then apply \cite[Theorem 8.20]{FreedHopkins}. 

\endofproof

\end{document}